\documentclass{report}

\input{diagrams}
\diagramstyle[PostScript=dvips,nohug]
\usepackage{epsfig}

\newtheorem{thm}{Theorem}[section]
\newtheorem{defn}[thm]{Definition}
\newtheorem{propn}[thm]{Proposition}
\newtheorem{lemma}[thm]{Lemma}
\newtheorem{cor}[thm]{Corollary}

\newcommand{\mcm}[3]{\newcommand{#1}[#2]{{\ensuremath{#3}}}}

\usepackage{latexsym}
\usepackage{amssymb}

\mcm{\blank}{0}{(\emptybk)}
\mcm{\dashbk}{0}{\mbox{---}}
\mcm{\emptybk}{0}{\:\:}
\mcm{\diagspace}{0}{\mbox{\hspace{2em}}}
\newcommand{\done}{\hfill\ensuremath{\Box}}
\newcommand{\bref}[1]{(\ref{#1})}
\newcommand{\ucontents}[2]{\addcontentsline{toc}{#1}{\numberline{}{#2}}}
\mcm{\cat}{1}{\mc{#1}}
\mcm{\fcat}{1}{\mb{#1}}
\mcm{\mc}{1}{\mathcal{#1}}
\mcm{\mr}{1}{\mathrm{#1}}
\mcm{\mb}{1}{\mathbf{#1}}
\mcm{\scat}{1}{\Bbb{#1}}
\mcm{\twid}{1}{\widetilde{#1}}
\mcm{\sub}{0}{\,\subseteq\,}
\mcm{\eqv}{0}{\,\simeq\,}
\mcm{\iso}{0}{\,\cong\,}
\mcm{\of}{0}{\raisebox{0.2mm}{\ensuremath{\scriptstyle\circ}}}
\mcm{\Eee}{0}{\cat{E}}
\newcommand{\epsln}{\varepsilon}
\mcm{\blob}{0}{\scriptscriptstyle{\bullet}}
\mcm{\Hom}{0}{\mr{Hom}}
\mcm{\op}{0}{\mr{op}}
\mcm{\Ab}{0}{\fcat{Ab}}
\mcm{\Alg}{0}{\fcat{Alg}}
\mcm{\HtyAlg}{0}{\fcat{HtyAlg}}
\mcm{\Cat}{0}{\fcat{Cat}}
\mcm{\One}{0}{\fcat{1}}
\mcm{\Set}{0}{\fcat{Set}}
\mcm{\integers}{0}{\mathbb{Z}}
\mcm{\range}{2}{#1,\,\ldots\,,#2}
\mcm{\bftuple}{2}{\tuplebts{\range{#1}{#2}}}
\mcm{\tuple}{3}{\tuplebts{\range{#1,#2}{#3}}}
\mcm{\pr}{2}{\tuplebts{#1,#2}}
\mcm{\triple}{3}{\tuplebts{#1,#2,#3}}
\mcm{\ehom}{3}{#1[#2,#3]}
\mcm{\ftrcat}{2}{[#1,#2]}
\mcm{\homset}{3}{#1(#2,#3)}
\mcm{\multihom}{3}{#1(#2;#3)} 
\mcm{\go}{0}{\rTo}
\mcm{\goby}{1}{\rTo^{#1}}
\mcm{\goesto}{0}{\,\longmapsto\,}
\mcm{\goiso}{0}{\goby{\diso}}
\mcm{\og}{0}{\lTo}
\mcm{\ogby}{1}{\lTo^{#1}}
\mcm{\oppair}{2}{\stackrel{\rTo^{#1}}{\lTo_{#2}}}
\newarrow{Equals}=====
\newarrow{Get}....>
\newarrow{Goesto}|---{->}
\newarrow{Monic}{vee}---{vee}
\newarrow{NT}===={=>}
\mcm{\diso}{0}{\sim}
\mcm{\nat}{0}{\mathbb{N}}	
\mcm{\atuplebts}{1}{\langle #1 \rangle}
\mcm{\tuplebts}{1}{(#1)}
\mcm{\Mod}{0}{\mb{Mod}}
\mcm{\GrMod}{0}{\mb{GrMod}}
\mcm{\GrAb}{0}{\mb{GrAb}}
\mcm{\cee}{0}{\cat{C}}
\mcm{\lwr}{1}{\mb{#1}}
\mcm{\upr}{1}{[#1]}
\newenvironment{proof}{\paragraph*{Proof}}{\paragraph*{}}
\newenvironment{sketchpf}{\paragraph*{Sketch Proof}}{\paragraph*{}}
\mcm{\Top}{0}{\fcat{Top}}
\mcm{\Topstar}{0}{\fcat{Top_*}}
\mcm{\ChCx}{0}{\fcat{ChCx}}
\mcm{\emm}{0}{\cat{M}}
\mcm{\Wtwo}{0}{\raisebox{-.2em}{\epsfig{file=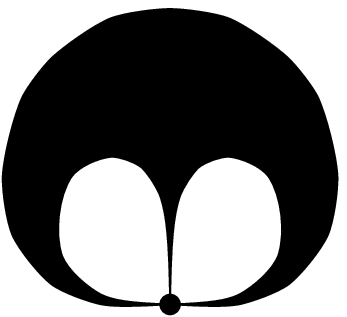,height=1.4em}}}
\mcm{\Wedgepic}{0}{\raisebox{-.2em}{\epsfig{file=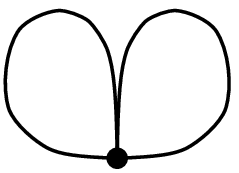,height=.7em}}}
\mcm{\Loops}{0}{\homset{\Topstar}{S^1}{B}}
\mcm{\wej}{0}{\vee}
\mcm{\Obj}{0}{\fcat{Obj}}
\mcm{\SObj}{0}{\fcat{SObj}}
\mcm{\Mon}{0}{\fcat{Mon}}
\mcm{\SMon}{0}{\fcat{SMon}}
\mcm{\Sem}{0}{\fcat{Sem}}
\mcm{\CMon}{0}{\fcat{CMon}}
\mcm{\CSem}{0}{\fcat{CSem}}
\mcm{\Pt}{0}{\fcat{Pt}}
\mcm{\SPt}{0}{\fcat{SPt}}
\mcm{\Sym}{0}{\fcat{Sym}}
\mcm{\Inv}{0}{\fcat{Inv}}
\mcm{\Act}{1}{\fcat{Act}_{#1}}
\mcm{\SAct}{1}{\fcat{SAct}_{#1}}
\mcm{\Lie}{0}{\fcat{Lie}}
\mcm{\GrLie}{0}{\fcat{GrLie}}
\mcm{\Ger}{0}{\fcat{Ger}}
\mcm{\Map}{1}{\fcat{Map}_{#1}}
\mcm{\ehomset}{3}{#1[#2,#3]}
\mcm{\veeh}{0}{\cat{V}}
\mcm{\fmon}{1}{\widehat{#1}}
\newarrow{Incl}C--->
\newcounter{bean}
\mcm{\smsh}{0}{\wedge}
\mcm{\und}{1}{|#1|}
\mcm{\ovln}{1}{\overline{#1}}
\mcm{\ob}{0}{\mr{ob}}	
\newcommand{\url}[1]{\mbox{\tt #1 }}
\mcm{\formal}{3}{#1\atuplebts{#2,\ldots,#3}}

\newlength{\gwidth}	
\newlength{\gvert}	
\newlength{\gdrop}	
\newlength{\gbaredrop}	
\newlength{\goffset}	
\newlength{\gtemp}	
\newcommand{\present}[1]{%
\makebox[1\gwidth]{%
\rule[-1\gdrop]{0ex}{1\gvert}%
\raisebox{-1\gbaredrop}{#1}}}

\newcommand{\cinitdims}[2]{%
\setlength{\unitlength}{1em}%
\setlength{\goffset}{.35\unitlength}%
\setlength{\gwidth}{#1\unitlength}%
\setlength{\gvert}{#2\unitlength}%
\setlength{\gdrop}{.5\gvert}%
\addtolength{\gdrop}{-1\goffset}%
\setlength{\gbaredrop}{1\gdrop}%
\addtolength{\gvert}{.6\unitlength}%
\addtolength{\gdrop}{.3\unitlength}}	

\newcommand{\abovepic}[1]{%
\settoheight{\gtemp}{\ensuremath{#1}}%
\addtolength{\gvert}{1\gtemp}%
\settodepth{\gtemp}{\ensuremath{#1}}%
\addtolength{\gvert}{1\gtemp}}

\newcommand{\belowpic}[1]{%
\settoheight{\gtemp}{\ensuremath{#1}}%
\addtolength{\gvert}{1\gtemp}%
\addtolength{\gdrop}{1\gtemp}%
\settodepth{\gtemp}{\ensuremath{#1}}%
\addtolength{\gvert}{1\gtemp}%
\addtolength{\gdrop}{1\gtemp}}

\newcommand{\cell}[4]{\put(#1,#2){\makebox(0,0)[#3]{\ensuremath{#4}}}}

\newcommand{\prectwo}[3]%
{\begin{picture}(4.2,3.4)(-0.1,-0.2)%
\cell{2}{3.2}{b}{#1}%
\cell{2}{-0.2}{t}{#2}%
\cell{2.2}{1.5}{l}{#3}%
\qbezier(0,2)(2,4)(4,2)%
\qbezier(0,1)(2,-1)(4,1)%
\put(4,2){\vector(1,-1){0}}%
\put(4,1){\vector(1,1){0}}%
\put(2,2.5){\vector(0,-1){2}}%
\end{picture}}

\mcm{\ctwo}{3}{%
\cinitdims{4.2}{3.4}%
\abovepic{#1}%
\belowpic{#2}%
\present{\prectwo{#1}{#2}{#3}}}

\newcommand{\precthree}[5]{%
\begin{picture}(4.2,5.4)(-0.1,-0.2)%
\cell{2}{5.2}{b}{#1}%
\cell{1}{2.7}{b}{#2}%
\cell{2}{-.2}{t}{#3}%
\cell{2.2}{3.75}{l}{#4}%
\cell{2.2}{1.25}{l}{#5}%
\qbezier(0,3)(2,7)(4,3)%
\qbezier(0,2)(2,-2)(4,2)%
\put(0,2.5){\vector(1,0){4}}%
\put(2,4.5){\vector(0,-1){1.5}}%
\put(2,2){\vector(0,-1){1.5}}%
\put(4,3){\vector(1,-3){0}}%
\put(4,2){\vector(1,3){0}}%
\end{picture}}

\mcm{\cthree}{5}{%
\cinitdims{4.2}{5.4}%
\abovepic{#1}%
\belowpic{#3}%
\present{\precthree{#1}{#2}{#3}{#4}{#5}}}

\newcommand{\prectwopar}[4]{%
\begin{picture}(4.2,3.4)(-0.1,-0.2)%
\cell{2}{3.2}{b}{#1}%
\cell{2}{-0.2}{t}{#2}%
\cell{1.6}{1.5}{r}{#3}%
\cell{2.4}{1.5}{l}{#4}%
\qbezier(0,2)(2,4)(4,2)%
\qbezier(0,1)(2,-1)(4,1)%
\put(4,2){\vector(1,-1){0}}%
\put(4,1){\vector(1,1){0}}%
\put(1.8,2.5){\vector(0,-1){2}}%
\put(2.2,2.5){\vector(0,-1){2}}%
\end{picture}}

\mcm{\ctwopar}{4}{%
\cinitdims{4.2}{3.4}%
\abovepic{#1}%
\belowpic{#2}%
\present{\prectwopar{#1}{#2}{#3}{#4}}}

\begin{document}

\title{Homotopy Algebras for Operads}
\author{Tom Leinster\\ \\
        \normalsize{Department of Pure Mathematics, University of
        Cambridge}\\ 
        \normalsize{Email: leinster@dpmms.cam.ac.uk}\\
        \normalsize{Web: http://www.dpmms.cam.ac.uk/$\sim$leinster}}

\date{\normalsize
	\mbox{}\vspace*{5mm}\\
	\textbf{Abstract}\\
	\vspace{3mm} \raggedright \setlength{\rightskip}{0pt}
	We present a definition of homotopy algebra for an operad, and
	explore its consequences.
	\newline\\
	\vspace{-2mm}
	The paper should be accessible to topologists, category theorists,
	and anyone acquainted with operads. After a review of operads and
	monoidal categories, the definition of homotopy algebra is given.
	Specifically, suppose that \emm\ is a monoidal category in which it
	makes sense to talk about algebras for some operad $P$. Then our
	definition says what a homotopy $P$-algebra in \emm\ is, provided
	only that some of the morphisms in \emm\ have been marked out as
	`homotopy equivalences'. 
	\newline\\
	\vspace{-2mm}
	The bulk of the paper consists of examples of homotopy algebras. We
	show 
	that any loop space is a homotopy monoid, and, in fact, that any
	$n$-fold loop space is an $n$-fold homotopy monoid in an appropriate
	sense. We try to compare weakened algebraic structures such as
	$A_\infty$-spaces, $A_\infty$-algebras and non-strict monoidal
	categories to our homotopy algebras, with varying degrees of
	success. We also prove results on `change of base', e.g.\ that the
	classifying space of a homotopy monoidal category is a homotopy
	topological monoid. Finally, we reflect on the advantages and
	disadvantages of our definition, and on how the definition really ought
	to be replaced by a more subtle $\infty$-categorical version.
	\newline\\
	\vspace{-2mm}
	This paper is long (100 pages), but a taste of it can be got from the
	introductory paper~\cite{UTHM} (8 pages), which does not use operads.
	\newline\\}

\maketitle

\addtocontents{toc}{\vspace*{-14mm}}
\tableofcontents

\chapter*{Introduction}
\ucontents{chapter}{Introduction}
\label{ch:intro}

A pressing concern in mathematics is to find a coherent theory of weakened
algebraic structures. In topology this need was apparent quite early on, with
the work of Boardman and Vogt (amongst others) on homotopy-invariant
algebraic structures; however, this only covered `algebraic structures' in
quite a narrow sense, and basically only in the context of topological
spaces. Another aspect is the recent push to develop a workable theory of
weak $n$-categories (perhaps most famously advocated in the Grothendieck
manuscript `Pursuing Stacks'~\cite{Gro}), which has so far resulted in a
multiplicity of proposed definitions whose relationships to one another
remain mysterious. In the last five years or so there has been a flood of
publications involving weakened or up-to-homotopy algebraic structures of
some sort: in algebraic geometry, topology, category theory, quantum algebra,
deformation theory, in the `operad renaissance', and even in mirror
symmetry; there have been far more developments than I can even attempt to
list. 

What follows is a contribution to the theory of weakened algebraic
structures. More specifically, it is a definition, in an appropriate context,
of a \emph{homotopy algebra} for an operad. In other words, given an operad,
one can consider its algebras; the definition given here allows one to
consider also its `weak algebras' or `algebras up to homotopy'. There have
been general definitions of this kind made before, as is more fully discussed
in Chapter~\ref{ch:thoughts}. There are also some very popular notions of
`homotopy algebra' for specific operads: for instance, $A_\infty$-algebras,
strong homotopy Lie algebras, and special $\Gamma$-spaces (to take a random
selection). But I think that the strength of the present definition lies in
its generality. Roughly speaking, suppose that $P$ is an operad and \emm\ a
monoidal category of some sort, so that it makes sense to talk about
$P$-algebras in \emm. Then the only extra ingredient we will need in order to
define homotopy $P$-algebras in \emm\ is the knowledge of which morphisms in
\emm\ are `homotopy equivalences'. So, for example, in order to talk about
homotopy topological algebras we only need to have before us the monoidal
category of spaces together with the knowledge of which continuous maps are
homotopy equivalences; to talk about homotopy differential graded algebras we
only need to know which maps of chain complexes are chain homotopy
equivalences; to talk about homotopy categorical algebras we only need to
know which functors are equivalences of categories. In particular, we do not
need to know what a `homotopy' between maps is, or anything about
resolutions, fibrations, cylinders, etc. Note also that the definition works
in any monoidal category, not just in those (like the category of spaces)
where the monoidal structure is given by cartesian product.

I hope that this paper will be accessible to both category theorists and
topologists, and in fact to anyone acquainted with operads. Although the main
examples come from topology and related areas, the spirit of this work is
fairly conceptual and category-theoretic. With luck, I have included enough
background material that no-one will be put off too rapidly. In particular,
there is a short introduction to operads, and this ought to give a rough idea
of what is going on to those who have not met them before. 

While on the subject of different readerships, I should say a couple of
things about loop spaces, which will mentioned frequently. Given a
topological space $B$ with basepoint, the 
\label{page:defn-loop}%
\emph{loop space} of $B$ is the space of all based loops in $B$; that is, it
is the set of basepoint-preserving maps from the circle $S^1$ into $B$,
endowed with a suitable topology. Thus the connected-components of the loop
space of $B$ are the elements of the fundamental group of $B$. Because loops
can be composed, and composition of loops is associative up to homotopy, any
loop space is a `topological monoid up to homotopy'. What exactly this
means---and there are subtleties concerning higher homotopies which I have
not mentioned---was the subject of a great deal of work by topologists. (A
summary can be found in~\cite{Ad}.) Two of the most popular ways
of saying `topological monoid up to homotopy' precisely are
`$A_\infty$-space' and `special $\Delta$-space', both of which will be
discussed later. For those who know all about loop spaces already, I should
add that the phrase `loop space' will be used to mean a space
\emph{homeomorphic} to the space of loops in some space $B$, not merely
homotopy-equivalent to such a space.

The idea behind this paper is very simple, so I am slightly embarrassed that
it has turned out at 100 pages. Maybe I can reassure the reader that this is
mostly for good reasons. For a start, the pace is slow and the margins
wide. Also, there is a lengthy preliminary chapter on the basics of operads
and monoidal categories, which can be skipped over by many readers. The main
definition (of homotopy algebra for an operad) is then made very quickly;
what takes up the bulk of the paper is that there are lots of examples. These
are the `good reasons' why the paper is long. A `bad reason' is that a
certain amount of unwieldy calculation is present, although this is mostly
sketched rather than done explicitly. As discussed in the final chapter, I
think that while this computational effort is inevitable with today's
technology, it may be possible to give a more conceptual account when the
theory of weak $n$-categories is better developed.

\subsection*{The origin of the idea}

The way I came to the idea behind this paper is as follows. Being a
`historical' explanation of a mathematical idea, it does not represent the
cleanest approach, but perhaps it will be helpful. Readers without a
background in topology, especially, should not let it put them off the rest
of this work.

I had been reading Segal's paper~\cite{SegCCT} defining $\Gamma$-spaces
(nowadays usually called `special $\Gamma$-spaces'), which are a precise
formulation of the idea of up-to-homotopy topological commutative monoids. A
(special) $\Gamma$-space is a contravariant functor from $\Gamma$ to the
category of topological spaces, with certain additional properties which need
not concern us for now; what intrigued me was the process by which the
concept of `commutative monoid' gave rise to the category $\Gamma$.

Explicitly, the objects of $\Gamma$ are all finite sets, and a map from $S$
to $T$ in $\Gamma$ is a function $\theta$ from $S$ to the set of subsets of
$T$, such that $\theta(i) \cap \theta(j) = \emptyset$ when $i \neq j$. The
idea is that the morphisms $\theta: \{\range{1}{m}\}\go \{\range{1}{n}\}$ in
$\Gamma$ are the maps $A^n \go A^m$ which exist for a `generic' commutative
monoid $A$; thus $\theta$ corresponds to the map
\[
\begin{array}{ccc}
A^n			&\go		&A^m	\\
\bftuple{a_1}{a_n}	&\goesto	&\bftuple{b_1}{b_m}
\end{array}
\]
with $b_j = \sum_{i\in\theta (j)} a_i$. But on closer inspection, this idea
is rather strange. It says, for instance, that a typical map $A^8 \go
A^3$ arising `purely because $A$ is a commutative monoid' is something like
\[
(a_1,a_2,a_3,a_4,a_5,a_6,a_7,a_8)
\goesto
(a_3 +a_4 +a_8, a_5, a_2 +a_6).
\]
Note that on the right-hand side, none of the terms $a_i$ are repeated
(because of the restriction $\theta(i) \cap \theta(j) = \emptyset$), but some
of the terms are omitted altogether (namely, $a_1$ and $a_7$). So each $a_i$
can be used 0 times or 1 time, but not 2 or more times. It is hard to see in
what context this would be reasonable. On the one hand, if we are discussing
commutative monoids in the category of sets or of topological spaces, then
the maps $A^n \go A^m$ for a generic commutative monoid $A$ are the $m\times
n$ matrices of natural numbers: an $m\times n$ matrix $X$ corresponds to the
map 
\[
\begin{array}{ccc}
A^n			&\go		&A^m		\\
\bftuple{a_1}{a_n}	&\goesto	&\bftuple{b_1}{b_m}
\end{array}
\]
with $b_j = \sum_{i=1}^{n} X_{ji} a_i$. In this context, each $a_i$ can be
used $k$ times for any $k\geq 0$. On the other hand, suppose we are
discussing commutative monoids in the category of abelian groups, of graded
abelian groups, or of topological abelian groups (i.e.\ commutative rings,
commutative graded rings, or commutative topological rings), so that we no
longer have product-projections $A^{\otimes n} \go A$. Then a map $A^{\otimes
n} \go A^{\otimes m}$ for a generic commutative monoid $A$ is simply a
function $\phi: \{\range{1}{n}\} \go \{\range{1}{m}\}$, corresponding to the
map
\[
\begin{array}{ccc}
A^{\otimes n} 			&\go		&A^{\otimes m}		\\
a_1 \otimes\cdots\otimes a_n	&\goesto	&b_1 \otimes\cdots\otimes b_m
\end{array}
\]
where $b_j = \sum_{i\in \phi^{-1}(j)} a_i$. In this case, each $a_i$ is used
precisely once. The category $\Gamma$ does not fit either situation: it is
neither fish nor fowl.

Categorical logic, which has plenty to say on how an algebraic theory gives
rise to a category, does not provide an immediate answer either. On the one
hand, a commutative monoid in a category \emm\ with finite products is
essentially a finite-product-preserving functor
\[
\fcat{Matr}_{\nat} \go \emm,
\]
where $\fcat{Matr}_{\nat}$ is the category whose objects are $0,1,2,\ldots$
and whose morphisms $n\go m$ are the $m\times n$ matrices of natural
numbers. On the other hand, a commutative monoid in a symmetric monoidal
category \emm\ is essentially a map
\[
\Phi \go \emm
\]
of symmetric monoidal categories, where $\Phi$ is the category of finite sets
and functions, with monoidal structure given by disjoint union (as the
tensor) and the empty set (as the unit). Neither of these categories,
$\fcat{Matr}_{\nat}$ or $\Phi$, is the same as $\Gamma$. 

So, at first it is rather difficult to understand how the category $\Gamma$
arises from the theory of commutative monoids. The answer to the puzzle comes
in the realization that a contravariant functor from $\Gamma$ to spaces can
be described in an alternative but equivalent way: namely, as a `colax
symmetric monoidal functor' from $\Phi$ to spaces. Later, `colax symmetric
monoidal functor' will be defined properly; for now, all that matters is that
it is one of the various possible notions of a map between symmetric monoidal
categories.

We now have two important facts:
\begin{enumerate}
\item 	\label{part:Phi}
the theory of commutative monoids naturally gives rise to the category $\Phi$
\item	\label{part:csmf}
a special $\Gamma$-space (i.e.\ a homotopy topological commutative monoid)
can be defined as a colax symmetric monoidal functor $\Phi\go\mbox{(spaces)}$
with certain additional properties.
\end{enumerate}
Both of these facts are ripe for generalization. In~\bref{part:Phi}, all that
matters about commutative monoids is that they are the algebras for a certain
operad; thus if $P$ is any operad, there is a symmetric monoidal category
\fmon{P} which plays the same role in relation to $P$-algebras as $\Phi$ does
in relation to commutative monoids. In~\bref{part:csmf}, there is nothing
special about the symmetric monoidal category of spaces except that certain
of its morphisms can be distinguished as `homotopy equivalences'; this is
what is used to define the `additional properties' referred to. So the scene
is set: given any operad $P$, and any symmetric monoidal category \emm\ in
which some of the morphisms are called `homotopy equivalences', we ought to
be able to define a \emph{homotopy $P$-algebra in \emm} as a colax symmetric
monoidal functor $\fmon{P} \go \emm$ with certain additional properties. This
is, in fact, what we will do.

\subsection*{Layout}

The paper is laid out as follows. In Chapter~\ref{ch:prelims}, Preliminaries,
we cover the basic facts of operads and monoidal categories, and how the two
concepts are connected. There are no new ideas here. The first `proper'
chapter,~\ref{ch:defn}, gives the definition of homotopy algebra for an
operad. Chapter~\ref{ch:hty-mons} explores homotopy monoids and homotopy
semigroups, including how these relate to special $\Gamma$-spaces and
$\Delta$-spaces, to monoidal categories in the traditional non-strict sense,
to $A_\infty$-spaces, and to $A_\infty$-algebras; this chapter also contains
a proof that any loop space is a homotopy topological monoid. In
Chapter~\ref{ch:other} we look at some other examples of homotopy algebras,
including the homotopy algebraic structure on an iterated loop space, and
make a closer examination of homotopy algebras in the category \Cat\ of
categories. Chapter~\ref{ch:change} is on `change of environment' or `change
of base', and includes such results as `the classifying space of a monoidal
category is a homotopy topological monoid', as well as a way of explaining
`why' the higher homotopy groups of a space are abelian. We finish with a
discussion (Chapter~\ref{ch:thoughts}) of various general issues arising
in the course of the paper, including homotopy invariance, the relation of
this work to higher-dimensional category theory, and the pros and cons of our
definition of homotopy algebra. There is also a glossary, listing most of the
notation we have introduced.

\subsection*{Acknowledgements}

Many people have helped me in the course of this research. I am especially
grateful to Martin Hyland, Jim Stasheff and Graeme Segal, each of whom has
given generously of his time and expertise. Elizabeth Mann, Ian Grojnowski,
Craig Snydal, Martin Markl, Paolo Salvatore, Ivan Le Creurer, Todd Trimble,
Nils Baas, Bj\o rn Ian Dundas, Rainer Vogt, Tim Porter and Ross Street have
also assisted me in assorted ways, and I thank them for it. I would also like
to thank the organizers of the very stimulating workshop on operads in
Utrecht, June 1999. At Utrecht, Morten Brun told me a little about PROPs, for
which I am grateful to him; he also told me that he had been thinking about
$\Gamma$-spaces in terms of colax monoidal functors, which is the approach of
this paper too. I am informed by Rainer Vogt that this way of thinking is
well-known amongst topologists, but the only reference I have is a
preprint~\cite{FF} of Floyd and Floyd, which I have not seen.

This work was supported by the Laurence Goddard Fellowship at St John's
College, Cambridge, and by an EPSRC grant. The diagrams without curves were
made with Paul Taylor's commutative diagrams package.

\chapter{Preliminaries}
\label{ch:prelims}

This paper makes fundamental use of two basic concepts: operads and monoidal
categories. In this preliminary chapter we review operads and monoidal
categories, and look at the connection between the two. To do this it is
useful, although not essential, to look also at multicategories (`coloured
operads'). We will also need to be able to do everything in the enriched
context---that is, in the context where operations form structures more
complex than mere sets---and so we review enrichment for both categories
and operads.

Everything in this chapter is old hat, and as such it might not be obvious why
it needed to be written. My reason is that although the topics covered
(operads, multicategories, monoidal categories and enrichment) form a natural
unit, not everyone who knows some of it knows all of it. Specifically, I hope
that this paper will be read both by some who would describe themselves as
topologists and some who self-define as category theorists; but in category
theory it is (sadly) not widely appreciated that an operad is a very natural
categorical structure, while topologists are perhaps not so conversant with
multicategories and enrichment. So I have included sketches of all these
ideas.

That said, experts will be able to skip lightly over this chapter, pausing
perhaps to pick up some notation and terminology. (Note, in particular, the
terminology concerning symmetric \emph{vs.}\ non-symmetric
operads~(\ref{sec:operads}), the context in which algebras are
taken~(\ref{sec:operads}), and the notation \fmon{P} for the free monoidal
category on an operad $P$~(\ref{sec:free-mon}).) There is a glossary at the
end of the paper containing the names of commonly-used categories and
operads. 

We start with a review of monoidal categories~(\ref{sec:prelims-mon-cats}),
operads~(\ref{sec:operads}), and
multicategories~(\ref{sec:multicats}). Enriched categories and operads are
covered in~\ref{sec:enr-cats}
and~\ref{sec:enr-ops}. Finally~(\ref{sec:free-mon}) we look at how to form
the free monoidal category on an operad.

\paragraph*{Miscellaneous notation}	\label{page:misc-notation}
If \cee\ is a category and $A$, $B$ are objects of \cee\ then
\homset{\cee}{A}{B} means \homset{\mr{Hom}_\cee}{A}{B}. The opposite (dual)
category of \cee\ is written $\cee^\op$; thus $\homset{\cee^\op}{B}{A} =
\homset{\cee}{A}{B}$. If
\[
\cat{C}\ctwo{F}{G}{\alpha}%
\cat{C'}\ctwo{F'}{G'}{\alpha'}%
\cat{C''}
\]
is a diagram of categories, functors and natural transformations, then we
write $\alpha' * \alpha$ for the
composite natural transformation
\[
\cat{C}\ctwo{F'\of F}{G'\of G}{}\cat{C''}.
\]

The symbol $\iso$ denotes isomorphism, whereas $\eqv$ means equivalence
(between categories, spaces, etc.). 

$0$ is an element of the natural numbers, \nat.

\section{Monoidal Categories}	\label{sec:prelims-mon-cats}

We will consider monoidal categories \emm, not necessarily symmetric, in
which the monoidal product is written $\otimes$ and the unit object is
written $I$. The associativity and unit isomorphisms will go nameless, as
will the symmetries $A\otimes B \go B\otimes A$ in symmetric monoidal
categories.

A strict monoidal category is one in which the associativity and unit
isomorphisms are actually identities. The coherence theorem for monoidal
categories states that any monoidal category is equivalent (in a suitable
sense) to a strict monoidal category (see \cite{JS}). This justifies leaving
out the brackets in expressions such as $A\otimes B\otimes C$, which we will
often do. In the symmetric case, the corresponding coherence result is that
every symmetric monoidal category is equivalent to a symmetric strict
monoidal category. Note that the word `strict' qualifies `monoidal' but not
`symmetric' in the term `symmetric strict monoidal category': we can force
the tensor product to satisfy strict associativity and unit laws, but not to
be strictly commutative.

\subsubsection*{Examples}

\begin{enumerate}
\item	\label{eg:mon-cat:cart}
Let \emm\ be any category in which finite products exist. By choosing a
particular product $A\times B$ for each pair \pr{A}{B} of objects, and a
particular terminal object $1$, one obtains a symmetric monoidal category
\triple{\emm}{\times}{1} in a natural way. A monoidal category arising in
this way is called a \emph{cartesian monoidal category}. 

\item	\label{eg:mon-cat:Cat}
\Cat\ is the category of (small) categories and functors. The usual
(cartesian) product $\times$ and the terminal category \One\ make \Cat\ into
a symmetric monoidal category.

\item	\label{eg:mon-cat:Set}
\Set, the category of sets and functions, with cartesian product $\times$ and
the one-element set $1$, forms a symmetric monoidal category.

\item	\label{eg:mon-cat:Mod}
$\Mod = \Mod_R$ is the category of left modules over a fixed commutative ring
$R$. It has a symmetric monoidal structure given by $\otimes$ and $R$. 

\item	\label{eg:mon-cat:GrMod}
$\GrMod = \GrMod_R$ is the category of $\integers$-graded $R$-modules with
the usual $\otimes$,
\[
(A\otimes B)_n = \bigoplus_{p+q=n} A_p \otimes B_q,
\]
and unit object $R$ (abusing notation slightly) given by 
\[
R_n = \left\{
\begin{array}{ll}
R	&\textrm{if }n=0		\\
0	&\textrm{if }n\neq 0.
\end{array}
\right.
\]
There are (at least) two possible symmetries on \GrMod: we can define
$\gamma: A\otimes B \go B\otimes A$ either by $\gamma(a\otimes b) = b\otimes
a$ or by $\gamma(a\otimes b) = (-1)^{pq} b\otimes a$, for $a\in A_p$ and
$b\in B_q$. We shall generally use the latter.

\item	\label{eg:mon-cat:ChCx}
$\ChCx = \ChCx_R$ is the category of \integers-graded chain complexes of
$R$-modules, that is, of diagrams
\[
\cdots \go A_1 \goby{d} A_0 \goby{d} A_{-1} \go \cdots
\]
of $R$-modules with $d\of d = 0$. The tensor and unit object (also denoted
$R$) are the usual ones, and the symmetry is $\gamma(a\otimes b) = (-1)^{pq}
b\otimes a$ for $a\in A_p, b\in B_q$. (This time there is no choice; the
other $\gamma$ mentioned for \GrMod\ isn't a chain map.)

\item	\label{eg:mon-cat:Top}
\Top\ is the category of topological spaces and continuous maps, with
symmetric monoidal structure given by cartesian product. At various points we
will actually need to use some cartesian closed version of the category of
spaces, i.e.\ a version where it is possible to form function spaces. In
these situations \Top\ will mean the category of compactly generated
Hausdorff spaces and continuous maps. This category carries a symmetric
monoidal structure by virtue of having products, but these products are not
the same as those in the category of all topological spaces. In what follows,
the issue of which version of \Top\ is appropriate is usually swept under the
carpet. 

\item	\label{eg:mon-cat:Topstar}
\Topstar\ is the category of based spaces, whose objects are topological
spaces with basepoint and whose maps are continuous basepoint-preserving
functions. As with \Top\ above, sometimes we will really mean `compactly
generated Hausdorff space' rather than `space'. Two (symmetric) monoidal
structures on \Topstar\ will be of interest to us. Firstly, there is that
given by the product $\times$ in \Topstar\ and the one-point space
$1$. Secondly, there is the wedge product \wej\ (join two spaces by their
basepoints), whose unit is also~$1$.

\item	\label{eg:mon-cat:Phi}	
$\Phi$ is the skeletal category of finite sets: its objects are the finite
sets $n=\{0,\ldots,n-1\}$ for each integer $n\geq 0$, and its maps are all
functions. This is equivalent to the category of all finite sets. Addition
(disjoint union) provides a monoidal product, with unit object the empty set
$0$. Then $\Phi$ becomes a symmetric monoidal category.

\item	\label{eg:mon-cat:Delta}	
So far all the examples of monoidal categories have had a symmetry on
them. This one does not. Let $\Delta$ be the category whose objects are the
finite sets $n$ as in \bref{eg:mon-cat:Phi} for $n\geq 0$, and whose maps are
the order-preserving functions with respect to the obvious total order on
each $n$. Thus $\Delta$ is equivalent to the category of all finite totally
ordered sets. Note that the empty set $0$ is an object of $\Delta$, so that
$\Delta$ is one object bigger than the category usually denoted $\Delta$ by
topologists. Our $\Delta$ becomes a monoidal category via $+$ and $0$, but
there is no symmetry. (The identity maps $m+n \go n+m$ don't provide one,
because they don't form a natural transformation.)

\end{enumerate}

There are various notions of a map between monoidal categories, and the
distinction is important in this work. I will therefore give a formal
definition of monoidal functor, and of monoidal transformation too. 

\begin{defn}	\label{defn:mon-ftr}
Let \cat{L} and \emm\ be monoidal categories. A \emph{monoidal functor}
$\pr{X}{\xi}: \cat{L} \go \emm$ consists of a functor $X: \cat{L} \go \emm$
together with isomorphisms
\[
\begin{array}{rccc}
\xi_{A,B}:	&X(A\otimes B)	&\go	&X(A) \otimes X(B)	\\
\xi_0:		&X(I)		&\go	&I,			
\end{array}
\]
the former natural in $A,B \in \cat{L}$, such that the following diagrams
commute for all $A, B, C \in \cat{L}$:
\[
\begin{diagram}
X(A\otimes B\otimes C)		&\rTo^{\xi_{A\otimes B, C}}
&X(A\otimes B) \otimes X(C)	\\
\dTo<{\xi_{A, B\otimes C}}	&
&\dTo>{\xi_{A,B} \otimes 1}	\\
X(A) \otimes X(B\otimes C)	&\rTo_{1\otimes \xi_{B,C}}
&X(A) \otimes X(B) \otimes X(C)	\\
\end{diagram}
\]\[
\begin{diagram}
X(I\otimes A)	&\rTo^{\xi_{I,A}}	&X(I) \otimes X(A)	\\
		&\rdTo<{\diso}		&\dTo>{\xi_0 \otimes 1}	\\
		&			&I \otimes X(A)		\\
\end{diagram}
\diagspace
\begin{diagram}
X(A) \otimes X(I)	&\lTo^{\xi_{A,I}}	&X(A\otimes I)	\\
\dTo<{1\otimes \xi_0}	&\ldTo>{\diso}		&		\\
X(A) \otimes I.		&			&		\\
\end{diagram}
\]
If \cat{L} and \emm\ are symmetric monoidal categories then a \emph{symmetric
monoidal functor} $\cat{L} \go \emm$ consists of an \pr{X}{\xi} as above,
satisfying the additional axiom that
\begin{diagram}
X(A\otimes B)	&\rTo^{\xi_{A,B}}	&X(A) \otimes X(B)	\\
\dTo<{\diso}	&			&\dTo>{\diso}		\\
X(B\otimes A)	&\rTo^{\xi_{B,A}}	&X(B) \otimes X(A)	\\
\end{diagram}
commutes for each $A,B\in\cat{L}$.
\end{defn}
Note that the maps $\xi_{A,B}$ and $\xi_0$ are required to be isomorphisms;
thus tensor and unit are preserved up to coherent isomorphism.

\begin{defn}	\label{defn:mon-transf}
Let \cat{L} and \emm\ be (symmetric or not) monoidal categories, and let
\[
\pr{W}{\omega}, \pr{X}{\xi}: \cat{L} \go \emm
\]
be (symmetric or not) monoidal functors. A \emph{monoidal transformation}
\[
\sigma: \pr{W}{\omega} \go \pr{X}{\xi}
\] 
is a natural transformation $\sigma: W\go X$ such that the following
coherence diagrams commute ($A,B\in \cat{L}$):
\[
\begin{diagram}
W(A\otimes B)		&\rTo^{\sigma_{A\otimes B}}	&X(A\otimes B)	\\
\dTo<{\omega_{A,B}}	&				&\dTo>{\xi_{A,B}}\\
W(A)\otimes W(B)	&\rTo_{\sigma_A \otimes \sigma_B}&X(A)\otimes X(B)\\
\end{diagram}
\diagspace
\begin{diagram}
W(I)		&\rTo^{\sigma_I}	&X(I)		\\
\dTo<{\omega_0}	&			&\dTo>{\xi_0}	\\
I		&\rEquals		&I.		\\
\end{diagram}
\]
\end{defn}

Thus if \cat{L} and \emm\ are monoidal categories, there is a category
\homset{\Mon}{\cat{L}}{\emm} of monoidal functors from \cat{L} to \emm\ and
monoidal transformations. Similarly, if \cat{L} and \emm\ are symmetric
monoidal categories then there's a category \homset{\SMon}{\cat{L}}{\emm} of
symmetric monoidal functors and monoidal transformations.

\section{Operads}
\label{sec:operads}

In this section are the definitions of operad and of algebra for an operad, with
numerous examples. The reader is probably very familiar with the definitions;
nevertheless, I have included them to establish my terminology for operads
and my context for algebras, both of which are slightly non-standard. 

\begin{defn}	\label{defn:operad}
\begin{enumerate}
\item
A \emph{non-symmetric operad} $P$ consists of a sequence
$(P(n))_{n\in\nat}$ of sets, together with an element $1\in P(1)$ and a
function 
\[
\begin{array}{ccc}
P(n) \times P(k_1) \times \cdots \times P(k_n) 
&\go 		&P(k_1 + \cdots + k_n)			\\
(\theta,\theta_1,\ldots,\theta_n)
&\goesto	&\theta\of(\theta_1,\ldots,\theta_n)
\end{array}
\]
for each $n, k_1, \ldots, k_n \geq 0$, satisfying unit and associativity
axioms.

\item
A \emph{symmetric operad} consists of a non-symmetric operad $P$ together
with a right action of the symmetric group $S_n$ on $P(n)$, for each $n$,
satisfying compatibility laws.
\end{enumerate}
\end{defn}

\subsubsection*{Points to note}

\begin{itemize}
\item
The exact axioms can be found in \cite{MayDOA}.

\item
By default, our operads will be operads of sets: that is, each $P(n)$ is a
set rather than a space or an abelian group, etc. Later (\ref{sec:enr-ops})
we will consider these more sophisticated kinds of $P(n)$.

\item
We give equal emphasis to symmetric operads (usually just called `operads' in
the literature) and non-symmetric operads (also called `non-$\Sigma$
operads'). We will generally regard the non-symmetric case as the more basic,
and the symmetric case as an elaboration of it. The term \emph{operad} on its
own will refer to both cases equally.

\item
Operads always have a unit element $1\in P(1)$. There is no requirement
(unlike in \cite[1.1]{MayGIL}) that $P(0)$ has only one element.
\end{itemize}

\begin{defn}	\label{defn:algebra}
\begin{enumerate}
\item	\label{defn:non-sym-alg}
Let $P$ be a non-symmetric operad and let \emm\ be a monoidal category. An
\emph{algebra for $P$ in \emm} (or a \emph{$P$-algebra in \emm}) consists of
an object $A$ of \emm\ together with a function
\[
P(n) \go \homset{\emm}{A^{\otimes n}}{A}
\]
for each $n$, written $\theta\goesto\ovln{\theta}$ and satisfying some
axioms. 
\item
Let $P$ be a symmetric operad and let \emm\ be a symmetric monoidal
category. An \emph{algebra for $P$ in \emm} consists of an algebra for $P$
in \emm\ in the sense of~\bref{defn:non-sym-alg}, satisfying further axioms
concerning symmetries. 
\end{enumerate}
\end{defn}

The axioms can be found in \cite{MayDOA}. There is an obvious notion of a map
of algebras, and we thus obtain the category \homset{\Alg}{P}{\emm} of
$P$-algebras in \emm. (We will always regard an operad $P$ as either being
symmetric or being non-symmetric; however, the notation \homset{\Alg}{P}{\emm}
does not reveal which.)

\subsubsection*{Examples}

\begin{enumerate}

\item	\label{eg:operad:Obj}
Let \Obj\ be the unique non-symmetric operad with
\[
\Obj(n) = \left\{
\begin{array}{ll}
1		&\textrm{if }n=1	\\
\emptyset	&\textrm{otherwise}.
\end{array}
\right.
\]
Then an algebra for \Obj\ in a monoidal category \emm\ is merely an object of
\emm, and in fact $\homset{\Alg}{\Obj}{\emm} \iso \emm$. Similarly, let
\SObj\ be the \emph{symmetric} operad defined by the same formula as \Obj;
then $\homset{\Alg}{\SObj}{\emm} \iso \emm$ for any symmetric monoidal
category \emm.

\item	\label{eg:operad:Mon}	
Let \Mon\ be the unique non-symmetric operad with $\Mon(n)=1$ for all $n\geq
0$. Then an algebra for \Mon\ in a monoidal category \emm\ is simply a
monoid in \emm: that is, an object $A$ of \emm\ equipped with maps
$m:A\otimes A \go A$ and $e:I\go A$ such that $m$ is associative and $e$ is a
unit for $m$. So:
\begin{itemize}
\item	
\homset{\Alg}{\Mon}{\Set} is the category of monoids in the usual sense
\item	
\homset{\Alg}{\Mon}{\Cat} is the category of \emph{strict} monoidal
categories 
\item	
\homset{\Alg}{\Mon}{\Mod_R} is the category of algebras over the commutative
ring $R$; when we speak of algebras (or graded algebras, etc.) over a ring,
we always mean \emph{unital} algebras
\item	
\homset{\Alg}{\Mon}{\GrMod_R} is the category of graded $R$-algebras
\item	
\homset{\Alg}{\Mon}{\ChCx_R} is the category of differential graded
$R$-algebras 
\item
\homset{\Alg}{\Mon}{\Top} is the category of topological monoids
\item
\homset{\Alg}{\Mon}{\Topstar} is also (isomorphic to) the category of
topological monoids.
\end{itemize}

\item	\label{eg:operad:Sem}	
Let \Sem\ be the unique non-symmetric operad with 
\[
\Sem(n) = \left\{
\begin{array}{ll}
1		&\textrm{if }n\geq 1	\\
\emptyset	&\textrm{if }n=0.
\end{array}
\right.
\]
Then an algebra for \Sem\ in a monoidal category \emm\ is a semigroup in
\emm, that is, an object $A$ of \emm\ equipped with an associative binary
operation $A\otimes A \go A$. So
\begin{itemize}
\item
\homset{\Alg}{\Sem}{\Set} is the category of semigroups in the usual sense
\item
\homset{\Alg}{\Sem}{\Mod_R} is the category of non-unital $R$-algebras, and
similarly \GrMod\ and \ChCx\ for graded and differential graded non-unital
algebras 
\item
\homset{\Alg}{\Sem}{\Top} is the category of topological semigroups 
\item
An object of \homset{\Alg}{\Sem}{\Topstar} is a topological semigroup with a
distinguished idempotent element; note that this is not necessarily a
topological monoid.
\end{itemize}

\item	\label{eg:operad:CMon}
Let \CMon\ be the unique symmetric operad with $\CMon(n)=1$ for all $n$ (this
being the symmetric analogue of the non-symmetric operad \Mon\ in
\bref{eg:operad:Mon}). Then a \CMon-algebra in a symmetric monoidal category
\emm\ is a commutative monoid in \emm; all but one of the examples of \emm\
given in \bref{eg:operad:Mon} can be repeated here, with the word
`commutative' inserted each time.  The exception is \Cat: a commutative
monoid in \Cat\ is a strict monoidal category in which $x\otimes y = y\otimes
x$ and $f\otimes g = g\otimes f$ for all objects $x,y$ and morphisms
$f,g$. Such `strictly symmetric strict monoidal categories' are very rare in
nature; see the comments on coherence at the start of
Section~\ref{sec:prelims-mon-cats}. 

\item	\label{eg:operad:CSem}
Let \CSem\ be defined by the same formula as \Sem\ in \bref{eg:operad:Sem},
but now regarded as a symmetric operad (in the only possible way). Then a
\CSem-algebra in a symmetric monoidal category \emm\ is a commutative
semigroup in \emm.

\item	\label{eg:operad:Pt}
Let \Pt\ be the unique non-symmetric operad defined by 
\[
\Pt(n) = \left\{
\begin{array}{ll}
1		&\textrm{if }n=0\textrm{ or }n=1	\\
\emptyset	&\textrm{if }n\geq 2.
\end{array}
\right.
\]
Then a \Pt-algebra in a monoidal category \emm\ is a pointed
object of \emm, i.e.\ an object $A$ of \emm\ together with a map $I\go
A$. So:
\begin{itemize}
\item
\homset{\Alg}{\Pt}{\Set} is the category of pointed sets (i.e.\ sets with a
distinguished element)
\item
$\homset{\Alg}{\Pt}{\Top} \iso \homset{\Alg}{\Pt}{\Topstar} \iso \Topstar$
\item
\homset{\Alg}{\Pt}{\Mod_R} is the category in which an object is an
$R$-module with a chosen element, and a morphism is a homomorphism preserving
chosen elements
\item
Similarly, an object of \homset{\Alg}{\Pt}{\GrMod} is a graded module $A$
together with a chosen element of $A_0$
\item
An object of \homset{\Alg}{\Pt}{\ChCx} is a chain complex $A$ together with a
chosen cycle in $A_0$
\item
An object of \homset{\Alg}{\Pt}{\Cat} is a (small) category with a chosen
object.
\end{itemize}
A symmetric operad \SPt\ can be defined by the same formula as \Pt, and 
\[
\homset{\Alg}{\SPt}{\emm} \iso \homset{\Alg}{\Pt}{\emm}
\]
for any symmetric monoidal category \emm.

\item	\label{eg:operad:Sym}
There is a symmetric operad \Sym\ in which $\Sym(n)$ is $S_n$, the $n$th
symmetric group. The unit for this operad is uniquely determined and the
symmetric group actions are by translation; a description of the composition
is more lengthy, but can be found, effectively, in
\cite[1.1.1(c)]{MayDOA}. If \emm\ is any symmetric monoidal category then
\homset{\Alg}{\Sym}{\emm} is the category of monoids in \emm, which category
we have also described via the non-symmetric operad \Mon\ in
\bref{eg:operad:Mon}.

\item	\label{eg:operad:Inv}
Let us say that a \emph{monoid with involution} in a symmetric monoidal
category \emm\ is a monoid $A$ in \emm\ together with a map $\blank^*: A\go
A$, satisfying commuting diagrams corresponding to the equations
\[
1^* = 1, \diagspace
(a\cdot b)^* = b^* \cdot a^*.
\]
(For instance, any group becomes a monoid with involution in \Set, by
defining $a^* = a^{-1}$.) There is an operad \Inv\ such that \Inv-algebras
are monoids with involution, in any symmetric monoidal category \emm. More
explicitly,
\[
\Inv(n) = C_2^n \times S_n
\]
where $C_2$ is the cyclic group of order 2; the description of the rest of
the operad structure is omitted.

\item	\label{eg:operad:Act-G}
Fix a monoid $G$ (in \Set). Then there is a non-symmetric operad \Act{G},
with
\[
\Act{G}(n) = \left\{
\begin{array}{ll}
G		&\textrm{if }n=1	\\
\emptyset	&\textrm{otherwise.}
\end{array}
\right.
\]
Composition and identity in \Act{G} are given by multiplication and
identity in $G$. An algebra for \Act{G} in a monoidal category
\emm\ is a \emph{left $G$-object} in \emm: that is, an object $A$ of \emm\
equipped with a map $g\cdot\dashbk: A\go A$ for each $g\in G$,
satisfying the usual axioms for an action. E.g.:
\begin{itemize}
\item
\homset{\Alg}{\Act{G}}{\Set} is the category of left $G$-sets
\item
\homset{\Alg}{\Act{G}}{\Mod_R} is the category of $R$-linear representations
of $G$.
\end{itemize}
There is also a \emph{symmetric} operad \SAct{G}, defined by the same formula
as \Act{G}. If \emm\ is a symmetric monoidal category (as in the two examples
just mentioned) then 
\[
\homset{\Alg}{\SAct{G}}{\emm} 
\iso
\homset{\Alg}{\Act{G}}{\emm}.
\]
\end{enumerate}

\section{Multicategories}
\label{sec:multicats}

We shall occasionally make passing reference to
multicategories. Multicategories are the same as `coloured operads' or `typed
operads', if we ignore the issue of whether or not there are symmetric group
actions, and they are a very natural generalization of operads.  As far as I
can tell, multicategories were actually invented a little earlier than
operads, their first applications being in logic, linguistics and computer
science rather than in topology. (See Lambek's original paper \cite{Lam}.)

So, a \emph{multicategory} $P$ (or `non-symmetric multicategory', for
emphasis) consists of 
\begin{itemize}
\item a collection $\ob(P)$ of objects
\item for each $n\geq 0$ and $a_1,\ldots, a_n, a \in \ob(P)$, a set
$
\multihom{P}{\range{a_1}{a_n}}{a},
$
whose elements are written $\theta:\range{a_1}{a_n}\go a$
\item for each $a\in \ob(P)$, an
`identity' element $1_a$ of \multihom{P}{a}{a} 
\item `composition' functions
\begin{eqnarray*}
\multihom{P}{\range{a_1}{a_n}}{a} \times
\multihom{P}{\range{a_1^1}{a_1^{k_1}}}{a_1}
\times \cdots \times
\multihom{P}{\range{a_n^1}{a_n^{k_n}}}{a_n}\\
\ \go
\multihom{P}{\range{a_1^1}{a_n^{k_n}}}{a}.
\end{eqnarray*}
\end{itemize}
Associative and unit laws must be obeyed. The exact definition can be found
in \cite[p.~103]{Lam}.

A \emph{symmetric multicategory} is a multicategory $P$ together with a
function
\[
\begin{array}{ccc}
\multihom{P}{\range{a_1}{a_n}}{a}
&\go&
\multihom{P}{\range{a_{\sigma(1)}}{a_{\sigma(n)}}}{a},	\\
\theta &\goesto& \theta.\sigma
\end{array}
\]
for each $a_1,\ldots, a_n, a \in \ob(P)$ and $\sigma\in S_n$, satisfying
further axioms.

A one-object multicategory is precisely an operad (in both the
non-symmetric and symmetric flavours): if the single object of the
multicategory $P$ is called $a$, then the corresponding operad $P'$ has
\[
P'(n) = \multihom{P}{\underbrace{\range{a}{a}}_{n}}{a}.
\]
We shall subsequently write the operad and the multicategory both as $P$. 

This provides one source of examples of multicategories. Another comes from
monoidal categories: if \emm\ is a monoidal category then there is a
multicategory $\ovln{\emm}$ with the same objects as \emm, and with
\[
\multihom{\ovln{\emm}}{\range{a_1}{a_n}}{a}
=
\homset{\emm}{a_1\otimes\cdots\otimes a_n}{a}.
\]
Alternatively, we could obtain a multicategory by taking only some of the
objects of \emm; if we take only one then we obtain the familiar
endomorphism operads. If \emm\ is a \emph{symmetric} monoidal category, then
$\ovln{\emm}$ becomes a symmetric multicategory.

We can now rephrase the definition of algebra. An algebra in a monoidal
category \emm\ for a non-symmetric operad $P$ is simply a map $P\go
\ovln{\emm}$ of multicategories; similarly for the symmetric
version. More generally, we can define an algebra in \emm\ for a
\emph{multicategory} $P$ as a multicategory map $P\go \ovln{\emm}$.

One more example of a multicategory will be referred to later
(Chapter~\ref{ch:thoughts}). Fix an operad $P$---say non-symmetric, for
simplicity. Then there is a 2-object multicategory (a `2-coloured operad')
\Map{P}, with the property that a \Map{P}-algebra in a monoidal category
\emm\ consists of a pair \pr{A_0}{A_1} of $P$-algebras in \emm\ together with
a $P$-algebra map $A_0 \go A_1$. Explicitly, define the objects of \Map{P} to
be $\{0,1\}$, and define the `hom-sets' of \Map{P} by
\begin{eqnarray*}
\multihom{\Map{P}}{\range{a_1}{a_n}}{0}
&=&
\left\{
\begin{array}{ll}
P(n)		&\textrm{if }a_1=\cdots=a_n=0	\\
\emptyset	&\textrm{otherwise,}
\end{array}
\right.	\\
\multihom{\Map{P}}{\range{a_1}{a_n}}{1}
&=&
P(n) \textrm{ for all } \range{a_1}{a_n}\in \{0,1\}.
\end{eqnarray*}
Composition and identities are defined as in $P$.

\section{Enriched Categories}	
\label{sec:enr-cats}

The core ideas of this paper do not depend at all on the idea of
\emph{enrichment}, in which `hom-sets' (and similar things) are not in fact
sets but some richer kind of structure. However, the core ideas (such as the
definition of homotopy algebra) can all be extended to the enriched setting;
we will repeatedly say `here's an idea; now here it is again in the enriched
setting'. (Indeed, the enriched setting is where some of the most interesting
examples occur.) This process of extension begins in the next two sections,
in which we discuss enriched categories and enriched operads.

Let \veeh\ be a monoidal category. A \emph{category enriched in \veeh}, or
\emph{\veeh-enriched category}, \cee, consists of a class $\ob(\cee)$ (the
\emph{objects} of \cee), an object \ehomset{\cee}{A}{B} of \veeh\ for each
$A, B \in \ob(\cee)$, and then morphisms in \veeh\ representing composition
and identities in \cee. The full definition is laid out in \cite[6.2.1]{Borx}.

\subsubsection*{Examples}

\begin{enumerate}

\item	
A \Set-enriched category is just a category.

\item	\label{eg:enr-cat:closed}
Let \veeh\ be any symmetric monoidal category which is closed, in the sense
that there is a functor
\[
[\dashbk,\dashbk] : \veeh^\op \times \veeh \go \veeh
\]
such that 
\[
\homset{\veeh}{U\otimes V}{W} \iso \homset{\veeh}{U}{[V,W]}
\]
naturally in $U,V,W \in \veeh$. Then we obtain a \veeh-enriched category,
which we also call \veeh, by putting $\ehomset{\veeh}{U}{V} = [U,V]$.

\item	\label{eg:enr-cat:Set-closed}
A particular example of~\bref{eg:enr-cat:closed} is \Set, where $[U,V]$ is
the set of functions from $U$ to $V$.

\item	\label{eg:enr-cat:Mod}
Another example of~\bref{eg:enr-cat:closed} is $\Mod_R$, where $[U,V]$ is the
usual module of homomorphisms $U\go V$.

\item	\label{eg:enr-cat:GrMod}
Another is \GrMod; this time $[U,V]_n$ is the module of degree $n$ maps from
$U$ to $V$ (i.e.\ families of homomorphisms $(U_k \go
V_{k+n})_{k\in\integers}$). 

\item	\label{eg:enr-cat:Top}
Another is \Top, recalling~(\ref{sec:prelims-mon-cats}\bref{eg:mon-cat:Top})
that we can form the function space $[U,V]$ because our spaces are assumed to
be compactly generated Hausdorff. (The square bracket notation does not mean
homotopy classes of maps, as sometimes it does in the literature.)

\item	\label{eg:enr-cat:Cat}
Another is \Cat, where $[U,V]$ is the usual category of functors from $U$ to
$V$ and natural transformations.

\item		\label{eg:enr-cat:ChCx}
Let $\veeh=\GrMod$. Then there is a \veeh-enriched category \ChCx, in which
the objects are chain complexes and if $U$ and $V$ are chain complexes,
$(\ehomset{\ChCx}{U}{V})_n$ is the module whose elements are the degree $n$
chain maps from $U$ to $V$. 

\end{enumerate}

Any \veeh-enriched category has an underlying (\Set-enriched)
category. Formally this is obtained by applying the `change of base'
$\homset{\veeh}{I}{\dashbk}: \veeh\go\Set$. Informally it's clear enough
what's going on in each of our examples of \veeh: e.g.\ if \cee\ is a
\Top-enriched category, one obtains an ordinary category simply by forgetting
the topology on each \ehomset{\cee}{A}{B} and regarding it as a mere set. In
the case $\veeh=\GrMod$ it's not quite so obvious; the analogous process is
to take each graded module \ehomset{\cee}{A}{B} and extract from it the set
\label{page:underlying-chain-complex}%
$(\ehomset{\cee}{A}{B})_0$ which is its degree 0 part. Thus the underlying
category of the \GrMod-enriched category \ChCx\ (see~\bref{eg:enr-cat:ChCx})
is the category we called \ChCx\ in
\ref{sec:prelims-mon-cats}\bref{eg:mon-cat:ChCx}. 

Now suppose that \veeh\ is a \emph{symmetric} monoidal category. Then it's
possible to define what a \emph{\veeh-enriched monoidal category} is, and
similarly a \emph{\veeh-enriched symmetric monoidal category}, a
\emph{\veeh-enriched (symmetric) monoidal functor}, and a
\emph{\veeh-enriched monoidal transformation}. (The details might be in the
encyclopaedic \cite{Kel}; I could not locate a copy.) A \Set-enriched
(symmetric) monoidal category is just a (symmetric) monoidal category, and
similarly functors and transformations. All the
examples~\bref{eg:enr-cat:closed}--\bref{eg:enr-cat:ChCx} of \veeh-enriched
categories are in fact \veeh-enriched symmetric monoidal categories in an
obvious way: in~\bref{eg:enr-cat:closed}, for instance, any symmetric
monoidal closed category is a symmetric monoidal category enriched in itself,
and in~\bref{eg:enr-cat:ChCx}, \ChCx\ is a \GrMod-enriched symmetric monoidal
category.

\section{Enriched Operads}
\label{sec:enr-ops}

We now move on to enrichment of operads. We could, more generally, talk about
enriched multicategories, but will not; it is in that context that the term
`enrichment' is most evidently appropriate.

So, let \veeh\ be a symmetric monoidal category. Then \emph{\veeh-enriched
operads} (symmetric and non-symmetric) are defined just as ordinary operads
were, except that the sets $P(n)$ are now objects of \veeh, cartesian product
$\times$ becomes $\otimes$ (the tensor product in \veeh), and the identity
element of $P(1)$ is now a map $I\go P(1)$ in \veeh. So a \Set-enriched
operad is an operad as defined above (\ref{defn:operad}), a \Top-enriched
operad is a `topological operad', a \ChCx-enriched operad is what is known as
a `differential graded operad', and so on; all of this in both symmetric and
non-symmetric flavours.

We can discuss algebras too. If $P$ is a \veeh-enriched non-symmetric operad
and \emm\ a \veeh-enriched monoidal category, then a \emph{$P$-algebra in
\emm} is an object $A$ of \emm\ together with a map $P(n) \go
\ehomset{\emm}{A^{\otimes n}}{A}$ in \veeh\ for each $n$, satisfying suitable
axioms. The symmetric case is similar. In both cases, the $P$-algebras in
\emm\ form a category $\homset{\Alg_\veeh}{P}{\emm}$ ---or just
$\homset{\Alg}{P}{\emm}$, for simplicity.

\subsubsection*{Examples}

\begin{enumerate}

\item
When $\veeh=\Set$, this is the definition of algebra given above
(\ref{defn:algebra}). 

\item	\label{eg:enr-op:Act-G}
Let $\veeh=\Top$ and let $G$ be a topological monoid: then there is a
\Top-enriched non-symmetric operad \Act{G}, defined by the formula of
Example~\ref{sec:operads}\bref{eg:operad:Act-G}. An \Act{G}-algebra in \Top\
is a space with a continuous left action by $G$. The same applies to the
symmetric version, \SAct{G}. Alternatively, we can use \Cat\ instead of \Top\
and take $G$ to be a strict monoidal category.

\item	\label{eg:enr-op:Act-G-alg}
Let $\veeh=\GrMod$ and $\emm=\ChCx$, as
in~\ref{sec:enr-cats}\bref{eg:enr-cat:ChCx}; let $G$ be a graded
algebra. Then, as in \bref{eg:enr-op:Act-G}, there is a \GrMod-enriched
non-symmetric operad \Act{G}. An algebra for \Act{G} in \ChCx\ is a chain
complex with a left action by the graded algebra $G$.

\item	\label{eg:enr-op:Lie}
Let \veeh\ be \Ab\ (abelian groups) and let \emm\ be $\Mod_R$, which is a
\veeh-enriched symmetric monoidal category in a natural way. There is an
\Ab-enriched symmetric operad \Lie\ with the property that
\homset{\Alg_\Ab}{\Lie}{\Mod_R} is the category of Lie algebras over
$R$. See~\cite[2.2]{Kap}, \cite[1.3.9]{GK}, \cite{MayOAM}, or \cite[1.5]{KSV}
for more details on \Lie.

\item	\label{eg:enr-op:GrLie}
Let $\veeh=\GrAb(=\GrMod_\integers)$ and let $\emm=\GrMod_R$. There is a
\GrAb-enriched symmetric operad \GrLie\ such that
\homset{\Alg_\GrAb}{\GrLie}{\GrMod_R} is the category of graded Lie algebras
over $R$. By a `graded Lie algebra' I mean a graded module $A$ together with
a binary operation of degree $-1$ ---that is, a family of homomorphisms
\[
[\dashbk,\dashbk]: A_p \otimes A_q \go A_{p+q-1}
\]
---satisfying suitable identities. See~\ref{subsec:GrLie} for further
   details.

(If we want the bracket to be of degree 0
then we can get away with taking \veeh\ to be the more simple category \Ab\
instead: just change \emm\ to $\GrMod_R$ in Example~\bref{eg:enr-op:Lie}.)

\item
Taking $\emm=\ChCx$ in \bref{eg:enr-op:GrLie}, \homset{\Alg}{\GrLie}{\ChCx}
is the category of differential graded Lie algebras.

\item	\label{eg:enr-op:Ger}
Let $\veeh=\GrAb$ and $\emm=\GrMod_R$. There's a certain \GrAb-enriched
symmetric operad \Ger, such that \homset{\Alg_\GrAb}{\Ger}{\GrMod_R} is the
category of Gerstenhaber algebras over $R$. A Gerstenhaber algebra (see
\cite{Vor}) is by definition a graded module which is both a
graded-commutative algebra and a graded Lie algebra, with the two structures
being compatible.

\end{enumerate}

\subsubsection*{Aside: the definition of algebra}
\small

Suppose we have a \veeh-enriched operad $P$, and wish to discuss $P$-algebras
in some kind of monoidal category \emm. In order to do this it isn't actually
necessary for \emm\ to be enriched in \veeh. For instance, if $\emm=\veeh$ is
the category of \emph{all} topological spaces, not necessarily compactly
generated Hausdorff, then one can still define a `$P$-algebra in \emm'
sensibly: it's an object $A$ of \emm\ together with a suitable family of maps
\[
P(n) \times A^n \go A.
\]
More generally, suppose that \emm\ is an (ordinary) monoidal category and
that \veeh\ `acts' on \emm, in the sense that there is a functor
\[
\begin{array}{ccc}
\veeh\times\emm		&\go		&\emm,		\\
\pr{V}{A}		&\goesto	&V\cdot A
\end{array}
\]
with suitable properties. Then one can define a $P$-algebra in \emm\ as an
object $A$ of \emm\ together with maps 
\[
P(n) \cdot A^{\otimes n} \go A
\]
satisfying axioms. For instance, any symmetric monoidal category \veeh\
(closed or not) acts on itself, and so one has a notion of $P$-algebras in
\veeh. If \veeh\ \emph{is} closed then the two notions of $P$-algebra in
\veeh, one given by enrichment and the other by action, coincide.

So we now have two possible contexts for forming a category of algebras:
\emm\ can either be enriched in \veeh\ or acted on by \veeh. How are we to
combine the two? An obvious answer is to stipulate that \emm\ is a monoidal
category and that there is a given functor
\[
H: \veeh^\op \times \emm^\op \times \emm \go \Set
\]
with suitable properties. If \emm\ is enriched in \veeh\ then $H$ arises as 
\[
\triple{V}{A}{B} \goesto \homset{\veeh}{V}{\ehomset{\emm}{A}{B}},
\]
and if \emm\ is acted on by \veeh\ then $H$ arises as
\[
\triple{V}{A}{B} \goesto \homset{\emm}{V\cdot A}{B}.
\]
For a general $H$, one ought to be able to define a category
\homset{\Alg}{P}{\emm} in a sensible way.

However, we do not take these thoughts any further in this work. By good
luck, and with the aid of devices such as compactly generated spaces,
enrichment suffices to cover all the examples that have come to mind.

\normalsize

\section{The Free Monoidal Category on an Operad}
\label{sec:free-mon}

This last preliminary section explains how an operad gives rise to a strict
monoidal category. This process was probably first described by Boardman and
Vogt; an account can also be found in the book of Adams. (See
\cite{BV} and \cite[p.~42]{Ad}.) 

Here we offer three different descriptions of the construction: the first
abstract, the last concrete, and the second somewhere in between. Then, after
some examples, we prove a crucial property of the construction, which 
provides an alternative description of algebras for an operad and is 
a conceptual stepping-stone to the definition of homotopy algebra. 

Our aim, then, is to take a (symmetric or non-symmetric) operad $P$ and
construct from it a (symmetric or not) strict monoidal category \fmon{P}. 

\paragraph*{First Description} Recall from~\ref{sec:multicats} that any
strict monoidal category $L$ has an `underlying' multicategory \ovln{L}. This
defines a functor
\[
\ovln{\blank}: (\textrm{strict monoidal categories})
\go (\textrm{multicategories}),
\]
which happens to have a left adjoint \fmon{\blank}. A non-symmetric operad
$P$ is just a one-object multicategory, and so we obtain from $P$ a strict
monoidal category \fmon{P}. The same goes in the symmetric case.

(This adjunction is discussed in more depth and generality
in~\cite[4.3]{GOM}.) 

\paragraph*{Second Description} Let $P$ be a non-symmetric operad; again, what
follows can be repeated with the obvious changes for symmetric operads. Then
\fmon{P} can be constructed as the free strict monoidal category containing a
$P$-algebra. Consider, by way of comparison, the category $\Delta$ of
Example~\ref{sec:prelims-mon-cats}\bref{eg:mon-cat:Delta}. $\Delta$ can be
described as the free strict monoidal category containing a monoid: for it is
generated (as a strict monoidal category) by the objects and morphisms
\[
0 \goby{e} 1 \ogby{m} 1+1,
\]
subject to the usual monoid axioms on $m$ and $e$. (See \cite[VII.5.1]{CWM}.)
To build our category \fmon{P}, we must first of all put into it an object
$A$, the underlying object of the $P$-algebra it contains; then, since
\fmon{P} is a monoidal category, it must contain the $n$th tensor power
$A^{\otimes n}$ for each $n \geq 0$; then, since $A$ is meant to be a
$P$-algebra in \fmon{P}, there must be a morphism $A^{\otimes n} \go A$ in
\fmon{P} for each element of $P(n)$; and so on.

In the Third Description below, which is the one we will actually use, the
object $A^{\otimes n}$ of \fmon{P} is written merely as $n$, and $\otimes$ is
therefore written as $+$.

\paragraph*{Third Description} Concretely, let $P$ be a non-symmetric operad,
and define a monoidal category \fmon{P} as follows.
The objects are the natural numbers $0, 1, \ldots$, and the monoidal
structure on the objects is addition, with unit $0$. The homsets are
given by
\[
\homset{\fmon{P}}{m}{n} = 
\coprod_{m_1 + \cdots + m_n = m}
P(m_1) \times \cdots \times P(m_n)
\]
for $m,n\in\nat$. Thus if $\theta_1 \in P(m_1), \ldots, \theta_n \in P(m_n)$
and $m_1 + \cdots + m_n = m$, there is an element
\bftuple{\theta_1}{\theta_n} of \homset{\fmon{P}}{m}{n}. The
tensor product of morphisms is defined by
\[
\bftuple{\theta_1}{\theta_n} \otimes \bftuple{\theta'_1}{\theta'_{n'}}
=
\tuplebts{\theta_1, \ldots, \theta_n, \theta'_1, \ldots \theta'_{n'}}.
\]
It remains to describe the identity and composition in \fmon{P} (and to check
all the axioms). The identity element of \homset{\fmon{P}}{m}{m} is
\bftuple{1}{1}, consisting of $m$ copies of the unit $1$ of $P$. For
composition, take $\phi_1 \in P(k_1), \ldots, \phi_m \in P(k_m)$ with $k_1 +
\cdots + k_m = k$, so that $\bftuple{\phi_1}{\phi_m} \in
\homset{\fmon{P}}{k}{m}$, and take $\bftuple{\theta_1}{\theta_n}
\in \homset{\fmon{P}}{m}{n}$ as above. Since $m_1 + \cdots + m_n
= m$, we may rewrite the sequence \bftuple{k_1}{k_m} as \tuplebts{k_1^1,
\ldots, k_1^{m_1}, \ldots, k_n^1, \ldots, k_n^{m_n}}, and similarly
\bftuple{\phi_1}{\phi_m} as \tuplebts{\phi_1^1, \ldots, \phi_1^{m_1}, \ldots,
\phi_n^1, \ldots, \phi_n^{m_n}}.  Thus $\phi_i^j \in P(k_i^j)$. We then have
\begin{eqnarray*}
\theta_1 \of \bftuple{\phi_1^1}{\phi_1^{m_1}} &\in& 
P(k_1^1 + \cdots + k_1^{m_1}), \ldots,\\
\theta_n \of \bftuple{\phi_n^1}{\phi_n^{m_n}} &\in& 
P(k_n^1 + \cdots + k_n^{m_n}),
\end{eqnarray*}
and the composite $\bftuple{\theta_1}{\theta_m} \of \bftuple{\phi_1}{\phi_n}$
is defined as 
\[
\bftuple{\theta_1 \of \bftuple{\phi_1^1}{\phi_1^{m_1}}}{\theta_n \of
\bftuple{\phi_n^1}{\phi_n^{m_n}}}.
\]
This is indeed an element of \homset{\fmon{P}}{k}{n}, since 
\[
(k_1^1 + \cdots + k_1^{m_1}) + \cdots + (k_n^1 + \cdots + k_n^{m_n})
=
k_1 + \cdots + k_m
= 
k.
\]

Similarly, let $P$ be a symmetric operad: then there is an associated
symmetric strict monoidal category, described shortly, which we also write as
\fmon{P}. This is possibly an abuse of language, because \fmon{P} is different
depending on whether $P$ is considered with or without its symmetric
structure, but I hope that context will always make things clear.

The objects of \fmon{P} are again the natural numbers, with monoidal
structure given by $+$ and $0$. The homsets are given by
\[
\homset{\fmon{P}}{m}{n} =
\coprod_{f\in \homset{\Phi}{m}{n}}
P(f^{-1}\{0\}) \times \cdots \times P(f^{-1}\{n-1\}).
\]
Here $\Phi$ is (a skeleton of) the category of finite sets, defined in
\ref{sec:prelims-mon-cats}\bref{eg:mon-cat:Phi}, and we write $P(S)$ to mean
$P(s)$ when $s$ is the cardinality of a finite set $S$. Thus if we replace
$\Phi$ by $\Delta$ in this formula, we obtain the definition of
\homset{\fmon{P}}{m}{n} in the non-symmetric version. If $f\in
\homset{\Phi}{m}{n}$ and if
\[
\theta_1 \in P(f^{-1}\{0\}), \ldots, \theta_n \in P(f^{-1}\{n-1\})
\] 
then the corresponding element of $\homset{\fmon{P}}{m}{n}$
is written $\tuplebts{f; \theta_1, \ldots, \theta_n}$; the tensor of
morphisms is defined by
\[
\tuplebts{f; \theta_1, \ldots, \theta_n}
\otimes
\tuplebts{f'; \theta'_1, \ldots, \theta'_{n'}}
=
\tuplebts{f+f'; \theta_1, \ldots, \theta_n, \theta'_1, \ldots, \theta'_{n'}}.
\]
The symmetry on \fmon{P} is given by the element \tuplebts{t_{m,n}; 1, \dots,
1} of \homset{\fmon{P}}{m+n}{n+m}, where $t_{m,n} \in
\homset{\Phi}{m+n}{n+m}$ adds $n$ to the first $m$ elements of $m+n$ and
subtracts $m$ from the last $n$, and where there are $m+n$ copies of $1=1_P$
after the semicolon. The identity morphism $1_m \in \homset{\fmon{P}}{m}{m}$
is given by the identities in $\Phi$ and in $P$; composition in \fmon{P} is
described more or less as in the non-symmetric version, but with some
slightly intricate manipulation of permutations. (A closely related but
slightly different construction is detailed in \cite[4.1]{MT}, and this gives
an impression of the method involved. In the formula there for
\homset{\widehat{\cee}}{m}{n}, the first $\prod$ should be a $\coprod$.)

\subsubsection*{Examples}

\begin{enumerate}

\item
Let \Mon\ be the non-symmetric operad whose algebras are monoids, as
in~\ref{sec:operads}\bref{eg:operad:Mon}. Then \fmon{\Mon} is $\Delta$, the
skeleton of the category of finite totally ordered sets (defined
in~\ref{sec:prelims-mon-cats}\bref{eg:mon-cat:Delta}), with $+$ and $0$ as
its monoidal structure.

\item	\label{eg:free-mon:Sem}
Similarly, take the non-symmetric operad \Sem\
of~\ref{sec:operads}\bref{eg:operad:Sem}: then \fmon{\Sem} is
$\Delta_{\mr{surj}}$, the subcategory of $\Delta$ consisting of all its
objects $n$ ($n\geq 0$) but only the \emph{surjective} order-preserving maps.

\item	\label{eg:free-mon:Pt}
As some kind of dual to~\bref{eg:free-mon:Sem}, $\fmon{\Pt}=\Delta_{\mr{inj}}$,
where \Pt\ is the non-symmetric operad
of~\ref{sec:operads}\bref{eg:operad:Pt} and $\Delta_{\mr{inj}}$ is the
subcategory of $\Delta$ made up from injective maps.

\item
Let \CMon\ be the symmetric operad for commutative monoids
(\ref{sec:operads}\bref{eg:operad:CMon}): then \fmon{\CMon} is $\Phi$, the
skeleton of the category of finite sets (defined
in~\ref{sec:prelims-mon-cats}\bref{eg:mon-cat:Phi}).

\item	\label{eg:free-mon:CSem}
As in~\bref{eg:free-mon:Sem} and~\bref{eg:free-mon:Pt},
$\fmon{\CSem}=\Phi_{\mr{surj}}$ and $\fmon{\SPt}=\Phi_{\mr{inj}}$, where
$\Phi_{\mr{surj}}$ and $\Phi_{\mr{inj}}$ are defined in the obvious way.

\item	\label{eg:free-mon:Act-G}
\setcounter{bean}{\value{enumi}}
If $G$ is a monoid and \Act{G} the non-symmetric operad whose algebras are
left $G$-objects (\ref{sec:operads}\bref{eg:operad:Act-G}) then
\fmon{\Act{G}} is the monoidal category with objects $0,1,\ldots$ and
\[
\homset{\fmon{\Act{G}}}{m}{n} = 
\left\{
\begin{array}{ll}
G^n		&\textrm{if }m=n	\\
\emptyset	&\textrm{otherwise.}
\end{array}
\right.
\]
Composition and identities are as in $G$, and tensor of morphisms is
juxtaposition. 

If we take the symmetric operad \SAct{G} instead, then
\[
\homset{\fmon{\SAct{G}}}{m}{n} = 
\left\{
\begin{array}{ll}
G^n \times S_n	&\textrm{if }m=n	\\
\emptyset	&\textrm{otherwise.}
\end{array}
\right.
\]

\end{enumerate}

In the Second Description we characterized \fmon{P} as `the free (symmetric)
monoidal category containing a $P$-algebra'. This was meant in a syntactic
sense: \fmon{P} is generated by certain objects and morphisms subject to
certain equations. But this characterization can also be interpreted as a
universal property: if \emm\ is any (symmetric) strict monoidal category then
$P$-algebras in \emm\ correspond one-to-one with (symmetric) strict monoidal
functors $\fmon{P} \go \emm$. In fact, the correspondence extends to the
non-strict situation, as stated in the following important result.

\begin{thm}	\label{thm:alt-alg}
\begin{enumerate}
\item	\label{part:non-sym}
Let $P$ be a non-symmetric operad and \emm\ a monoidal category. Then there
is an equivalence of categories 
\[
\homset{\Alg}{P}{\emm} \eqv \homset{\Mon}{\fmon{P}}{\emm}.
\] 
\item	\label{part:sym}
Let $P$ be a symmetric operad and \emm\ a symmetric monoidal category. Then
there is an equivalence of categories 
\[
\homset{\Alg}{P}{\emm} \eqv \homset{\SMon}{\fmon{P}}{\emm}.
\]
\end{enumerate}
\end{thm}

\begin{sketchpf}
For~\bref{part:non-sym}, take a $P$-algebra in \emm, consisting of an
object $A$ of \emm\ and a map $\ovln{\theta}: A^{\otimes n} \go A$ for each
$\theta\in P(n)$. Then there arises a functor $X: \fmon{P} \go \emm$ given by
setting $X(n)=A^{\otimes n}$ on objects, and by setting
\[
X(\theta_1, \ldots, \theta_n) =
\ovln{\theta_1} \otimes\cdots\otimes \ovln{\theta_n}:
A^{\otimes m_1} \otimes\cdots\otimes A^{\otimes m_n}
\go A^{\otimes n}
\]
for any $\theta_1\in P(m_1), \ldots, \theta_n\in P(m_n)$ making up a map $m_1
+\cdots + m_n \go n$ in \fmon{P}. This functor $X$ has a natural monoidal
structure.

Conversely, take a monoidal functor $\pr{X}{\xi}: \fmon{P}\go\emm$. Put
$A=X(1)$. The components of $\xi$ fit together (in exactly one way) to
produce an isomorphism $\xi^{(n)}: X(n) \goiso X(1)^{\otimes n}$, for any
given $n$. Thus if $\theta\in P(n)$ then we may define $\ovln{\theta}$ to be
the composite
\[
A^{\otimes n} \goby{(\xi^{(n)})^{-1}}
X(n) \goby{X(\theta)}
X(1) = A,
\]
where we are regarding $\theta$ as an element of
\homset{\fmon{P}}{n}{1}. This defines an algebra structure on $A$. 

We have to prove that two categories are equivalent, and have shown how to
pass from an object of either category to an object of the other. These
processes extend to morphisms in a straightforward way, and the two functors
so defined are mutually inverse up to natural isomorphism.

Part~\bref{part:sym} is just a more elaborate version
of~\bref{part:non-sym}. The trickiest moment is in obtaining the functor
$X:\fmon{P} \go \emm$ arising from a $P$-algebra \emm: one must use some
permutations to define $X$ on morphisms, and then check that $X$ really is a
functor.
\done
\end{sketchpf}
\subsubsection*{Examples}
\begin{enumerate}
\setcounter{enumi}{\value{bean}}
\item
In the case of the non-symmetric operad \Mon, the Theorem says that the
category of monoids in a monoidal category \emm\ is equivalent to the
category of monoidal functors $\Delta\go\emm$. This is very well-known: see
\cite[VII.5.1]{CWM}. 

\item
Similarly, in the symmetric case, taking $P=\CMon$ tells us that a
commutative monoid in a symmetric monoidal category \emm\ is essentially the
same thing as a symmetric monoidal functor $\Phi\go\emm$.

\end{enumerate}

We finish by observing that the theory above can be generalized in two
directions. Firstly, `operad' can be replaced by `multicategory' ($=$
`coloured operad') in the Theorem, provided that \fmon{P} is defined
correctly (for which see the First Description above). We shall not need this
generalization.

Secondly, we can extend to the situation where all the operads and monoidal
categories concerned are enriched in a suitable symmetric monoidal category
\veeh. Thus if \cat{L} and \emm\ are \veeh-enriched monoidal categories,
there is an (ordinary) category \homset{\Mon}{\cat{L}}{\emm} of
\veeh-enriched monoidal functors $\cat{L}\go\emm$ and \veeh-enriched monoidal
transformations, and similarly \homset{\SMon}{\cat{L}}{\emm} in the symmetric
case: see~\ref{sec:enr-cats}. Now assume that \veeh\ has finite coproducts
and that $\otimes$ distributes over them, as is the case for the \veeh\ in
each of
Examples~\ref{sec:enr-cats}\bref{eg:enr-cat:Set-closed}--\bref{eg:enr-cat:Cat}.
Then for any \veeh-enriched operad $P$, a \veeh-enriched (symmetric) monoidal
category \fmon{P} can be defined just as in the non-enriched version above
(in the Third Description), replacing $\times$ by $\otimes$ and $1$ by $I$,
and reading $\coprod$ as coproduct in \veeh. Theorem~\ref{thm:alt-alg}
follows, with $P$ and \emm\ both \veeh-enriched. Thus a $P$-algebra in \emm\
is essentially the same thing as a \veeh-enriched (symmetric) monoidal
functor $\fmon{P} \go \emm$.

\chapter{The Definition of Homotopy Algebra}
\label{ch:defn}

The path to defining homotopy algebras is now clear. The key is
Theorem~\ref{thm:alt-alg}, which gave an alternative description of an
algebra for an operad: namely, if $P$ is an operad and \emm\ a monoidal
category, then a $P$-algebra in \emm\ is a functor $X: \fmon{P} \go \emm$
together with a coherent family of isomorphisms
\[
X(m+n) \go X(m)\otimes X(n),
\diagspace
X(0) \go I.
\]
To define `homotopy $P$-algebra in \emm', we simply take this description and
change the word `isomorphisms' to `homotopy equivalences': and that is our
definition.

In order for this to make sense, we must of course have some notion of what a
`homotopy equivalence' in \emm\ is. Since in a naked monoidal category there
is no \emph{a priori} notion of homotopy, this is tagged on as extra
structure. Thus we consider a monoidal category \emm\ equipped with a class
of its morphisms, called the `homotopy equivalences'; and \emm\ equipped with
this extra structure is a suitable environment in which to define homotopy
algebras for an operad. (Naturally enough, we insist that the class of
homotopy equivalences obeys a few axioms such as closure under composition,
so that the definition behaves reasonably.)

This method of capturing the notion of homotopy is very crude, and
consequently the definition of homotopy algebra is a crude one. Taking \emm\
to be the monoidal category of topological spaces, for instance, we have
simply recorded which continuous maps are homotopy equivalences. We have no
notion of what it means for one map to be homotopy-inverse to another, or for
two maps to be homotopic, or for two homotopies to be homotopic, and so
on. Thus we are missing out a vast amount of the homotopy theory of
spaces---all we have is the `1-dimensional trace' of a whole
$\infty$-category of information---and we should therefore expect our
definition of homotopy algebra to have certain shortcomings.

On the other hand, the definition has some advantages. Simplicity is
one. Another is historical precedent: a homotopy topological commutative
monoid will turn out to be exactly the same as a (special) $\Gamma$-space in
the sense of Segal's paper \cite{SegCCT}. (I call this an `advantage' because
$\Gamma$-spaces have already been well explored, thus providing a firm
attachment between the present definition and established topology.)
Moreover, despite the fact that our definition is only a dim reflection of
the fully glorious (and as yet unformulated) $\infty$-categorical definition,
it can at least be seen as the first rung on a ladder leading up to this
ideal.

A further discussion of how the definition fits into the big picture,
including some more on $\infty$-categories, can be found in the conclusion,
Chapter~\ref{ch:thoughts}.  

This chapter is laid out as follows. In~\ref{sec:env} we make precise the
notion of a monoidal category with a class of equivalences (the `environment'
in which homotopy algebras are taken), and run through some
examples. Section~\ref{sec:defn} consists of the definition of homotopy
algebra. In~\ref{sec:brief} we take a first look at some examples: both some
rather trivial ones, and sketches of some more substantial ones which are the
subject of later chapters. All of this so far is for the non-enriched case
(i.e.\ for operads $P$ in which each $P(n)$ is a mere set, with no extra
structure); but in~\ref{sec:enr-defn} we extend the definition to the
enriched setting.

\section{The Environment}
\label{sec:env}

\begin{defn}	\label{defn:equivs}
A \emph{monoidal category with equivalences} is a monoidal category \emm\
equipped with a subclass \Eee\ of the morphisms in \emm, whose elements are
called \emph{equivalences} or \emph{homotopy equivalences}, such that the
following properties hold:
\begin{description}
\item[E1]
any isomorphism is an equivalence
\item[E2]
if $h=g\of f$ is a composite of morphisms in \emm, and if any two of
$f, g, h$ are equivalences, then so is the third
\item[E3]
if $A\goby{f}B$ and $A'\goby{f'}B'$ are equivalences then so is
$A\otimes A' \goby{f\otimes f'} B\otimes B'$.
\end{description}
If \emm\ is a symmetric monoidal category, then \emm\ together with \Eee\
forms a \emph{symmetric monoidal category with equivalences}. In both cases,
we call \Eee\ a \emph{class of equivalences in \emm}. 
\end{defn}

\subsubsection*{Examples}

\begin{enumerate}

\item	\label{eg:env:isos}
Let \emm\ be any monoidal category and let \Eee\ be the class of isomorphisms
in \emm. This is (by \textbf{E1}) the smallest possible class
of equivalences in \emm.

\item
Dually, if \emm\ is any monoidal category then taking \Eee\ to be all
morphisms in \emm\ gives the largest possible class of equivalences in \emm. 

\item	\label{eg:env:Cat}
\emm\ is \triple{\Cat}{\times}{\One}, and equivalences are equivalences of
categories: that is, those functors $G$ for which there exists some functor
$F$ with $F\of G \iso 1$ and $G\of F \iso 1$. ($F$ is called a
\emph{pseudo-inverse} to $G$). 

\item	\label{eg:env:Top-times}
\emm\ is \triple{\Top}{\times}{1}, and equivalences are homotopy
equivalences. 

\item	\label{eg:env:Topstar-times}
\emm\ is \triple{\Topstar}{\times}{1}, and equivalences are homotopy
equivalences relative to basepoints.

\item
Example~\bref{eg:env:Topstar-times} can be repeated with the wedge product
\wej\ in place of $\times$. 

\item	\label{eg:env:ChCx}
\emm\ is \triple{\ChCx}{\otimes}{R}
(see~\ref{sec:prelims-mon-cats}\bref{eg:mon-cat:ChCx}); equivalences are
chain homotopy equivalences.

\item
The reader might be wondering whether, in~\bref{eg:env:ChCx}, we could have
taken quasi-isomorphisms in place of chain homotopy equivalences. (A chain
map is called a \emph{quasi-isomorphism} if it induces an isomorphism on each
homology group.) Axioms \textbf{E1} and \textbf{E2} are easily verified, but
\textbf{E3} is more demanding. Consider the commutative square
\begin{diagram}
H_{\blob}(A) \otimes H_{\blob}(A')	&
\rTo^{f_* \otimes f'_*}			&
H_{\blob}(B) \otimes H_{\blob}(B')	\\
\dTo	&	&\dTo			\\
H_{\blob}(A\otimes A')			&
\rTo_{(f\otimes f')_*}			&
H_{\blob}(B\otimes B')			\\
\end{diagram}
of graded modules, in which the vertical maps are the natural ones. We know
that the map along the top is an isomorphism, and would like to conclude that
the map along the bottom is an isomorphism. This will be true if the vertical
maps are also isomorphisms, which in turn is true if the ground ring $R$ is a
field (by the K\"{u}nneth Theorem, \cite[3.6.3]{Wei}). So if $R$ is a field
then the quasi-isomorphisms form a class of equivalences in $\ChCx_R$.

However, this gives us little or nothing more than
Example~\bref{eg:env:ChCx}: for when one is working over a field,
quasi-isomorphisms are (almost?)\ the same thing as chain homotopy
equivalences. More precisely, I am informed that they are the same thing if
either the field is of characteristic $0$, or if the complexes concerned are
$0$ in negative degrees; I do not know if the statement is true in complete
generality. Whatever the truth, quasi-isomorphisms will not be mentioned
again. 

\item
Similarly, in~\bref{eg:env:Top-times} and~\bref{eg:env:Topstar-times} we can
take weak homotopy equivalences rather than homotopy equivalences as long as
we work over a field, but this does not seem to provide a significant
generalization. 

\item
If \Eee\ is a class of equivalences in a monoidal category \emm, then \Eee\
is also a class of equivalences in the opposite category $\emm^\op$.

\item	\label{eg:env:mon-2-cat}
\setcounter{bean}{\value{enumi}}
Let \cat{N} be a monoidal 2-category---that is, a \Cat-enriched monoidal
category. Then \cat{N} consists of \emph{0-cells} (or \emph{objects}),
\emph{1-cells}, and \emph{2-cells}, together with various ways of composing
them, a tensor product, and a unit object. There is a notion of a 1-cell $G:
A\go B$ in \cat{N} being an \emph{equivalence}: namely, if there exists a
1-cell $F: B\go A$, an invertible 2-cell between $G\of F$ and $1_B$, and an
invertible 2-cell between $F\of G$ and $1_A$. The underlying (1-)category
\emm\ of \cat{N}, formed by the 0-cells and 1-cells, is then a monoidal
category with equivalences. A typical example is \Cat\ itself:
see~\bref{eg:env:Cat}. In fact,
Examples~\bref{eg:env:Top-times}--\bref{eg:env:ChCx} all arise in this way
too, as we shall see in Sections~\ref{sec:A-spaces} and~\ref{sec:A-algs}.

\end{enumerate}

\section{The Definition}
\label{sec:defn}

In order to make the definition of a homotopy algebra for an operad, we will
need a notion of map between monoidal categories which is more general than
the notion of monoidal functor (\ref{defn:mon-ftr}). The notion is that of a
\emph{colax monoidal functor}, and the definition is obtained from
Definition~\ref{defn:mon-ftr} simply by replacing the word `isomorphisms' with
`maps'. We get the definition of \emph{colax symmetric monoidal functor}
from~\ref{defn:mon-ftr} in the same way.

Thus in a (symmetric) colax monoidal functor $\pr{X}{\xi}: \cat{L} \go \emm$,
we have coherence maps
\[
\xi_0: X(I) \go I, 
\diagspace
\xi_{A,B}: X(A\otimes B) \go X(A) \otimes X(B)
\]
($A,B\in \cat{L}$). We shall also refer in passing to 
\label{page:defn-lax}%
\emph{lax (symmetric) monoidal functors}, in which these maps $\xi$ go in the
opposite direction. One explanation of the terminology is that a lax monoidal
functor from \cat{L} to \emm\ induces a functor from the category of monoids
in \cat{L} to the category of monoids in \cat{M}, whereas a \emph{co}lax
monoidal functor induces a functor between the categories of \emph{co}monoids.

Note that a monoidal functor is a special kind of colax monoidal functor, not
vice-versa. Thus the role of the adjective (`colax') is contrary to normal
English usage. Note also that the definition of monoidal transformation
(\ref{defn:mon-transf}) makes sense for colax monoidal functors in general. 

We now present the main definition of this paper. 

\begin{defn}	\label{defn:hty-alg}
\begin{enumerate}

\item	\label{defn:non-sym-hty-alg}
Let $P$ be a non-symmetric operad and let \emm\ be a monoidal category with
equivalences. A \emph{homotopy $P$-algebra in \emm} is a colax monoidal
functor $\pr{X}{\xi}: \fmon{P}\go\emm$ in which $\xi_0$ and each
$\xi_{m,n}$ ($m,n\in\nat$) are equivalences. (Here \fmon{P} is the monoidal category of Section~\ref{sec:free-mon}.)

\item	\label{defn:sym-hty-alg}
Let $P$ be a symmetric operad and let \emm\ be a symmetric monoidal category
with equivalences. A \emph{homotopy $P$-algebra in \emm} is a colax symmetric
monoidal functor $\pr{X}{\xi}: \fmon{P}\go\emm$ in which $\xi_0$ and each
$\xi_{m,n}$ ($m,n\in\nat$) are equivalences. (Here \fmon{P} is the
\emph{symmetric} monoidal category of~\ref{sec:free-mon}.)

\end{enumerate}
\end{defn}

In both symmetric and non-symmetric cases, a \emph{map of homotopy
$P$-algebras} is a monoidal transformation, and the homotopy $P$-algebras in
\emm\ thus form a category \homset{\HtyAlg}{P}{\emm}.

\section{Brief Examples}
\label{sec:brief}

Most of the rest of this paper consists of examples of homotopy
algebras. Each non-trivial example takes a while to explain, so for now we
just present the trivial cases and briefly sketch out the more substantial
examples. 

\begin{enumerate}

\item	\label{eg:hty-alg:isos}
Suppose that the only equivalences in \emm\ are the isomorphisms, as
in~\ref{sec:env}\bref{eg:env:isos}. Then a homotopy $P$-algebra in \emm\ is
essentially just a $P$-algebra in \emm, by Theorem~\ref{thm:alt-alg}. In
symbols,
\[
\homset{\Alg}{P}{\emm} \eqv
\homset{\Mon}{\fmon{P}}{\emm} =
\homset{\HtyAlg}{P}{\emm}.
\]
So in an \emm\ with `no interesting homotopy', homotopy algebras are just
algebras. This holds in both the symmetric and the non-symmetric case.

\item	\label{eg:hty-alg:real-to-hty}
Take any $P$ and \emm\ as in
Definition~\ref{defn:hty-alg}~(\bref{defn:non-sym-hty-alg}
or~\bref{defn:sym-hty-alg}). Then by axiom \textbf{E1} for a class of
equivalences, any $P$-algebra is a homotopy $P$-algebra. More precisely,
there is an inclusion as shown:
\[
\homset{\Alg}{P}{\emm} \eqv
\homset{\Mon}{\fmon{P}}{\emm} \rIncl
\homset{\HtyAlg}{P}{\emm}.
\]

\item	\label{eg:hty-alg:Obj}
Let $P=\Obj$ (see~\ref{sec:operads}\bref{eg:operad:Obj}) and let \emm\ be any
monoidal category with equivalences. An \Obj-algebra in \emm\ is just an
object of \emm; what is a homotopy \Obj-algebra? Roughly speaking, when
$\emm=\Top$, for instance, a homotopy \Obj-algebra consists of a space $A$
together with a homotopy model for each power $A^n$ of $A$.

In detail, \fmon{\Obj} is the discrete category \nat\ (all morphisms are
identities), so a colax monoidal functor $\fmon{\Obj} \go \emm$ consists of a
sequence $X(0), X(1), \ldots$ of objects of \emm, together with maps
\[
\xi_{m,n}: X(m+n) \go X(m) \otimes X(n),
\diagspace
\xi_0: X(0) \go I
\]
satisfying coherence axioms. These axioms guarantee that for each sequence
$k_1, \ldots k_n$ (with $n\geq 0$, $k_i \geq 0$), there is a unique map
\[
\xi_{k_1, \ldots, k_n}:
X(k_1 + \cdots k_n) \go
X(k_1) \otimes\cdots\otimes X(k_n)
\]
built up from the $\xi_{m,n}$'s and $\xi_0$. (The notations $\xi_0$ and
$\xi_{k_1, \ldots, k_n}$ conflict, but this should not cause serious
problems.) In particular, taking all the $k_i$'s to be $1$ yields a canonical
map
\[
\xi^{(n)}: X(n) \go X(1)^{\otimes n}.
\]

A homotopy \Obj-algebra is a colax monoidal functor \pr{X}{\xi} as above,
with the property that $\xi_0$ and each $\xi_{m,n}$ are equivalences. Using
the axioms on equivalences (\ref{defn:equivs}), this property can be restated
in two different ways: that each $\xi_{k_1, \ldots, k_n}$ is an equivalence,
or that each $\xi^{(n)}: X(n) \go X(1)^{\otimes n}$ is an equivalence. 


Throughout this work we think of $X(1)$ as the `base object' of a homotopy
algebra \pr{X}{\xi}, and in fact we have already encountered this idea in the
context of genuine algebras (see the proof of~\ref{thm:alt-alg}). 
If $P$ is any operad, \emm\ a (symmetric) monoidal category, and $A$ an
object of \emm, we will write `\emph{$A$ is a homotopy $P$-algebra}' to mean
that there is a homotopy $P$-algebra \pr{X}{\xi} with $X(1)\iso A$.

\item
A similar analysis can be made of homotopy algebras for the symmetric operad
\SObj\ (defined in~\ref{sec:operads}\bref{eg:operad:Obj}). In this case the
maps $\xi_{m,n}$ are compatible with the symmetries in \emm, in the sense of
Definition~\ref{defn:mon-ftr}. Hence the maps $\xi^{(n)}: X(n) \go
X(1)^{\otimes n}$ are also compatible with the symmetries, in the obvious
sense.

\item	\label{eg:hty-alg:hty-mon-cat}
Let $\emm=\Cat$ and let
$P=\Mon$~(\ref{sec:prelims-mon-cats}\bref{eg:mon-cat:Cat}
and~\ref{sec:operads}\bref{eg:operad:Mon}). A $P$-algebra in \emm\ is a
monoid in \Cat, that is, a strict monoidal category. Homotopy algebras are
meant to be some weakened version of genuine algebras, so a homotopy monoid
($=$ homotopy \Mon-algebra) in \Cat\ ought to be something comparable to a
(non-strict) monoidal category. Similarly, a homotopy \CMon-algebra in \Cat\
should be something along the lines of a (non-strict) symmetric monoidal
category. We look at these weakened notions of monoidal category in
Sections~\ref{sec:hty-mon-cats} and~\ref{sec:inside-Cat}. In particular, we
will see that a homotopy commutative monoid in \Cat\ is exactly what Segal
called a $\Gamma$-category in \cite{SegCCT}. (We will call these things
\emph{special} $\Gamma$-categories instead, following the more popular
terminology.) 

\item
A prime example of something which ought to be a homotopy monoid is a loop
space. It is, as is proved in~\ref{sec:loops}. More precisely, we prove that
for any based space $B$ there is a homotopy \Mon-algebra \pr{X}{\xi} in
\triple{\Top}{\times}{1} with $X(1)$ isomorphic to the space of based loops
in $B$; cf.\ the remarks at the end of Example~\bref{eg:hty-alg:Obj} above.

\item	
Taking \emm\ to be \triple{\Top}{\times}{1} again, we might well expect a
homotopy monoid in \emm\ to be something like an $A_\infty$-space (as defined
in \cite{HAHI}). Both concepts are, after all, meant to provide an
up-to-higher-homotopy version of topological
monoid. Section~\ref{sec:A-spaces} provides a partial comparison between the
two. More accurately, the comparison is between $A_\infty$-spaces and
homotopy \emph{semigroups} in \Topstar\ (the category of based spaces
(\ref{sec:env}\bref{eg:env:Topstar-times})): there is a slightly delicate
issue concerning spaces with or without basepoint and semigroups with or
without unit, which is explained there.

We will also see in Section~\ref{sec:Gamma} that a homotopy commutative
monoid in \Top\ is precisely a $\Gamma$-space (as defined in \cite{SegCCT}),
or `special $\Gamma$-space' in the alternative terminology. It is from this
re-definition of (special) $\Gamma$-space that the general definition of
homotopy algebra descended.

\item
Let \emm\ be the category \ChCx\ of chain complexes, with usual tensor and
homotopy equivalences, as in~\ref{sec:env}\bref{eg:env:ChCx}. Then a monoid
in \emm\ is a differential graded algebra, so a homotopy monoid in \emm\
should be something comparable to an $A_\infty$-algebra (as defined in
\cite{HAHII}). A comparison of sorts is made in Section~\ref{sec:A-algs}.

\item	\label{eg:hty-alg:Act-G}
The following example suggests that the definition of homotopy algebra does
not encompass as much as we might like. Fix a monoid $G$, let $\Act{G}$ be
the non-symmetric operad of Example~\ref{sec:operads}\bref{eg:operad:Act-G},
and let \emm\ be any monoidal category. An $\Act{G}$-algebra in \emm\ is an
object of \emm\ equipped with a left action by $G$, and we might therefore
expect a homotopy $\Act{G}$-algebra to be an object with an `action up to
homotopy', so that laws like $g\cdot (g'\cdot x) = (gg')\cdot x$ only hold in
some weak sense. This is not the case. For by the description of
\fmon{\Act{G}} in~\ref{sec:free-mon}\bref{eg:free-mon:Act-G}, a homotopy
$\Act{G}$-algebra consists of a sequence $X(0), X(1), \ldots$ of objects of
\emm, with a \emph{strict} action of $G^n$ on $X(n)$ for each $n$, and
homotopy equivalences $\xi_{m,n}, \xi_0$ (as in~\bref{eg:hty-alg:Obj}) which
preserve the $G^n$-actions. In particular, the `base' object $X(1)$ has a
strict action by $G$, which might be a disappointment in cases such as \Cat\
and \Top. This matter is discussed further in Chapter~\ref{ch:thoughts},
together with the related matter of homotopy invariance.

\end{enumerate}

By way of advertisement, one can perform various `changes of environment':
for instance, the classifying-space functor $B: \Cat\go\Top$ induces a map
from homotopy monoids in \Cat\ to homotopy monoids in \Top. This is the
subject of Chapter~\ref{ch:change}.

\section{The Enriched Version}	\label{sec:enr-defn}

Homotopy algebras can be defined in the enriched context too
(\ref{sec:enr-cats},~\ref{sec:enr-ops}). 

First of all, we need an enriched version of `monoidal categories with
equivalences'. Fix a symmetric monoidal category \veeh. Then a \emph{class of
equivalences} in a \veeh-enriched monoidal category \emm\ is simply a class
of equivalences in the underlying monoidal category $\und{\emm}$ of \emm.

\subsubsection*{Examples}

\begin{enumerate}
\item 
Let $\veeh=\Ab$. Then a \veeh-enriched monoidal category with equivalences is
an \Ab-enriched monoidal category \emm\ together with a sub\emph{set} (not
necessarily a sub\emph{group}) of \ehom{\emm}{A}{B}, for each $A,B\in\emm$,
whose elements are called the `equivalences' from $A$ to $B$, satisfying the
axioms \textbf{E1}--\textbf{E3} of Definition~\ref{defn:equivs}.

\item 
\setcounter{bean}{\value{enumi}}
Let $\veeh=\GrMod$ and $\emm=\ChCx$, as
in~\ref{sec:enr-cats}\bref{eg:enr-cat:ChCx}. The underlying (ordinary)
monoidal category \und{\emm} of \emm\ is the category also denoted by \ChCx,
and the (degree $0$) chain homotopy equivalences provide a class of
equivalences in \und{\emm} (by~\ref{sec:env}\bref{eg:env:ChCx}) and hence in
\emm. 
\end{enumerate}

I would now like to define homotopy algebras in the enriched setting. It may
be that the reader has no head or stomach for the niceties of enriched
category theory, in which case he should jump straight to the
examples. Otherwise, take a \veeh-enriched operad $P$ (symmetric or not) and
a \veeh-enriched (symmetric) monoidal category with equivalences, \emm. A
homotopy $P$-algebra in \emm\ will be defined as a \veeh-enriched colax
(symmetric) monoidal functor
\[
\pr{X}{\xi}: \fmon{P} \go \emm,
\]
where \fmon{P} is as at the end of~\ref{sec:free-mon}, satisfying certain
conditions. These conditions should say `the components of $\xi$ are
equivalences'. Now, the components of $\xi$ are maps
\[
\begin{array}{rccc}
\xi_{m,n}:	&I	&\go	
&\ehom{\emm}{X(M+n)}{X(m)\otimes X(n)},	\\
\xi_0:		&I	&\go
&\ehom{\emm}{X(0)}{I}
\end{array}
\]
in \veeh, which means exactly that $\xi_{m,n}$ is a map $X(m+n) \go
X(m)\otimes X(n)$ in \und{\emm}, and similarly $\xi_0$ is a map $X(0) \go I$
in \und{\emm}. So the following makes sense:
\begin{defn}	\label{defn:enr-hty-alg}
Let \veeh, $P$ and \emm\ be as above. A \emph{homotopy $P$-algebra in \emm}
is a \veeh-enriched colax (symmetric) monoidal functor
\[
\pr{X}{\xi}: \fmon{P} \go \emm
\]
in which the components $\xi_{m,n}$, $\xi_0$ of $\xi$ are equivalences ($m, n
\in \nat$). 
\end{defn}
A \emph{map of homotopy $P$-algebras} is a \veeh-enriched monoidal
transformation, and we thus obtain a category \homset{\HtyAlg}{P}{\emm}.

\subsubsection*{Examples}

\begin{enumerate}
\setcounter{enumi}{\value{bean}}
\item When $\veeh=\Set$, this reduces to the ordinary, non-enriched
definition of homotopy algebras. 

\item Suppose that the only equivalences in \emm\ are the isomorphisms. Then
by an enriched version of Theorem~\ref{thm:alt-alg}, a homotopy $P$-algebra
in \emm\ is essentially just a $P$-algebra in \emm. Compare
Example~\ref{sec:brief}\bref{eg:hty-alg:isos}. 

\item Take any \veeh, $P$ and \emm\ as in
Definition~\ref{defn:enr-hty-alg}. Then, just as in
Example~\ref{sec:brief}\bref{eg:hty-alg:real-to-hty}, any genuine $P$-algebra
is a homotopy $P$-algebra.

\item Referring back to Examples~\ref{sec:enr-ops}\bref{eg:enr-op:Act-G}
and~\ref{sec:brief}\bref{eg:hty-alg:Act-G}, if $\veeh=\emm=\Top$ and $G$ is a
topological monoid then a homotopy \Act{G}-algebra in \emm\ gives rise to a
\emph{strict} continuous action of $G$ on a space. The same goes for the
other \veeh's, \emm's and $G$'s in
Examples~\ref{sec:enr-ops}\bref{eg:enr-op:Act-G},~\bref{eg:enr-op:Act-G-alg}. 

\item
In~\ref{sec:enr-ops}\bref{eg:enr-op:GrLie} we defined a \GrAb-enriched operad
\GrLie, and observed that a \GrLie-algebra in \ChCx\ is a differential graded
Lie algebra. A homotopy \GrLie-algebra in \ChCx\ might therefore be called a
`homotopy d.g.\ Lie algebra'. It is natural to want to compare this
definition with that of $L_\infty$-algebras (also known as strong homotopy
Lie algebras---see \cite{KSV} and \cite{LM}); however, I have not made such a
comparison. We come back to homotopy Lie algebras in~\ref{subsec:GrLie}. 

\item
Similarly, a homotopy \Ger-algebra in \ChCx\ is some kind of `homotopy
Gerstenhaber algebra'. \Ger\ is defined
in~\ref{sec:enr-ops}\bref{eg:enr-op:Ger}, and homotopy Gerstenhaber algebras
are returned to in~\ref{subsec:Ger}.

\end{enumerate}

\chapter{Homotopy Monoids and Semigroups}
\label{ch:hty-mons}

Just about the simplest algebraic theory is the theory of monoids, and just
about the simplest operad is \Mon. In this chapter we look at homotopy
monoids---that is, homotopy \Mon-algebras---and at homotopy commutative
monoids, homotopy semigroups and homotopy commutative semigroups. (Recall
that a monoid is by definition a semigroup with unit.) Despite the simplicity
of monoids, homotopy monoids and homotopy semigroups provide some of the most
interesting and important examples of homotopy algebras.

The first section~(\ref{sec:Gamma}) is devoted to showing that homotopy
topological commutative monoids are exactly the same as special
$\Gamma$-spaces, and various similar results. For a topologist of a certain
kind, this should help to put the ideas of this paper back onto home
turf. But I have tried to write this paper for both topologists and category
theorists; and in terms of category theory the result is perhaps not all that
interesting. I will now sketch the ideas of Section~\ref{sec:Gamma}, and if
the reader judges them not worthy of further attention then she can
ignore~\ref{sec:Gamma} altogether, without causing a problem in understanding
later sections.

So, the basic idea of~\ref{sec:Gamma} runs as follows. Let $\Phi$ be the
skeletal category of finite
sets~(\ref{sec:prelims-mon-cats}\bref{eg:mon-cat:Phi}). It just so happens
that there is a category $\Gamma$ such that for any category \emm\ with
finite products,
\[
\mbox{colax monoidal functors }
\triple{\Phi}{+}{0} \go \triple{\emm}{\times}{1}
\]
correspond one-to-one with functors $\Gamma^\op \go \emm$. This is merely a
representability result and is not surprising or particularly interesting
from a categorical viewpoint (although the fact that $\Gamma$ has a simple
direct description is not so obvious). A homotopy \CMon-algebra in \emm\ is a
special kind of colax monoidal functor from \triple{\Phi}{+}{0} to
\triple{\emm}{\times}{1}, supposing now that \emm\ is equipped with a class
of equivalences; in other words, it is a special kind of functor $\Gamma^\op
\go \emm$. Such special functors are called `special $\Gamma$-spaces', and
this formulation of the notion of homotopy topological commutative monoid has
been used since around 1970 (see \cite{SegCCT}, \cite{Ad}, \cite{And}). 

The second section~(\ref{sec:loops}) sets out our first major example of a
homotopy algebra: any loop space is a homotopy topological monoid. In
Chapter~\ref{ch:other} we will also exhibit iterated loop spaces as
homotopy-algebraic structures, but this is left alone for now. (See
page~\pageref{page:defn-loop} for the definition of loop space.)

The last three sections (\ref{sec:hty-mon-cats}, \ref{sec:A-spaces},
\ref{sec:A-algs}) each provide a comparison between other notions of weakened
or up-to-homotopy algebraic structure and the present definition of homotopy
algebra. Respectively:
\begin{itemize}
\item any homotopy monoid in \Cat\ gives rise to a monoidal category
\item any homotopy semigroup in \Topstar\ gives rise to an $A_4$-space
\item any homotopy semigroup in \ChCx\ gives rise to an $A_4$-algebra
\end{itemize}
In the last two, I conjecture that `$A_4$' can be replaced by
`$A_\infty$'. In all three, I would like to understand how or whether we
might also pass in the opposite direction ($A_\infty$-algebras giving rise to
homotopy semigroups in \ChCx, for instance), but at present I do not. In a
later chapter~(\ref{sec:inside-Cat}) we describe a converse process for the
first case, \Cat, but there are still many unanswered questions there.

The strategy for these last three sections is to do the hard work only
once. Having done (or rather, asserted that we could do) many tedious
calculations to show that a homotopy monoid in \Cat\ gives rise to a monoidal
category, we can observe that the proof is repeatable in any monoidal
2-category; and from this the results on $A_4$-spaces and $A_4$-algebras
follow quite quickly. It is perhaps only the limitations of
higher-dimensional category theory today which prevent the $A_4$'s from
becoming $A_\infty$'s.

For those unfamiliar with $A_n$-spaces and $A_n$-algebras, the definitions
were made by Stasheff in his 1963 papers~\cite{HAHI}
and~\cite{HAHII}. $A_n$-spaces are to be thought of as topological monoids up
to higher homotopy, and $A_n$-algebras as differential graded algebras up to
higher homotopy; the higher $n$ is, the higher the levels of homotopy go.

\section{$\Gamma$-Objects}
\label{sec:Gamma}

This section is divided into three. The first part concerns homotopy
commutative monoids and special $\Gamma$-objects, and the last part is the
non-symmetric version of this, on homotopy monoids and special simplicial
objects. In the middle is an `aside' outlining a result of category theory
which makes the proofs in the first and last parts much easier. Fans of
abstract methods might like to read the aside first, but others can safely
ignore it.

Those who do not already know what a $\Gamma$-space is are warned (as above)
that they might not find this section very interesting, and might prefer to
jump to~\ref{sec:loops}.

\subsection*{The symmetric case: $\Gamma$-objects}

Let $\Gamma$ be the category defined in \cite{SegCCT}. It is most easily
described by saying that $\Gamma^{\op}$ is (a skeleton of) the category of
finite based sets: thus its objects are $\upr{n} = \{0, 1, \ldots, n \}$ for
$n\geq 0$, and a map $\upr{m} \go \upr{n}$ in $\Gamma^{\op}$ is a function
$g: \upr{m} \go \upr{n}$ such that $g(0)=0$.

Given categories \cat{L} and \emm, we write \ftrcat{\cat{L}}{\emm} for the
usual functor category. If \triple{\cat{L}}{\otimes}{I}
and \triple{\emm}{\otimes}{I} are symmetric monoidal categories then 
\[
\homset{\fcat{SColax}}{\triple{\cat{L}}{\otimes}{I}}{\triple{\emm}{\otimes}{I}}
\]
denotes the category of colax symmetric monoidal functors from
\triple{\cat{L}}{\otimes}{I} to \triple{\emm}{\otimes}{I} and monoidal
transformations. 

In this section, the objects of $\Phi$ will be written as $\lwr{0}, \lwr{1},
\ldots$ rather than $0,1,\ldots$, for clarity. Thus \lwr{n} and \upr{n-1} are
both $n$-element sets. 

\begin{propn}	\label{propn:Gamma-obj}
Let \emm\ be a category with finite products. Then there is an isomorphism of
categories
\[
\homset{\fcat{SColax}}{\triple{\Phi}{+}{\lwr{0}}}{\triple{\emm}{\times}{1}}
\iso
\ftrcat{\Gamma^{\op}}{\emm}.
\]
\end{propn}

\begin{proof}
This is a direct corollary of the general category-theoretic
Proposition~\ref{propn:Kleisli-sym}, under `Aside' below. Alternatively, a
direct argument can be used, as follows.

Given a colax symmetric monoidal functor
\[
\pr{X}{\xi}: \triple{\Phi}{+}{\lwr{0}} \go \triple{\emm}{\times}{1},
\]
we must define a functor $Y: \Gamma^{\op} \go \emm$. On objects, take
$Y\upr{n} = X(\lwr{n})$. To define $Y$ on morphisms, first note the
following:
\begin{itemize}
\item
for each $m\geq 0$, there is a map $\eta_m: \lwr{m} \go \lwr{1+m}$ in
$\Phi$ given by $\eta_m(i) = 1+i$
\item
for each map $g: \upr{m} \go \upr{n}$ in $\Gamma^\op$, there is a
corresponding map $g: \lwr{1+m} \go \lwr{1+n}$ in $\Phi$
\end{itemize}
Now if $g: \upr{m} \go \upr{n}$ is a map in $\Gamma^{\op}$, let $Y(g)$
be the composite
\[
X(\lwr{m}) 	\goby{X(\eta_m)} 
X(\lwr{1+m}) 	\goby{X(g)}
X(\lwr{1+n})	\goby{\xi^2_{1,n}}
X(\lwr{n})
\]
where $\xi^2_{1,n}$ is the second component of $\xi_{1,n}$. This defines a
functor $Y: \Gamma^\op \go \emm$. 

Conversely, take a functor $Y: \Gamma^{\op} \go \emm$. Define a functor $X:
\Phi \go \emm$ by $X(\lwr{n}) = Y\upr{n}$ on objects; if $f: \lwr{m} \go
\lwr{n}$ is a morphism in $\Phi$ then define
\[
X(f) = (Y\upr{m} \goby{Y\upr{g}} Y\upr{n})
\]
where the map $\upr{g}: \upr{m} \go \upr{n}$ in $\Gamma^{\op}$ is given by
\[
g(i) =
\left\{
\begin{array}{ll}
0		&\textrm{if }i=0		\\
1+f(i-1)	&\textrm{if }1\leq i\leq m.	
\end{array}
\right.
\]
To define $\xi$, obviously $\xi_0$ is the unique map $X(\lwr{0}) \go 1$. For
the $\xi_{m,n}$'s, first define maps
\[
\upr{m} \ogby{\pi^1_{m,n}} \upr{m+n} \goby{\pi^2_{m,n}} \upr{n}
\]
in $\Gamma^{\op}$ (for $m,n \geq 0$) by
\begin{eqnarray*}
\pi^1_{m,n}(i)	&=&
\left\{
\begin{array}{ll}
i	&\textrm{if }0\leq i\leq m	\\
0	&\textrm{if }m+1\leq i\leq m+n,	
\end{array}
\right.						\\
\pi^2_{m,n}(i)	&=&
\left\{
\begin{array}{ll}
0	&\textrm{if }0\leq i\leq m	\\
i-m	&\textrm{if }m+1\leq i\leq m+n.	
\end{array}
\right.	
\end{eqnarray*}
Then define
\[
\xi_{m,n}: X(\lwr{m+n}) \go X(\lwr{m}) \times X(\lwr{n})
\]
to be
\[
\pr{Y(\pi^1_{m,n})}{Y(\pi^2_{m,n})}: Y\upr{m+n} \go Y\upr{m} \times Y\upr{n}.
\]

After performing all the checks we see that the two processes are mutually
inverse, and that they can be extended to apply to transformations too. Thus
we obtain the required isomorphism of categories. \done
\end{proof}

We are really only interested in those colax symmetric monoidal functors
\[
\pr{X}{\xi}: \triple{\Phi}{+}{\lwr{0}} \go \triple{\emm}{\times}{1}
\]
in which the components of $\xi$ are equivalences---i.e.\ the homotopy
commutative monoids in \triple{\emm}{\times}{1}. The following result says
what the corresponding condition is on functors $\Gamma^\op \go \emm$. In its
statement, the maps $\pi^{1}_{m,n}$ and $\pi^{2}_{m,n}$ are as in the proof
of Proposition~\ref{propn:Gamma-obj}, and if $0\leq j < n$ then the map
\[
\rho_j^n: \upr{n} \go \upr{1}
\]
is given by
\[
\rho_j^n (i) = 
\left\{
\begin{array}{ll}
1	&\mbox{if }i=j	\\
0	&\mbox{otherwise.}
\end{array}
\right.
\]
Recall from~\ref{sec:prelims-mon-cats}\bref{eg:mon-cat:cart} that a
`cartesian monoidal category' is one in which the monoidal structure is given
by cartesian product and terminal object. 
\begin{propn}	\label{propn:Gamma-special}
Let \triple{\emm}{\times}{1} be a cartesian monoidal category
with equivalences. Let
\[
\pr{X}{\xi}: \triple{\Phi}{+}{\lwr{0}} \go \triple{\emm}{\times}{1}
\]
be a colax symmetric monoidal functor, and let
\[
Y: \Gamma^\op \go \emm
\]
be the functor corresponding to \pr{X}{\xi} under
Proposition~\ref{propn:Gamma-obj}. The following conditions are equivalent:
\begin{enumerate}
\item 	\label{part:Gamma-hty-comm-mon}
\pr{X}{\xi} is a homotopy commutative monoid
\item 	\label{part:Gamma-mn}
for each $m,n\geq 0$, the map
\[
\pr{Y(\pi^{1}_{m,n})}{Y(\pi^{2}_{m,n})}:
Y\upr{m+n} \go Y\upr{m} \times Y\upr{n}
\]
is an equivalence, and so is the unique map $Y\upr{0} \go 1$
\item  	\label{part:Gamma-n}
for each $n\geq 0$, the map
\[
\bftuple{Y(\rho_0^n)}{Y(\rho_{n-1}^n)}: Y\upr{n} \go Y\upr{1}^n
\]
is an equivalence.
\end{enumerate}
\end{propn}
\begin{proof}
\bref{part:Gamma-hty-comm-mon} $\Leftrightarrow$ \bref{part:Gamma-mn} is
immediate from the second half of the proof of
Proposition~\ref{propn:Gamma-obj}. An easy induction, again using this half
of the proof, shows that
\[
\bftuple{Y(\rho_0^n)}{Y(\rho_{n-1}^n)}: Y\upr{n} \go Y\upr{1}^n
\]
equals
\[
\xi^{(n)}: X(\lwr{n}) \go X(\lwr{1})^n:
\]
so by the comments in~\ref{sec:brief}\bref{eg:hty-alg:Obj}, we have
\bref{part:Gamma-hty-comm-mon} $\Leftrightarrow$ \bref{part:Gamma-n}.  \done
\end{proof}

A functor $Y: \Gamma^\op \go \emm$ satisfying the equivalent
conditions~\bref{part:Gamma-mn} and~\bref{part:Gamma-n} will be called a
\emph{special $\Gamma$-object} in \emm, and the category of such, with
natural transformations as morphisms, will be written 
\[
\homset{\fcat{Special}}{\Gamma^{\op}}{\emm}.
\]
(A \emph{$\Gamma$-object} in \emm\ is any old functor from $\Gamma^\op$ to
\emm.)  We then have:
\begin{cor}
Let \triple{\emm}{\times}{1} be a cartesian monoidal category with
equivalences. Then there is an isomorphism of categories
\[
\homset{\HtyAlg}{\CMon}{\emm} \iso 
\homset{\fcat{Special}}{\Gamma^{\op}}{\emm}.
\]
\done
\end{cor}
In the original paper \cite{SegCCT}, the cases $\emm=\Top$ and $\emm=\Cat$
are considered. (I should re-iterate that what are called $\Gamma$-spaces and
$\Gamma$-categories in \cite{SegCCT} are called special $\Gamma$-spaces and
special $\Gamma$-categories here.) In these cases we have:
\begin{cor}
\begin{enumerate}
\item homotopy topological commutative monoids are the same as special
$\Gamma$-spaces 
\item homotopy symmetric monoidal categories are the same as special
$\Gamma$-cat\-e\-gories. 
\end{enumerate}
\done
\end{cor}

\subsection*{Aside: a general result}
\small

Proposition~\ref{propn:Gamma-obj} is in fact a special case of the following
category-theoretic result. Anyone who understands the statement will probably
have little trouble in supplying a proof; I hope to write it up
separately. For the definition of Kleisli category see~\cite[VI.5]{CWM}.

\begin{propn}	\label{propn:Kleisli-sym}
Let \triple{\cat{L}}{\otimes}{I} be a symmetric monoidal category with a
terminal object $1$, and let $\cat{L}_{1\otimes\dashbk}$ be the Kleisli
category for the monad $1\otimes\dashbk$ on \cat{L}. For any category \emm\
with finite products, there is an isomorphism of categories
\[
\homset{\fcat{SColax}}{\triple{\cat{L}}{\otimes}{I}}{\triple{\emm}{\times}{1}}
\iso
\ftrcat{\cat{L}_{1\otimes\dashbk}}{\emm}.
\]
\done
\end{propn}
Taking $\triple{\cat{L}}{\otimes}{I} = \triple{\Phi}{+}{\lwr{0}}$ we obtain
Proposition~\ref{propn:Gamma-obj} immediately: for
$\cat{L}_{1\otimes\dashbk}$ is the skeletal category of finite based sets,
$\Gamma^{\op}$.

There is a non-symmetric version of Proposition~\ref{propn:Kleisli-sym} too,
as follows. Here \fcat{Colax} denotes the category of colax monoidal
functors.
\begin{propn}	\label{propn:Kleisli-non-sym}
Let \triple{\cat{L}}{\otimes}{I} be a monoidal category with a terminal
object $1$, and let $\cat{L}_{1\otimes\dashbk\otimes 1}$ be the Kleisli
category for the monad $1\otimes\dashbk\otimes 1$ on \cat{L}. For any
category \emm\ with finite products, there is an isomorphism of categories
\[
\homset{\fcat{Colax}}{\triple{\cat{L}}{\otimes}{I}}{\triple{\emm}{\times}{1}}
\iso
\ftrcat{\cat{L}_{1\otimes\dashbk\otimes 1}}{\emm}.
\]
\done
\end{propn}
To apply this to the case of (non-commutative) monoids, take
$\triple{\cat{L}}{\otimes}{I} = \triple{\Delta}{+}{\lwr{0}}$. Then the
Kleisli category $\Delta_{1+\dashbk +1}$ is a skeleton of the category of
\emph{finite strict intervals}, that is, of finite totally ordered sets with
distinct greatest and least elements. But it is well-known that this category
is isomorphic to $(\Delta^+)^{\op}$ (defined below), and this gives us
Proposition~\ref{propn:Delta-obj}. (The isomorphism
\[
(\Delta^+)^{\op} \iso \textrm{(finite strict intervals)}
\]
can be written in either direction as \homset{\Hom}{\dashbk}{\upr{1}}, i.e.\
\upr{1} is a `schizophrenic object'. A few more details can be found
in~\cite{Joy}.) 

We have extracted maximum use from Propositions~\ref{propn:Kleisli-sym}
and~\ref{propn:Kleisli-non-sym}, in the following sense: we can't
apply~\ref{propn:Kleisli-sym} with $\cat{L} = \fmon{P}$ for any $P$ other
than \CMon, since \cat{L} is required to have a terminal object, and
similarly \Mon\ for Proposition~\ref{propn:Kleisli-non-sym}. On the other
hand, some general categorical arguments guarantee that whether or not
\cat{L} has a terminal object, there exists a category \cat{L'} with the
property that $\cat{L}_{1\otimes\dashbk}$ (or
$\cat{L}_{1\otimes\dashbk\otimes 1}$) has in
Proposition~\ref{propn:Kleisli-sym} (or~\ref{propn:Kleisli-non-sym}). (I hope
to explain this properly elsewhere.)  So, for any operad $P$, there is
\emph{some} category $(\fmon{P})'$ playing the role that $\Gamma^{\op}$ did
for \CMon\ and that $(\Delta^+)^{\op}$ did for \Mon. However, $(\fmon{P})'$
might not always be easy to describe.
\normalsize

\subsection*{The non-symmetric case: simplicial objects}

Returning to the main exposition, let 
\label{page:Delta-plus}%
$\Delta^+$ be the ``topologists' simplicial category'', that is, the category
whose objects are $\upr{n} = \{0,1,\ldots,n\}$ for $n\geq 0$, and whose
morphisms are order-preserving functions. So $\Delta^+$ is obtained from
$\Delta$ by removing the object $0$ and renaming the objects which
remain. For emphasis I will now write the objects of $\Delta$ as $\lwr{0},
\lwr{1}, \ldots$ rather than $0,1,\ldots$; hence \lwr{n} and \upr{n-1} are
both $n$-element ordered sets.

By definition, a simplicial object in a category \emm\ is a functor from
$(\Delta^+)^\op$ to \emm. The next result gives an alternative definition of
simplicial object.

\begin{propn}	\label{propn:Delta-obj}
Let \emm\ be a category with finite products. Then there is an isomorphism of
categories
\[
\homset{\fcat{Colax}}{\triple{\Delta}{+}{\lwr{0}}}{\triple{\emm}{\times}{1}}
\iso
\ftrcat{(\Delta^+)^{\op}}{\emm}.
\]
\end{propn}
\paragraph*{Remark} \fcat{Colax} denotes the category of colax monoidal
functors and monoidal transformations. 
\begin{sketchpf}
A colax monoidal functor
\[
\pr{X}{\xi}: \triple{\Delta}{+}{\lwr{0}} \go \triple{\emm}{\times}{1}
\]
consists of a functor $X: \Delta \go \emm$ together with maps
\[
X(\lwr{m}) \ogby{\xi^1_{m,n}} X(\lwr{m+n}) \goby{\xi^2_{m,n}} X(\lwr{n})
\]
for each $m,n\geq 0$, satisfying certain axioms. These axioms imply that all
the $\xi^i_{m,n}$'s can be built up as composites of maps
\[
\xi^1_{k,1}: X(\lwr{k+1}) \go X(\lwr{k}),
\diagspace
\xi^2_{1,k}: X(\lwr{1+k}) \go X(\lwr{k}).
\]
Hence a colax monoidal functor from \triple{\Delta}{+}{\lwr{0}} to
\triple{\emm}{\times}{1} is a functor $X: \Delta \go \emm$ together with a
pair of maps
\[
X(\lwr{k+1}) \pile{\rTo\\ \rTo} X(\lwr{k})
\]
for each $k\geq 0$, satisfying certain axioms. This data can be depicted as
\begin{diagram}
X(\lwr{0})	&
\pile{\lGet \\ \rTo \\ \lGet}	&
X(\lwr{1})	&
\pile{\lGet\\ \rTo\\ \lTo\\ \rTo\\ \lGet}	&
X(\lwr{2})	&
\cdots		\\
\end{diagram}
where the solid lines are the image under $X$ of the face and degeneracy maps
in $\Delta$, and the dotted lines are $\xi^1_{k,1}$ and $\xi^2_{1,k}$. But
this diagram looks just like the usual picture of a simplicial object: so we
define a functor $Y: (\Delta^+)^{\op} \go \emm$ by 
\begin{eqnarray*}
Y\upr{n}	&=&	X(\lwr{n})	\\
Y(\sigma_i^n)	&=&	X(\delta_i^{n-1})	\\
Y(\delta_i^n)	&=&
\left\{
\begin{array}{ll}
\xi^1_{n-1,1}		&\textrm{if }i=0		\\
X(\sigma_{i-1}^{n-1})	&\textrm{if }1\leq i\leq n-1	\\
\xi^2_{1,n-1}		&\textrm{if }i=n		\\	
\end{array}
\right.
\end{eqnarray*}
where $\sigma_i^n$ are the degeneracy maps in $\Delta$, and $\delta_i^n$ the
face maps, with notation as in \cite[VII.5]{CWM}.  \done
\end{sketchpf}

I want to emphasize the difference between the categories $\Delta$ and
$\Delta^+$. For our purposes, the only interesting relation between the two
is that stated in the Proposition: a colax monoidal functor from $\Delta$ to
a cartesian monoidal category \emm\ is the same as a functor from
$(\Delta^+)^\op$ to \emm. The fact that $\Delta^+$ is $\Delta$ with one
object removed can be regarded as nothing more than a distracting
coincidence.

We are mostly interested in a certain subset of the colax monoidal
functors from $\Delta$ to \emm: the homotopy monoids. The next result
shows what property of a functor $\Delta^+ \go \emm$ corresponds to the colax
monoidal functor being a homotopy monoid.

To state it we need some more notation. For $m,n\geq 0$, define maps 
\[
\upr{m} \goby{\alpha^1_{m,n}} \upr{m+n} \ogby{\alpha^2_{m,n}} \upr{n}
\]
in $\Delta^+$ by $\alpha^1_{m,n}(i) = i$ and $\alpha^2_{m,n}(i) = m+i$. For
$0\leq j<n$, define a map
\[
\beta^n_j: \upr{1} \go \upr{n}
\]
by $\beta^n_j (i) = i+j$ for $i=0,1$.

\begin{propn}	\label{propn:Delta-special}
Let \emm\ be a cartesian monoidal category with equivalences. Let
\[
\pr{X}{\xi}: \triple{\Delta}{+}{\lwr{0}} \go \triple{\emm}{\times}{1}
\]
be a colax monoidal functor, and let
\[
Y: (\Delta^+)^{\op} \go \emm
\]
be the functor corresponding to \pr{X}{\xi} under
Proposition~\ref{propn:Delta-obj}. The following conditions are equivalent:
\begin{enumerate}
\item 	\label{part:special:hty-mon}
\pr{X}{\xi} is a homotopy monoid
\item 	\label{part:special:colax}
for each $m,n\geq 0$, the map 
\[
\pr{Y(\alpha^1_{m,n})}{Y(\alpha^2_{m,n})}:
Y\upr{m+n} \go Y\upr{m} \times Y\upr{n}
\]
is an equivalence, and so is the unique map $Y\upr{0} \go 1$
\item 	\label{part:special:powers}
for each $n\geq 0$, the map
\[
\bftuple{Y(\beta^n_0)}{Y(\beta^n_{n-1})}:
Y\upr{n} \go Y\upr{1}^n
\]
is an equivalence.
\end{enumerate}
\end{propn}

\begin{proof}
\bref{part:special:hty-mon} $\Leftrightarrow$ \bref{part:special:colax}: From
the definition of $Y$ in the proof of Proposition~\ref{propn:Delta-obj}, it
is easy to show that $Y(\alpha^1_{m,n}) = \xi^1_{m,n}$, and similarly
$\alpha^2_{m,n}$. Hence 
\[
\pr{Y(\alpha^1_{m,n})}{Y(\alpha^2_{m,n})} = \xi_{m,n}.
\]
Also, $\xi_0$ is the unique map from $X(\lwr{0}) = Y\upr{0}$ to
$1$.

\bref{part:special:hty-mon} $\Leftrightarrow$ \bref{part:special:powers}:
Similarly, an easy induction on $n$ shows that
\[
\bftuple{Y(\beta^n_0)}{Y(\beta^n_{n-1})} = \xi^{(n)}: X(\lwr{n}) \go
X(\lwr{1})^n, 
\]
where $\xi^{(n)}$ is the map defined
in~\ref{sec:brief}\bref{eg:hty-alg:Obj}. But by the comments
in~\ref{sec:brief}\bref{eg:hty-alg:Obj}, \pr{X}{\xi} is a homotopy monoid if
and only if each $\xi^{(n)}$ is an equivalence.
\done
\end{proof}

A functor $Y: (\Delta^+)^{\op} \go \emm$ satisfying the equivalent
conditions~\bref{part:special:colax} and~\bref{part:special:powers} will be
called a \emph{special simplicial object} in \emm. (When $\emm=\Top$ the
usual name is `special $\Delta$-space'; the `$\Delta$' used in this name is
our $\Delta^+$.)  Write 
\[
\homset{\fcat{Special}}{(\Delta^+)^{\op}}{\emm}
\]
for the category of special simplicial objects in \emm\ and natural
transformations. Then we immediately have:
\begin{cor}
Let \triple{\emm}{\times}{1} be a cartesian monoidal category with
equivalences. There is an isomorphism of categories
\[
\homset{\HtyAlg}{\Mon}{\emm} 
\iso 
\homset{\fcat{Special}}{(\Delta^+)^{\op}}{\emm}.
\]
\done
\end{cor}

In particular, this holds for $\emm=\Cat$ and $\emm=\Top$, in which contexts
special simplicial objects are best known (see~\cite{And}, for example). 

We have spent this section showing that homotopy commutative monoids are the
same as special $\Gamma$-objects, and homotopy monoids the same as special
simplicial objects. \emph{But this is only true when \emm\ is cartesian.}  If
the tensor $\otimes$ in a symmetric monoidal category \emm\ is not the
cartesian (categorical) product, or if the unit $I$ is not the terminal
object, then the two categories
\[
\ftrcat{\Gamma^\op}{\emm},
\diagspace
\homset{\fcat{SColax}}{\triple{\Phi}{+}{\lwr{0}}}{\triple{\emm}{\otimes}{I}}
\]
are in general nowhere near equivalent. Similarly, the two categories
\[
\ftrcat{(\Delta^+)^\op}{\emm},
\diagspace
\homset{\fcat{Colax}}{\triple{\Delta}{+}{\lwr{0}}}{\triple{\emm}{\otimes}{I}}
\]
are quite different. Moreover, the definition of `special' (Propositions
\ref{propn:Gamma-special}\bref{part:Gamma-mn},\bref{part:Gamma-n} and
\ref{propn:Delta-special}\bref{part:special:colax},\bref{part:special:powers})
depend on the product in \emm\ being the cartesian product, and there is no
obvious way to express it for an arbitrary monoidal structure on \emm. 

In Section~\ref{sec:A-algs}, on $A_\infty$-algebras, we will see our first
significant example of homotopy algebras in a non-cartesian monoidal
category.

\section{Loop Spaces}
\label{sec:loops}

Our main result is that any loop space is a homotopy monoid. More exactly:
\begin{thm}	\label{thm:loop}
There is a functor
\[
\Omega: \Topstar \go \homset{\HtyAlg}{\Mon}{\Top}
\]
which sends a based space $B$ to a homotopy monoid $\Omega B$, with 
\[
(\Omega B)(1) \iso \homset{\Topstar}{S^1}{B}.
\]
\end{thm}

The heart of the proof is that the circle $S^1$ is a `homotopy comonoid':
\begin{lemma}	\label{lemma:S1-comonoid}
There is a homotopy monoid
\[
\pr{W}{\omega}: \triple{\Delta}{+}{0} \go 
\triple{\Top_*^\op}{\wej}{1}
\]
with $W(1)=S^1$.
\end{lemma}
\begin{proof}
For each $n\geq 0$, let $(\Delta^n/\sim)$ denote the standard $n$-simplex
$\Delta^n$ with its $n+1$ vertices collapsed to a single point, and this
point declared the basepoint. Informally, \pr{W}{\omega} can be described as
follows:
\begin{itemize}
\item $W(n)=\Delta^n/\sim$, e.g.\ $W(0)$ is a single point, $W(1)=S^1$, and
$W(2)$ looks like \Wtwo
\item $W$ is defined on morphisms by the standard face and degeneracy maps of
simplices
\item $\omega$ is defined by face maps, e.g.\ 
\[
\omega_{1,1}: W(1) \wej W(1) \go W(2)
\]
is the evident inclusion
\[
\Wedgepic \rIncl \Wtwo,
\]
which is a homotopy equivalence. 
\end{itemize}

Formally, it's easiest to employ the description of homotopy monoids given in
Proposition~\ref{propn:Delta-special}. So, first consider the usual functor
from $\Delta^+$ to \Top, mapping \upr{n} to $\Delta^n$ and defined on
morphisms by face and degeneracy maps. Since all the face and degeneracy maps
take vertices to vertices, this functor induces another functor
\[
Y: \Delta^+ \go \Topstar
\]
with $Y\upr{n} = \Delta^n/\sim$. We may also view $Y$ as a functor
\[
Y: (\Delta^+)^\op \go \Top_*^\op;
\]
$\wej$ and $1$ are respectively coproduct and initial object in \Topstar, so
they are product and terminal object in $\Top_*^\op$. So by
Proposition~\ref{propn:Delta-obj}, $Y$ corresponds to a colax monoidal
functor
\[
\pr{W}{\omega}: \triple{\Delta}{+}{0} \go \triple{\Top_*^\op}{\wej}{1}
\]
with $W(1) = Y\upr{1} = S^1$. To see that \pr{W}{\omega} is a homotopy monoid
we must first of all check that the unique map $Y\upr{0} \go 1$ is a homotopy
equivalence, and then check that the map
\[
(\Delta^m/\sim) \wej (\Delta^n/\sim) \go \Delta^{m+n}/\sim 
\]
induced by the two maps
\[
\upr{m} \ogby{\alpha^1_{m,n}} \upr{m+n} \goby{\alpha^2_{m,n}} \upr{n}
\]
in $\Delta^+$ is also a homotopy equivalence (see
Proposition~\ref{propn:Delta-special}\bref{part:special:colax}). The first
check is trivial, and for the second it is easy to construct a homotopy
inverse. (From the conceptual angle, note however that there is no canonical
choice of a homotopy inverse: see the picture of $\omega_{1,1}$ above, for
instance.)  \done
\end{proof}

To prove Theorem~\ref{thm:loop}, observe that if $B$ is a based space then 
\[
\homset{\Topstar}{\dashbk}{B}:
\triple{\Top_*^{\op}}{\wej}{1} \go \triple{\Top}{\times}{1}
\]
is a monoidal functor, since $\wej$ is the coproduct in $\Topstar$ and $1$ is
the initial object. Observe also that \homset{\Topstar}{\dashbk}{B}
preserves homotopy equivalences. Hence the composite
\[
\triple{\Delta}{+}{0} \goby{\pr{W}{\omega}}
\triple{\Top_*^{\op}}{\wej}{1} \goby{\homset{\Topstar}{\dashbk}{B}}
\triple{\Top}{\times}{1}
\]
is a homotopy monoid, which we call $\Omega B$. Moreover, any map $f: B\go
B'$ induces a monoidal transformation
\[
\homset{\Topstar}{\dashbk}{B} \go \homset{\Topstar}{\dashbk}{B'}
\]
and therefore a map $\Omega B \go \Omega B'$; this makes $\Omega$ into a
functor 
\[
\Topstar \go \homset{\HtyAlg}{\Mon}{\Top}.
\]
Finally, 
\[
(\Omega B)(1) = \homset{\Topstar}{W(1)}{B} = \homset{\Topstar}{S^1}{B},
\]
as required. 

(In~\ref{sec:change-egs}\bref{eg:change:hom} we'll see a neater way to
express the proof of the Theorem: \homset{\Topstar}{\dashbk}{B} is a `change
of environment'. The proof that a loop space is a special simplicial object
was apparently first found by Segal.)

Let us pause for a moment to examine the homotopy monoid structure with which
we have just equipped loop spaces, and in particular at how the composition
of two loops is handled. 

Fix a space $B$, and write \pr{X}{\xi} for the homotopy monoid $\Omega B$.
We have
\[
X(2) = \homset{\Topstar}{\Wtwo}{B},
\]
and the pieces of \pr{X}{\xi} relevant to binary composition are the maps
\begin{diagram}[size=2.5em]
	&				&X(2)	&		&	\\
	&\ldTo<{\xi_{1,1}}>{\eqv}	&	&\rdTo>{X(!)}	&	\\
X(1)^2	&				&	&		&X(1).	\\
\end{diagram}
This diagram is
\begin{equation}	\label{eqn:canonical}
\begin{diagram}[size=2.5em]
	&		&\homset{\Topstar}{\Wtwo}{B}	&	&	\\
	&\ldTo<{\eqv}	&				&\rdTo	&	\\
\homset{\Topstar}{\Wedgepic}{B}	= \Loops^2&	&	&	&\Loops,\\
\end{diagram}
\end{equation}
where the map on the left is restriction to the two inner circles and the map
on the right is restriction to the outer circle. All of the data
making up \pr{X}{\xi}, and in particular the maps in~\bref{eqn:canonical},
is constructed canonically from $B$: no arbitrary choices have been made. In
contrast, there is no canonical map
\[
\Loops^2 \go \Loops
\]
defining `composition': although the obvious and customary choice is to use
the map described by the instruction `travel each loop at double speed', this
appears to have no particular advantage or special algebraic status compared
to any other choice. Since the usual formulation of the idea of homotopy
topological monoid, $A_\infty$-spaces, does entail this arbitrary choice of a
composition law for loops, one might regard this as an ideological virtue of
the definition presented here.

In Section~\ref{sec:A-spaces} below we make a down-to-earth comparison
between $A_\infty$-spaces and homotopy topological monoids.

This section closes with three further remarks on loop spaces. Firstly, we
have shown that any loop space is a homotopy topological monoid in the sense
of being a homotopy algebra for the non-symmetric operad \Mon. But the
symmetric operad \Sym\ (\ref{sec:operads}\bref{eg:operad:Sym}) also has the
property that algebras for it (in any symmetric monoidal category, such as
\Top) are monoids, so we might also try to show that any loop space is a
homotopy monoid in the sense of being a homotopy
\label{page:loop-hty-sym}%
\Sym-algebra. I believe this to be true, using a colax symmetric monoidal
functor \pr{W}{\omega} with $W(n) = \Delta^n/\sim$ again, but I am not
sure. It seems that a homotopy \Mon-algebra is not automatically a homotopy
\Sym-algebra; some extra input is required. 

Secondly, recall from~\ref{sec:Gamma} that a special $\Gamma$-space is the
same thing as a homotopy topological commutative monoid. The relationship
between special $\Gamma$-spaces and infinite loop spaces (spectra) is
well-explored (see \cite{SegCCT}, \cite{Ad}, \cite{MT} etc.), and we may say
`any infinite loop space is a homotopy commutative monoid'. So we have a
homotopy algebra structure on $n$-fold loop spaces for $n=1$ and
$n=\infty$. Thirdly, then, how about $1<n< \infty$? The answer is that any
$n$-fold loop space is an `$n$-fold homotopy monoid'---but for an explanation
of that, the reader will have to wait until Section~\ref{sec:iterated}.

\section{Monoidal Categories}
\label{sec:hty-mon-cats}

Another major example of homotopy algebras is that of homotopy \Mon-algebras
in \Cat, which I shall call \emph{homotopy monoidal categories}. This example
was introduced in~\ref{sec:brief}\bref{eg:hty-alg:hty-mon-cat}. We have three
different kinds of structure in front of us: strict monoidal categories,
(non-strict) monoidal categories, and homotopy monoidal categories. In this
section a partial comparison is made between the last two.

The philosophical difference between monoidal categories and homotopy monoidal
categories can be put like this. In an (ordinary) monoidal category such as
\label{page:Ab}%
\triple{\Ab}{\otimes}{\integers}, tensor is an operation with definite and
precise values: that is, if $A$ and $B$ are abelian groups then there is
assigned an abelian group $A\otimes B$, not just defined up to isomorphism,
but with actual, specific, elements. Of course, we do not care precisely what
these elements are; there are various definitions of $\otimes$ which are all
naturally isomorphic, and this is all that matters in practice. So there is
some degree of artificial choice in the definition of tensor, but it is at
least an honest functor. Similar comments apply to the cartesian product
$\times$ of sets. In contrast, homotopy monoidal categories allow a certain
amount of fuzziness: there is not an actual operation `tensor' (as we shall
see), but the substitute for `tensor' is intended to avoid artificial
choices. Applied to \Ab, for instance, the rough idea is to record all
quadruples $\tuplebts{A,B,u,C}$ where $u: A\times B \go C$ is a bilinear map
with the usual universal property, but without choosing any preferred $C$ for
each $A$ and $B$. Thus one could speak of $C$ as `a tensor' of $A$ and $B$,
but one would never speak of `the tensor'.

This idea is well known (if not deeply understood) in higher-dimensional
category theory: see, for instance, \cite{HerCSU}, \cite{Baez},
\cite{BaDoHDA3}, \cite{HMP}, \cite{Joy}. Within category theory it was
perhaps first taken seriously by Makkai in his study of anafunctors (see
\cite{Mak}), which arise implicitly throughout the present work. Homotopy
monoidal categories resemble closely what Makkai called anamonoidal
categories. The idea was also exploited in topology by Segal.

In this section we show how a homotopy monoidal category gives rise to a
monoidal category, and similarly for symmetric monoidal categories. We then
notice that this result generalizes effortlessly to an arbitrary monoidal
2-category (in place of \Cat). This is very trivial, but will enable us later
to read off results about $A_n$-spaces and $A_n$-algebras. 

There is also a converse process in the case of \Cat: any monoidal category
gives rise to a homotopy monoidal category. This is more mysterious and less
well-developed than the process the first way round, and the overall
comparison between monoidal categories and homotopy monoidal categories is
therefore incomplete. We do not describe this converse process here, but
instead leave it until a later section (\ref{sec:inside-Cat}) where it is
done in more generality.

\subsection*{How a homotopy monoidal category give rise to a monoidal
category}

A homotopy monoidal category consists of a functor $C:\Delta \go \Cat$
(previously called $X$) together with equivalences of categories
\[
\begin{array}{rccc}
\xi_{m,n}:	&C(m+n) 	&\go 	&C(m) \times C(n)	\\
\xi_{0}:	&C(0) 		&\go 	&\One
\end{array}
\]
($m,n\geq 0$) fitting together nicely.  We regard $C(1)$ as the `base
category' of the homotopy monoidal category \pr{C}{\xi}, recalling
Example~\ref{sec:brief}\bref{eg:hty-alg:Obj}.

\begin{propn}	\label{propn:hty-mon-cat}
A homotopy monoidal category gives rise to a monoidal category.
\end{propn}
\begin{proof}
Take a homotopy monoidal category \pr{C}{\xi} in \Cat, and construct from it
a monoidal category as follows.

\begin{description}
\item[Underlying category:] $C(1)$.

\item[Tensor:] What we want to define is a functor
\[
\otimes: C(1) \times C(1) \go C(1);
\]
what we actually have are functors
\begin{diagram}[size=2.5em]
		&			&C(2)	&		&	\\
		&\ldTo<{\xi_{1,1}}>{\eqv}&	&\rdTo>{C(!)}	&	\\
C(1)\times C(1)	&			&	&		&C(1)	\\
\end{diagram}
where $!$ is the unique map $2\go 1$ in $\Delta$. So for each $m$ and $n$,
choose (arbitrarily) a pseudo-inverse $\psi_{m,n}$ to $\xi_{m,n}$, and
define $\otimes$ as the composite
\[
C(1) \times C(1) \goby{\psi_{1,1}} C(2) \goby{C(!)} C(1).
\]

\item[Associativity isomorphisms:] The next piece of data we need is a
natural isomorphism between $\otimes\of(\otimes\times 1)$ and
$\otimes\of(1\times\otimes)$. To see why such an isomorphism should exist,
consider what would happen if the $\psi_{m,n}$'s were \emph{genuine}
inverses to the $\xi_{m,n}$'s. Then the $\psi_{m,n}$'s would satisfy the
same coherence and naturality axioms as the $\xi_{m,n}$'s (with the arrows
reversed), and this would guarantee that all sensible diagrams built up out
of $\psi_{m,n}$'s commuted. Hence $\otimes$ would be strictly
associative. As it is, $\psi_{m,n}$ is only inverse to $\xi_{m,n}$ up to
isomorphism, and correspondingly $\otimes$ is associative up to isomorphism.

In practice, let us choose (at random) natural isomorphisms
\[
\eta_{m,n}: 1 \goiso \xi_{m,n}\of\psi_{m,n},
\diagspace
\epsln_{m,n}: \psi_{m,n}\of\xi_{m,n} \goiso 1
\]
for each $m$ and $n$. Then a natural isomorphism
\[
\alpha: \otimes\of(\otimes\times 1) \goiso \otimes\of(1\times\otimes)
\]
can be built up from the $\eta_{m,n}$'s and $\epsln_{m,n}$'s. The exact formula
for $\alpha$ is rather complicated, and only included for the record. Most
readers will therefore want to skip the next paragraph and ignore
Figure~\ref{fig:associator}.

\begin{figure}

\begin{diagram}
C(1)^3			&\rTo^{\psi_{1,1}\times 1}	&
C(2)\times C(1)		&\rTo^{C(!)\times 1}		&C(1)^2	\\
\dTo<{1\times\psi_{1,1}}	&	&
\dTo>{\psi_{2,1}}		&	&\dTo>{\psi_{1,1}}	\\
C(1)\times C(2)		&\rTo^{\psi_{1,2}}		&
C(3)			&\rTo^{C(\sigma_0)}		&C(2)	\\
\dTo<{1\times C(!)}		&	&
\dTo>{C(\sigma_1)}		&	&\dTo>{C(!)}		\\
C(1)^2			&\rTo_{\psi_{1,1}}		&
C(2)			&\rTo_{C(!)}			&C(1)	\\
\end{diagram}

\[
\begin{array}{rl}
&\otimes\of(\otimes\times 1)						\\
=
&C(!) \of \psi_{1,1} \of [C(!)\times 1] \of [\psi_{1,1}\times 1] 	\\
\goby{1*\eta_{2,1}*1}
&C(!) \of \psi_{1,1} \of [C(!)\times 1] \of \xi_{2,1} \of \psi_{2,1}
\of [\psi_{1,1}\times 1] 						\\
=
&C(!) \of \psi_{1,1} \of \xi_{1,1} \of C(\sigma_0) \of \psi_{2,1}
\of [\psi_{1,1}\times 1] 						\\
\goby{1*\epsln_{1,1}*1}
&C(!) \of C(\sigma_0) \of \psi_{2,1} \of [\psi_{1,1}\times 1] 	\\
\goby{1*[1\times \eta_{1,1}]}
&C(!) \of C(\sigma_0) \of \psi_{2,1} \of [\psi_{1,1}\times 1]
\of [1\times \xi_{1,1}] \of [1\times \psi_{1,1}] 			\\
\goby{1*\eta_{1,2}*1}
&C(!) \of C(\sigma_0) \of \psi_{2,1} \of [\psi_{1,1}\times 1]
\of [1\times \xi_{1,1}] \of \xi_{1,2} \of \psi_{1,2}
\of [1\times \psi_{1,1}] 						\\
=
&C(!) \of C(\sigma_0) \of \psi_{2,1} \of [\psi_{1,1}\times 1]
\of [\xi_{1,1}\times 1] \of \xi_{2,1} \of \psi_{1,2}
\of [1\times \psi_{1,1}] 						\\
\goby{1*[\epsln_{1,1}\times 1]*1}
&C(!) \of C(\sigma_0) \of \psi_{2,1} \of \xi_{2,1} \of \psi_{1,2}
\of [1\times \psi_{1,1}] 						\\
\goby{1*{\epsln_{2,1}}*1}
&C(!) \of C(\sigma_0) \of \psi_{1,2} \of [1\times \psi_{1,1}] 	\\
=
&C(!) \of C(\sigma_1) \of \psi_{1,2} \of [1\times \psi_{1,1}] 	\\
\goby{1*\epsln_{1,1}^{-1}*1}
&C(!) \of \psi_{1,1} \of \xi_{1,1} \of C(\sigma_1) \of \psi_{1,2} 
\of [1\times \psi_{1,1}] 						\\
=
&C(!) \of \psi_{1,1} \of [1\times C(!)] \of \xi_{1,2}\of \psi_{1,2} 
\of [1\times \psi_{1,1}] 						\\
\goby{1*\eta_{1,2}^{-1}*1}
&C(!) \of \psi_{1,1} \of [1\times C(!)] \of [1\times \psi_{1,1}] 	\\
=
&\otimes\of (1\times\otimes).
\end{array}
\]

\caption{Formula for the associativity isomorphism $\alpha$}
\label{fig:associator}
\end{figure}

In order to define $\alpha$, consider the diagram at the top of
Figure~\ref{fig:associator}. The composites around the outside are
$\otimes\of(\otimes\times 1)$ and $\otimes\of (1\times\otimes)$, so we must
find natural isomorphisms inside each of the four inner squares. The
bottom-right square, in which $\sigma_0, \sigma_1$ are the two surjections
$3 \go 2$ in $\Delta$, is genuinely commutative. In each of the
other three squares, imagine taking each arrow labelled by a $\psi$,
reversing its direction, and changing the $\psi$ to a $\xi$. The imaginary
square would then be genuinely commutative in each case, which means that the
actual square is commutative up to isomorphism. This is the thought behind
the formula for $\alpha$ given in the rest of
Figure~\ref{fig:associator}. (For the usage of $*$, see
page~\pageref{page:misc-notation}.)

\item[Pentagon:] We must now check that the associativity isomorphism just
defined satisfies the famous pentagon coherence axiom. This asserts the
commutativity of a certain diagram built up from components of $\alpha$, that
is, built up from $\eta_{m,n}$'s and $\epsln_{m,n}$'s. However, this diagram
does \emph{not} commute, which is perhaps unsurprising since $\eta_{m,n}$ and
$\epsln_{m,n}$ were chosen independently.

But all is not lost: for recall the result that if
\[
\begin{array}{ll}
F: A \go B,			&G: B\go A,		\\
\sigma: 1 \goiso G\of F,	&\tau: F\of G \goiso 1
\end{array}
\]
is an equivalence of categories, then $\tau$ can be exchanged for another
natural isomorphism $\tau'$ so that $(F,G,\sigma,\tau')$ is both an
adjunction and an equivalence (see \cite[IV.4.1]{CWM}). So when we chose the
natural isomorphisms $\eta_{m,n}$ and $\epsln_{m,n}$ above, we could have
done it so that $(\psi_{m,n},\xi_{m,n},\eta_{m,n},\epsln_{m,n})$ was an
adjunction. Assume that we did so. Then this being an adjunction says that
certain basic diagrams involving $\eta_{m,n}$ and $\epsln_{m,n}$ commute
(namely, the diagrams for the triangle identities \cite[IV.1(9)]{CWM}): and
that is enough to ensure that the pentagon commutes.

Any reader who followed the construction of $\alpha$ will see that the
pentagon involves 40 terms of the form $\eta_{m,n}$ or
$\epsln_{m,n}$. Checking that it commutes is therefore an appreciable task,
but in the absence of higher technology there is no alternative. 

\item[Unit:] So far we have only mentioned binary tensor, and not units. To
construct the unit object $I$ of $C(1)$, choose a pseudo-inverse $\psi_0$
to the equivalence of categories $\xi_0: C(0) \go \One$ (in other words, pick
an object of $C(0)$), and define $I\in C(1)$ as (the image of) the composite
\[
\One \goby{\psi_0} C(0) \goby{C(!)} C(1).
\]

\item[Unit isomorphisms:] We need left and right unit isomorphisms 
\[
\lambda_a: I\otimes a \goiso a,
\diagspace
\rho_a:
a\otimes I \goiso a
\]
natural in $a\in C_0$. To define them, choose natural isomorphisms
\[
\eta_0: 1 \goiso \xi_0 \of \psi_0,
\diagspace
\epsln_0: \psi_0\of \xi_0 \goiso 1
\]
in such a way that \tuplebts{\psi_0,\xi_0,\eta_0,\epsln_0} is an adjunction
($\psi_0$ left adjoint to $\xi_0$); this is possible by the result referred
to under `Pentagon' above. (In fact, it's not strictly necessary to use that
general result, since the involvement of the category \One\ makes the
situation trivial; but the argument from general principles is conceptually 
cleaner.) Then $\lambda$ and $\rho$ can be built up from $\eta_{m,n}$'s,
$\epsln_{m,n}$'s, $\eta_0$ and $\epsln_0$. For the record only, $\lambda$ is
defined in Figure~\ref{fig:unitor}, which can be explained in the same way
that Figure~\ref{fig:associator} was.

\begin{figure}

\begin{diagram}
\One \times C(1)	&\rTo^{\psi_0\times 1}	&
C(0) \times C(1)	&\rTo^{C(!)\times 1}		&
C(1)\times C(1)		\\
			&\rdTo<{i}			&
\dTo>{\psi_{0,1}}	&				&
\dTo>{\psi_{1,1}}	\\
			&				&
C(1)			&\rTo^{C(\delta_0)}		&
C(2)			\\
			&				&
			&\rdTo<{1}			&
\dTo>{C(!)}		\\
&&&&C(1)		\\
\end{diagram}

\[
\begin{array}{rl}
&\otimes \of (I\times 1)						\\
=
&C(!) \of \psi_{1,1} \of [C(!) \times 1] \of \xi_{0,1} \of \psi_{0,1}
\of [\psi_0 \times 1]							\\
\goby{1*\eta_{0,1}*1}
&C(!) \of \psi_{1,1} \of [C(!) \times 1] \of \xi_{0,1} \of \psi_{0,1}
\of \xi_{0,1} \of \psi_{0,1} \of [\psi_0 \times 1]			\\
=
&C(!) \of \psi_{1,1} \of \xi_{1,1} \of C(\delta_0) \of \psi_{0,1}
\of \xi_{0,1} \of \psi_{0,1} \of [\psi_0 \times 1]			\\
\goby{1*\epsln_{1,1}*1}
&C(!) \of C(\delta_0) \of \psi_{0,1} \of [\psi_0 \times 1]		\\
=
&\psi_{0,1} \of [\psi_0 \times 1]					\\
=
&\psi_{0,1} \of [\psi_0 \times 1] \of [\xi_0 \times 1] \of \xi_{0,1}
\of i									\\
\goby{1*[\epsln_0 \times 1]*1}
&\psi_{0,1} \of \xi_{0,1} \of i					\\
\goby{\epsln_{0,1}*1}
&i
\end{array}
\]

\caption{Formula for the left unit isomorphism $\lambda$. The functor $i$ is
the canonical isomorphism, and $\delta_0$ is the map from $1$ to $2$ with
image $\{1\}$.}
\label{fig:unitor}
\end{figure}

\item[Triangle:] The final check is that the triangle axiom holds; this is
the `other' coherence axiom for monoidal categories, along with the
pentagon. It is built up out of $\lambda$, $\rho$ and $\alpha$, hence out of
$\eta_{m,n}$'s, $\epsln_{m,n}$'s, $\eta_0$ and $\epsln_0$, and commutes for
the same reason that the pentagon commutes. 
\done

\end{description}
\end{proof}

The statement of the Proposition is rather vague. What the proof actually
consists of is a construction, involving arbitrary choices, of a monoidal
category from a homotopy monoidal category. Soon
(Theorem~\ref{thm:can-hty-mon-cat}) we will give an exact statement
capturing what we have done, and at the same time we will
see that making the arbitrary choices differently only affects the resulting
monoidal category up to isomorphism. To achieve this we consider the
functoriality of the construction.

The category \homset{\HtyAlg}{P}{\emm} was defined with `strict' maps of
homotopy algebras as its morphisms. Thus if \pr{X}{\xi} and \pr{X'}{\xi'} are
homotopy $P$-algebras in \emm, and $\sigma$ a map $\pr{X}{\xi} \go
\pr{X'}{\xi'}$ in \homset{\HtyAlg}{P}{\emm}, then the diagram
\begin{diagram}
X(m)			&\rTo^{\sigma_m}	&X'(m)			\\
\dTo<{X(\Theta)}	&			&\dTo>{X'(\Theta)}	\\
X(n)			&\rTo_{\sigma_n}	&X'(n)			\\
\end{diagram}
in \emm\ commutes (for $\Theta\in\homset{\fmon{P}}{m}{n}$), as do the
squares of Definition~\ref{defn:mon-transf}. If \emm\ is a category like
\Top, \ChCx\ or \Cat\ where it is meaningful to speak of a diagram
commuting `up to homotopy' or `up to isomorphism', then one can consider a
more relaxed kind of morphism of homotopy algebras. But, of course, this is
not meaningful for a general \emm\ in our theory, since all we know about
\emm\ is which maps in it are `homotopy equivalences'. The general
point about weak morphisms of homotopy algebras is returned to in
Chapter~\ref{ch:thoughts}.

The maps in \homset{\HtyAlg}{\Mon}{\Cat} should, therefore, not be considered
as being as weak or lax as ordinary monoidal functors. Nevertheless, any map
in \homset{\HtyAlg}{\Mon}{\Cat} certainly ought to give rise to a monoidal
functor, and this is what the next result says.

\begin{propn}	\label{propn:hty-mon-ftr}
A map of homotopy monoidal categories gives rise to a monoidal functor. That
is, if \pr{C}{\xi} and \pr{C'}{\xi'} are homotopy monoids in \Cat, and if
$C(1)$ and $C'(1)$ are monoidal categories constructed from them as in
Proposition~\ref{propn:hty-mon-cat}, then any map $\pr{C}{\xi} \go
\pr{C'}{\xi'}$ induces (canonically) a monoidal functor $C(1) \go C'(1)$.
\end{propn}

\begin{proof}
Let $\sigma: \pr{C}{\xi} \go \pr{C'}{\xi'}$ be a monoidal transformation. Let
$\psi_{m,n}$ and $\psi_0$ be the (arbitrarily-chosen) functors used in
the construction of \pr{C}{\xi}, and $\eta_{m,n}, \eta_0, \epsln_{m,n},
\epsln_0$ the natural transformations, and similarly $\psi'_{m,n}$ etc.\
for \pr{C'}{\xi'}.

We now construct a monoidal functor from $C(1)$ to $C'(1)$. The functor part
is $\sigma_1: C(1) \go C'(1)$. For the rest of the structure we need
isomorphisms $\sigma_1 (a\otimes b) \goiso \sigma_1 (a) \otimes' \sigma_1
(b)$ (natural in $a,b\in C(1)$) and $\sigma_1 (I) \goiso I'$. The first can
be extracted from the diagram
\begin{diagram}
C(1)\times C(1)	&\rTo^{\psi_{1,1}}	&
C(2)		&\rTo^{C(!)}		&
C(1)		\\
\dTo<{\sigma_1 \times \sigma_1}	&	&
\dTo>{\sigma_2}			&	&
\dTo>{\sigma_1}	\\
C'(1)\times C'(1)&\rTo_{\psi'_{1,1}}	&
C'(2)		&\rTo_{C'(!)}		&
C'(1),		\\
\end{diagram}
in which the right-hand square commutes and the left-hand square commutes up
to isomorphism. The second arises similarly from the diagram
\begin{diagram}
\One	&\rTo^{\psi_0}&C(0)		&\rTo^{C(!)}	&C(1)	\\
\dEquals&		&\dTo>{\sigma_0}&		&\dTo>{\sigma_1}\\
\One	&\rTo_{\psi'_0}&C'(0)		&\rTo_{C'(!)}	&C'(1).	\\
\end{diagram}
Once these coherence isomorphisms have been written down explicitly, it is
just a matter of checking the axioms.
\done
\end{proof}

This construction preserves composites and identities, and so by making a
large number of non-canonical choices we obtain a functor
\[
\fcat{HtyMonCat} \go \fcat{MonCat}.
\]
Here and in what follows,
\[
\label{page:HtyMonCat}
\fcat{HtyMonCat} = \homset{\HtyAlg}{\Mon}{\Cat}
\]
and 
\label{page:MonCat}%
\fcat{MonCat} is the category of (small) monoidal categories and monoidal
functors. 

In order to state the result more precisely, and to get a \emph{canonical}
functor, let us introduce a new category, 
\label{page:HtyMonCat-twiddle}%
\twid{\fcat{HtyMonCat}}. An object
of \twid{\fcat{HtyMonCat}} is a homotopy monoidal category \pr{C}{\xi}
together with a functor
\[
\psi_{m,n}: C(m) \times C(n) \go C(m+n)
\]
and natural isomorphisms 
\[
\eta_{m,n}: 1 \goiso \xi_{m,n}\of\psi_{m,n}
\diagspace
\epsln_{m,n}: \psi_{m,n}\of\xi_{m,n} \go 1
\]
obeying the triangle identities, for each $m$ and $n$, and similarly
$\psi_0, \eta_0, \epsln_0$. (Thus
\tuplebts{\psi_{m,n},\xi_{m,n},\eta_{m,n},\epsln_{m,n}} is an adjoint
equivalence, as is \tuplebts{\psi_0,\xi_0,\eta_0,\epsln_0}.) A map 
\[
\tuplebts{C,\xi,\psi,\eta,\epsln}
\go
\tuplebts{C',\xi',\psi',\eta',\epsln'}
\]
in \twid{\fcat{HtyMonCat}} is just a monoidal transformation $\pr{C}{\xi} \go
\pr{C'}{\xi'}$.

There is a canonical functor
\[
\twid{\fcat{HtyMonCat}} \go \fcat{HtyMonCat}
\]
(forget $\psi$, $\eta$ and $\epsln$) which is full, faithful and surjective
on objects. Hence \twid{\fcat{HtyMonCat}} and \fcat{HtyMonCat} are
equivalent. Our main result can now be stated as follows:

\begin{thm}	\label{thm:can-hty-mon-cat}
There is a canonical functor
\[
\twid{\fcat{HtyMonCat}} \go \fcat{MonCat}.
\]
\done
\end{thm}

The \emph{canonical} functors which have entered our discussion can be
arranged in a diagram,
\begin{diagram}
			&\twid{\fcat{HtyMonCat}}	&		\\
\ldTo(1,2)<{\eqv}	&				&\rdTo(1,2)	\\
\fcat{HtyMonCat}	&				&\fcat{MonCat}.	\\
\end{diagram}
By choosing a pseudo-inverse to the left-hand functor, one obtains a functor
from \fcat{HtyMonCat} to \fcat{MonCat}, as we had in
Proposition~\ref{propn:hty-mon-cat}. But there is no \emph{canonical}
pseudo-inverse, and no \emph{canonical} such functor. In the language of
\cite{Mak}, this diagram depicts an anafunctor from \fcat{HtyMonCat} to
\fcat{MonCat}.

Theorem~\ref{thm:can-hty-mon-cat} has the following corollary, which
says that although the construction in Proposition~\ref{propn:hty-mon-cat}
involves arbitrary choices, these choices do not affect the outcome
significantly.
\begin{cor}
Let \pr{C}{\xi} be a homotopy monoidal category, let $D$ be a monoidal
category arising from \pr{C}{\xi} as in the proof of
Proposition~\ref{propn:hty-mon-cat}, and let $D'$ be another monoidal
category arising in this way via different choices. Then $D$ and $D'$ are
isomorphic objects of \fcat{MonCat}.
\end{cor}
\begin{proof}
Let \triple{\psi}{\eta}{\epsln} be the choices for $D$, and
\triple{\psi'}{\eta'}{\epsln'} those for $D'$. Then
\tuplebts{C,\xi,\psi,\eta,\epsln} and
\tuplebts{C,\xi,\psi',\eta',\epsln'} are isomorphic objects of
\twid{\fcat{HtyMonCat}}, so their images in \fcat{MonCat} under the functor
of Theorem~\ref{thm:can-hty-mon-cat} are also isomorphic. In other
words, $D\iso D'$ in \fcat{MonCat}.
\done
\end{proof}

I want to finish this part with two remarks. First of all, in a brutally
honest world the Propositions above should be called conjectures: I have not
checked every detail. Secondly, the entire theory above can be repeated for
\emph{symmetric} monoidal categories, using homotopy algebras in \Cat\ for
the symmetric operad \CMon. This extension should be absolutely
straightforward. In fact, we have already seen (\ref{sec:Gamma}) that
homotopy symmetric monoidal categories are the same as the
$\Gamma$-categories defined in \cite{SegCCT}, which we call special
$\Gamma$-categories here.

\subsection*{A Generalization}

At no point in our discussion of monoidal categories have we mentioned their
objects and morphisms. To be a little more accurate, we \emph{have} mentioned
them now and then (e.g.\ the `$a$' in $\lambda_a$ and $\rho_a$ in the proof
of Proposition~\ref{propn:hty-mon-cat}), but only as a matter of linguistic
convenience. Essentially, the discussion took place purely in terms of
categories, functors, natural transformations and products. Indeed, only the
purely formal properties of products were used---we did not exploit their
universal property at all.

It follows almost instantly that all the results above hold in any monoidal
2-category, not merely in \Cat. What exactly this means will be clear to
readers experienced in 2-categories, but I will explain it now.

The term `monoidal 2-category' is defined in
Example~\ref{sec:env}\bref{eg:env:mon-2-cat}; for an account of 2-categories
in general, see \cite{KeSt} or \cite{KV}. Recall also
from~\ref{sec:env}\bref{eg:env:mon-2-cat} that any monoidal 2-category
\triple{\cat{N}}{\otimes}{I} has an underlying monoidal category with
equivalences, which we will call \und{\cat{N}}. So, on the one hand, we have
homotopy monoids in \und{\cat{N}}. On the other hand, we have the concept of
a \emph{weak monoid} (also known as a
\label{page:pseudo-monoid}%
\emph{pseudo-monoid}) in \cat{N}. Weak monoids are
defined so that a weak monoid in \Cat\ is a monoidal category in the
traditional sense: thus a weak monoid in \cat{N} consists of 
\begin{itemize}
\item
an object $A$ of \cat{N}
\item
1-cells $t: A\otimes A \go A$, $i: I \go A$
\item
invertible 2-cells
\[
\begin{array}{c}
\alpha: t\of (1\times t) \go t\of(t\times 1),	\\
\lambda: t\of (i\times 1) \go 1, 
\diagspace
\rho: t\of (1\times i) \go 1
\end{array}
\]
satisfying pentagon and triangle axioms. 
\end{itemize}
\emph{Weak maps} of weak monoids are defined in a similar style, so that a
weak map of weak monoids in \Cat\ is a monoidal functor.

The arguments concerning homotopy monoidal categories then give us at once:
\begin{propn}	\label{propn:hty-monoids}
Let \cat{N} be a monoidal 2-category and \und{\cat{N}} the associated
monoidal category with equivalences. Then there is a (non-canonical) functor
\[
\homset{\HtyAlg}{\Mon}{\und{\cat{N}}} 
\go
(\textup{weak monoids and weak maps in }\cat{N}).
\]
\done
\end{propn}
Again, this result can be rephrased to eliminate the element of arbitrary
choice. 

There is just one point where the generalization might not be quite clear,
and this concerns adjoint equivalences in \cat{N}. In \cite[IV.4.1]{CWM} it
is shown that any equivalence of categories `might as well' be an adjoint
equivalence, but it is not obvious from the proof that this is in fact a
general 2-categorical result.

To state the result we need some definitions. Take 0-cells $A$ and $B$ of
\cat{N}, 1-cells
\[
A \oppair{f}{g} B,
\]
and 2-cells
\[
\eta: 1 \go g\of f, 
\diagspace
\epsln: f\of g \go 1. 
\]
Then \tuplebts{f,g,\eta,\epsln} is called an \emph{adjunction} in \cat{N} if
the triangle identities (\cite[IV.1(9)]{CWM}) are satisfied, an
\emph{equivalence} if $\eta$ and $\epsln$ are both invertible, and an
\emph{adjoint equivalence} if both an adjunction and an equivalence. These
definitions have the usual meaning when $\cat{N}=\Cat$; in the case of
adjunctions, $f$ is left adjoint to $g$.

\begin{lemma}	\label{lemma:adjt-equiv}
Let \cat{N} be a 2-category and let 
\[
A \oppair{f}{g} B,
\diagspace
\eta: 1 \goiso g\of f, 
\diagspace
\epsln: f\of g \goiso 1
\]
be an equivalence in \cat{N}. Then there is a 2-cell $\epsln':f\of g \goiso
1$ such that \tuplebts{f,g,\eta,\epsln'} is an adjoint equivalence. 
\end{lemma}
\paragraph*{Remark} When \cat{N} is the 2-category of topological spaces,
continuous maps, and homotopy classes of homotopies, this result is known as
Vogt's Lemma (see \cite[IV.1.14]{KP} and \cite{Vogt}). The general result is
probably due to Street, and has existed as folklore (at least) since the
1970s.
\paragraph*{Proof}
Take $\epsln'$ to be the composite
\[
fg \goby{\epsln^{-1}fg} fgfg \goby{f\eta^{-1}g} fg \goby{\epsln} 1
\]
and check the triangle identities (a tricky but elementary exercise).
\done
\paragraph*{}
Thus we have the result we need on adjoint equivalences, and
obtain Proposition~\ref{propn:hty-monoids}. As usual, the same (probably)
goes for the commutative case.

One final observation will be useful later. By leaving out all mention of
$i$, $\lambda$ and $\rho$ in the definition of weak monoid, we obtain the
definition of \emph{weak semigroup}. By leaving out all mention of $i$,
$\lambda$ and $\rho$ in the proof of~\ref{propn:hty-monoids}, we also obtain:
\begin{propn}	\label{propn:hty-semigps}
Let \cat{N} be a monoidal 2-category and \und{\cat{N}} the associated
monoidal category with equivalences. Then there is a (non-canonical) functor
\[
\homset{\HtyAlg}{\Sem}{\und{\cat{N}}} 
\go
(\textup{weak semigroups in }\cat{N}).
\]
\done
\end{propn}
This result will be employed in the next two sections, on $A_\infty$-spaces
and $A_\infty$-algebras.

\section{$A_\infty$-Spaces}
\label{sec:A-spaces}

The main result of this section is that any homotopy semigroup in $\Topstar$
gives rise to an $A_4$-space. This is an almost immediate consequence of the
results of the previous section. There are also some conjectures on how
this result might extend to give $A_\infty$-algebras, and on related
matters. But before I address any of this, there is a small matter of
basepoints which needs clearing up.

Recall that a semigroup is a set with an associative binary operation, and a
monoid is a semigroup with a two-sided unit; recall also that \Top\ is the
category of spaces and \Topstar\ the category of spaces with basepoint. One
might casually imagine that a semigroup in \Topstar\ is the same thing as a
monoid in \Top: after all, it's just a matter of whether the special point is
regarded as part of the topological data (the basepoint) or the algebraic
data (the unit). But this is not the case, as we saw in
Example~\ref{sec:operads}\bref{eg:operad:Sem}: a semigroup in \Topstar\ is a
topological semigroup together with a distinguished idempotent element, which
need not be a unit.

If we look at the original definition of an $A_\infty$-space (\cite{HAHI})
then we can see that conceptually, an $A_\infty$-space is an up-to-homotopy
version of a semigroup in \Topstar, rather than of a monoid in \Top. This
manifests itself in several ways. Firstly, any semigroup $A$ in \Topstar\ is
naturally an $A_\infty$-space (with trivial structure in dimensions $3$ and
above); the basepoint of $A$ does not need to be a unit. Secondly, it is only
homotopy associativity which is considered (think of the title of
\cite{HAHI}!), and not homotopy unit laws. Thirdly, in the definition of
monoidal category one has both the pentagon and the triangle coherence laws,
whereas the associahedra used to define `$A_\infty$-space' only include the
pentagon, again revealing the semigroupal flavour of $A_\infty$-spaces. 

It therefore seems appropriate to compare $A_\infty$-spaces with
homotopy semigroups in \Topstar, rather than with homotopy monoids in \Top:
and this is what we do here. 

\begin{thm}
Any homotopy semigroup \pr{X}{\xi} in \Topstar\ gives rise to an $A_4$-space,
whose underlying space is $X(1)$. 
\end{thm}
\begin{proof}
There is a 2-category \Topstar\ whose objects and 1-cells are the same as
those of the (1-)category \Topstar, and whose 2-cells are homotopy classes of
homotopies. An equivalence in this 2-category is just a homotopy
equivalence. (All homotopies mentioned here must respect basepoints.)
Moreover, cartesian product $\times$ and the one-point space $1$ make
\triple{\Topstar}{\times}{1} into a monoidal 2-category. So by
Proposition~\ref{propn:hty-semigps}, a homotopy semigroup \pr{X}{\xi} in
\Topstar\ gives rise to a weak semigroup in \Topstar\ with underlying space
$X(1)$. 

Now we only have to see that a weak semigroup in the monoidal 2-category
$\Topstar$ gives rise to an $A_4$-space. (In fact they are more or less the
same thing, as the argument reveals.) I shall only do this informally. A weak
semigroup in \Topstar\ consists of a based space \pr{A}{a_0}, a
multiplication
\[
m_2=t: A\times A \go A
\]
with $m_2(a_0,a_0)=a_0$, and a 2-cell $\alpha$ between the two maps
\[
A^3 \pile{\rTo^{m_2\of(m_2\times 1)}\\
		\rTo_{m_2\of(1\times m_2)}} A,
\]
such that $\alpha$ satisfies the pentagon axiom. Now a 2-cell in \Topstar\ is
a homotopy class of homotopies, so we may pick a representative $m_3$ of
$\alpha$. The pentagon axiom then says that a certain pair of homotopies
$h,h'$ (built up from $m_3$'s) belong to the same homotopy class. Choose a
homotopy $m_4$ between $h$ and $h'$: then $m_4$ is essentially a map
\[
K_4 \times A^4 \go A
\]
where $K_4$ is the solid pentagon (as in \cite{HAHI}), and the data
\tuplebts{A,a_0,m_2,m_3,m_4} thus describes an $A_4$-space. 
\done
\end{proof}

This proof uses 2-category theory, but as far as I know the number 2 has only
one special property: it is the largest value of $n$ for which $n$-category
theory is currently well-understood. I hope that in the near future (weak)
$n$-category theory will be viable, and it will be possible to show
\begin{itemize}
\item that a homotopy semigroup in an $\infty$-category \cat{N} gives rise to
a `weak semigroup' in \cat{N} (where `weak' is meant in an
$\infty$-categorical sense), and
\item that a weak semigroup in the $\infty$-category \Topstar\ is (more or
less) an $A_\infty$-space.
\end{itemize}
So I conjecture that in the Theorem, `4' can be replaced by `$\infty$'. (Of
course, the conjecture might be provable without use of higher-dimensional
categories.) Naturally this should extend to morphisms too: a map of homotopy
semigroups ought to give rise to an $A_\infty$-map.

As for the converse---$A_\infty$-spaces giving rise to homotopy
semigroups---I do not know; see~\ref{sec:inside-Cat} for a similar question
with categories in place of spaces.

\section{$A_\infty$-Algebras}
\label{sec:A-algs}

In the previous section we showed that a homotopy semigroup in \Topstar\
gives rise to an $A_4$-space, by
\begin{itemize}
\item considering \Topstar\ as a monoidal 2-category,
\item employing the result~(\ref{propn:hty-semigps}) that a homotopy
semigroup in a monoidal 2-category \cat{N} gives rise to a weak semigroup in
\cat{N}, then
\item seeing that a weak semigroup in \Topstar\ is an $A_4$-space.
\end{itemize}
In this section we do exactly the same for $A_4$-algebras, replacing
\Topstar\ by \ChCx. So, a homotopy semigroup in \ChCx\ gives rise to an
$A_4$-algebra, and I conjecture that this process can be extended to
give an $A_\infty$-algebra. (See~\cite{HAHII} for the definitions of
$A_n$-algebra and $A_\infty$-algebra.)

This section is laid out as follows. First we look at how \ChCx\ can be made
into a monoidal 2-category, in which the (2-categorical) equivalences are
just the chain homotopy equivalences. Then we examine weak semigroups in
\ChCx, and show them to be essentially the same thing as
$A_4$-algebras. Proposition~\ref{propn:hty-semigps} then applies (as it did
for spaces), so we have a proof that homotopy differential graded non-unital
algebras give rise to $A_4$-algebras. We close with a concrete description of
this process.

Our first task is to describe the monoidal 2-category \ChCx. Objects are
chain complexes and 1-cells are (degree 0) chain maps. 2-cells are homotopy
classes of chain homotopies: but what does that mean? Take chain complexes
$A$ and $B$, chain maps $f,g:A\go B$, and chain homotopies $s,t:f\go g$, as
shown: 
\[
A%
\ctwopar{f}{g}{s}{t}%
B.
\]
Then a \emph{homotopy} (or \emph{secondary homotopy}) $\gamma:s\go t$ consists
of a map 
\[
\gamma: A_{p-2} \go B_p
\]
for each $p\in\integers$, such that 
\[
d\gamma - \gamma d = t-s.
\]
(Note the sign on the left-hand side.) We then say that the homotopies $s$
and $t$ are \emph{homotopic}, and being homotopic is an equivalence
relation. Later (page~\pageref{page:secondary}) we will address the
question of why this is a reasonable definition of secondary homotopy. 

We've now defined the 0-cells, 1-cells and 2-cells of the prospective
monoidal 2-category \ChCx. Composition of 1-cells is done in the usual way,
`vertical' composition of 2-cells
\[
A%
\cthree{f}{g}{h}{s}{t}%
B
\]
is by addition (i.e.\ $t\of s = t+s$), and identities work similarly. The
horizontal composite $s'*s$ of 2-cells
\[
A%
\ctwo{f}{g}{s}%
A'%
\ctwo{f'}{g'}{s'}%
A''
\]
is defined by
\[
(s'*s)(a) = s'g(a) + f's(a)
\]
for $a\in A_p$. The interchange law for 2-categories, which says that
\[
(t'\of s') * (t\of s) = (t'*t)\of(s'*s)
\]
for all suitable $s,s',t,t'$, does hold, but only because we have used
\emph{homotopy classes} of chain homotopies; it does not hold at the level of
ordinary (`primary') homotopies. (A related issue is that $s'*s$ could
equally well have been defined by
\[
(s'*s)(a) = g's(a) + s'f(a);
\]
the choice of one over the other is quite non-canonical.)

So we now have a 2-category \ChCx, and the next step is to endow it with a
monoidal structure. The tensor $\otimes$ of chain complexes and the unit
chain complex $R$ are as usual
(see~\ref{sec:prelims-mon-cats}\bref{eg:mon-cat:ChCx}), and the tensor of
chain maps
\[
A \goby{f} B,
\diagspace
A' \goby{f'} B'
\]
is given by the obvious formula. The tensor of chain homotopies
\[
A%
\ctwo{f}{g}{s}%
B,
\diagspace
A'%
\ctwo{f'}{g'}{s'}%
B'
\]
is given by
\[
(s\otimes s')(a\otimes a') =
s(a) \otimes f'(a') + (-1)^p g(a) \otimes s'(a')
\]
for $a\in A_p$ and $a'\in A'_{p'}$. Once again, this is one of two equally
appropriate formulae, but up to secondary homotopy they are the same.

Finally, then, we have a monoidal 2-category \ChCx. It is clear that an
equivalence inside this 2-category (as defined before
Lemma~\ref{lemma:adjt-equiv}) is just a chain homotopy equivalence.  The
usual symmetry~(\ref{sec:prelims-mon-cats}\bref{eg:mon-cat:ChCx}) in fact
makes \ChCx\ into a \emph{symmetric} monoidal 2-category, but we shall not
need to use this fact.

\subsection*{Aside: secondary chain homotopies}
\label{page:secondary}
\small

Earlier I promised to explain why the definition of homotopy of chain
homotopies is a reasonable one. The situation can be viewed as follows. Let
$A$ and $B$ be topological spaces and let $U$ be the unit interval
$[0,1]$. Then a homotopy between maps $f,g:A\go B$ is, of course, a map
\[
s: U\times A \go B
\]
with $s(0,\dashbk)=f$ and $s(1,\dashbk)=g$. Next, let $s$ and $t$ be two
homotopies between $f$ and $g$; a homotopy between the homotopies $s$ and $t$
is a map
\[
\gamma:U\times U\times A \go B
\]
satisfying
\[
\begin{array}{ll}
\gamma(0,\dashbk,\dashbk)=s,	&
\gamma(1,\dashbk,\dashbk)=t,	\\
\gamma(k,0,\dashbk)=f,		&
\gamma(k,1,\dashbk)=g
\end{array}
\]
for all $k\in U$.

The point is that both of these descriptions can be expressed diagrammatically
in the monoidal category \triple{\Top}{\times}{1}, with the aid of the maps
\[
i_0, i_1: 1\go U, 
\diagspace
j:U\go 1,
\]
where $i_0$ and $i_1$ have respective values $0$ and $1$. Now let's mimic
this in the monoidal category \triple{\ChCx}{\otimes}{R}: take $U$ to be
the chain complex
\begin{diagram}[height=2em]
\cdots	&0&\rTo	&R		&\rTo	&R\oplus R	&\rTo &0 &\cdots \\
	& &	&r		&\rGoesto&(-r,r)	& & & \\
\end{diagram}
with $R$ in degree 1 and $R\oplus R$ in degree 0, define $i_0, i_1: R \go
U$ by the first and second inclusions of $R$ into $R\oplus R$, and define $j:
U \go R$ by the addition map from $R\oplus R$ to $R$. A \emph{homotopy}
between chain maps can then be defined as a suitable map $U\otimes A \go B$,
as in the topological case, and this turns out to be equivalent to the usual
definition of chain homotopy. A \emph{secondary homotopy} can similarly be
defined as a map
\[
U\otimes U\otimes A \go B
\]
satisfying suitable `boundary conditions', and, with a significant amount of
calculation, this turns out to be equivalent to the very simple definition
given originally.

Further thoughts of this kind are laid out in \cite[III.3]{KP}.

\subsection*{}
\normalsize

Returning to the main flow, we have exhibited \ChCx\ as a monoidal 2-category
and now wish to look at weak semigroups in it. Such a thing consists of a
chain complex $A$, a (degree 0) chain map 
\[
m_2=t:A^{\otimes 2} \go A,
\]
and a homotopy class $\alpha$ of chain homotopies 
\begin{equation}	\label{eq:hty-class}
A^{\otimes 3}%
\ctwo{m_2\of(m_2\otimes 1)}{m_2\of(1\otimes m_2)}{}%
A,
\end{equation}
satisfying the pentagon axiom. Choose a representative $m_3$ of the class
$\alpha$. Then $m_3$ is a family of homomorphisms
\[
(A^{\otimes 3})_p \go A_{p+1}
\]
($p\in\integers$), and the fact that $m_3$ is a homotopy between the two maps
in~\bref{eq:hty-class} says that
\[
\begin{array}{l}
d(m_3(a_1,a_2,a_3))		\\
\mbox{} + m_3(da_1,a_2,a_3) + (-1)^{p_1}m_3(a_1,da_2,a_3) 
+ (-1)^{p_1+p_2}m_3(a_1,a_2,da_3)	\\
\ \ = \ \ -m_2(m_2(a_1,a_2),a_3) + m_2(a_1,m_2(a_2,a_3))
\end{array}
\]
for all $a_1 \in A_{p_1}$, $a_2 \in A_{p_2}$ and $a_3 \in A_{p_3}$. Finally,
the fact that $\alpha$ satisfies the pentagon identity means that there is a
secondary homotopy $m_4$ between a certain pair of homotopies built up as
composites of $m_3$'s. Thus $m_4$ is a family of homomorphisms
\[
(A^{\otimes 4})_p \go A_{p+2}
\]
($p\in\integers$), and the equation `$d\gamma - \gamma d = t-s$' in the
definition of secondary homotopy says that for $a_1\in A_{p_1}$, $a_2\in
A_{p_2}$, $a_3\in A_{p_3}$ and  $a_4\in A_{p_4}$, 
\begin{eqnarray}	
\lefteqn{d(m_4(a_1,a_2,a_3,a_4))
-m_4(da_1,a_2,a_3,a_4) 
-(-1)^{p_1}m_4(a_1,da_2,a_3,a_4)}		\nonumber\\
\lefteqn{\mbox{} -(-1)^{p_1+p_2}m_4(a_1,a_2,da_3,a_4) 
-(-1)^{p_1+p_2+p_3}m_4(a_1,a_2,a_3,da_4)}	\nonumber\\
	&=	&
-m_3(m_2(a_1,a_2),a_3,a_4)
+m_3(a_1,m_2(a_2,a_3),a_4)
-m_3(a_1,a_2,m_2(a_3,a_4))			\nonumber\\
	&	&
\mbox{} +m_2(m_3(a_1,a_2,a_3),a_4)
+(-1)^{p_1}m_2(a_1,m_3(a_2,a_3,a_4)).\label{eq:m4}
\end{eqnarray}
So the structure \tuplebts{A,m_2,m_3,m_4} is precisely an $A_4$-algebra. We
therefore have:
\begin{thm}
Any homotopy differential graded non-unital algebra \pr{X}{\xi} gives rise to
an $A_4$-algebra, whose underlying chain complex is $X(1)$. \done
\end{thm}

Let us now look more directly at how a homotopy semigroup
\[
\pr{X}{\xi}: 
\triple{\Delta_{\mr{surj}}}{+}{0} \go \triple{\ChCx}{\otimes}{R}
\]
leads to an $A_4$-algebra $A$. (See~\ref{sec:free-mon}\bref{eg:free-mon:Sem}
for the definition of $\Delta_{\mr{surj}}$.) First of all,
$A=X(1)$. Secondly, choose a chain homotopy inverse $\psi_{m,n}$ to each
$\xi_{m,n}$, as shown: 
\[
X(m+n)
\pile{\rTo^{\xi_{m,n}}\\ \lTo_{\psi_{m,n}}}
X(m) \otimes X(n).
\]
Then $m_2: A^{\otimes 2} \go A$ is defined as the composite
\[
X(1) \otimes X(1) \goby{\psi_{1,1}} X(2) \goby{X(!)} X(1)
\]
where $!$ is the unique map $2\go 1$ in $\Delta_{\mr{surj}}$. To describe
$m_3$, consider the diagram at the top of Figure~\ref{fig:associator}
(page~\pageref{fig:associator}), with $C$'s changed to $X$'s and $\times$'s
to $\otimes$'s. For the same reasons as given then, each inner square of the
diagram---hence the whole diagram---commutes up to homotopy. This says that
$m_2\of (m_2\otimes 1)$ and $m_2\of(1\otimes m_2)$ are chain-homotopic, and
indeed we can construct a particular such homotopy, $m_3$.

Using the formula in~\ref{fig:associator}, $m_3$ is a sum of 8 terms. This
means that the right-hand side of equation~\bref{eq:m4} is a sum of 40 terms,
and finding an $m_4$ to satisfy it would be an enormous task if attempted
from cold. However, the general method of Proposition~\ref{propn:hty-semigps}
provides an $m_4$ automatically.

Just as for $A_n$-spaces, I know of no reason why the process should
stop at $A_4$, and I conjecture that there are also maps $m_5$, $m_6$,
\ldots\ making $X(1)$ into an $A_\infty$-algebra. Similarly, it is plausible
that maps of homotopy semigroups in \ChCx\ give rise to $A_\infty$-maps.
Bearing in mind that \ChCx\ is a \emph{symmetric} monoidal category, it might
also be possible to do the same things for homotopy d.g.\ \emph{commutative}
non-unital algebras: such a structure might give rise to a $C_4$-algebra, and
perhaps a $C_\infty$-algebra.

\chapter{Other Examples of Homotopy Algebras}
\label{ch:other}

The previous chapter covered homotopy monoids and homotopy semigroups in some
detail. In this chapter we look at various other examples of homotopy
algebras, and develop some further theory concerning homotopy algebras in
general. 

Section~\ref{sec:misc-egs} is an assortment of examples of homotopy algebras:
homotopy graded Lie and Gerstenhaber algebras (both of which are homotopy
algebras in the enriched sense), homotopy monoids-with-involution, and
`homotopy homotopy algebras' (e.g.\ homotopy $L_\infty$-algebras). The whole
section raises more questions than it answers, and in particular poses a
concrete question concerning Hochschild cochains~(\ref{subsec:Ger}).

Section~\ref{sec:iterated} fulfils a promise made in
Chapter~\ref{ch:hty-mons}: to put a natural homotopy-algebraic structure on
an $n$-fold loop space, for any $1\leq n < \infty$. In order to do this we
have to develop a theory of `homotopy \bftuple{P_1}{P_n}-algebras' for any
family $(P_i)$ of operads, which we do briefly. An $n$-fold loop space is
then an `$n$-fold homotopy monoid', that is, a homotopy
$(\underbrace{\Mon,\ldots, \Mon}_{n})$-algebra.

The final two sections,~\ref{sec:mon-2-cat} and~\ref{sec:inside-Cat}, are
about comparing different notions of weakened algebraic structure. Any
monoidal 2-category has an underlying monoidal category with equivalences (as
we saw in Example~\ref{sec:env}\bref{eg:env:mon-2-cat}), one example being
\Cat. In~\ref{sec:mon-2-cat} we formulate a notion of `weak $P$-algebra' in
any monoidal 2-category, and extend the method of Chapter~\ref{ch:hty-mons}
to show that any homotopy $P$-algebra gives rise to a weak $P$-algebra. One
naturally wants to know whether it is possible to go in the opposite
direction too (weak algebras giving homotopy algebras); I cannot provide an
answer in general, but~\ref{sec:inside-Cat} shows how this is possible in the
case of \Cat. 

\section{Miscellaneous Examples}
\label{sec:misc-egs}

\subsection{Graded Lie algebras}
\label{subsec:GrLie}

So far we have not paid very much attention to the enriched setting: we have
said how to \emph{define} homotopy algebras for an enriched
operad~(\ref{sec:enr-defn}), but have not done much by way of examples. As
compensation, we now examine in detail the definition of a homotopy graded
Lie algebra. That is, we examine homotopy \GrLie-algebras in \ChCx.

So, take the \GrAb-enriched symmetric operad
\GrLie~(\ref{sec:enr-ops}\bref{eg:enr-op:GrLie}), which is generated by an
element
\[
[\emptybk,\emptybk] \in (\GrLie(2))_{-1}
\]
subject to equations
\begin{eqnarray*}
[\emptybk,\emptybk]
+[\emptybk,\emptybk].\tau 			&=	&0	\\
\mbox{}
[[\emptybk,\emptybk], \emptybk]
+[[\emptybk,\emptybk], \emptybk].\sigma
+[[\emptybk,\emptybk], \emptybk].\sigma^2	&=	&0
\end{eqnarray*}
where $\tau\in S_2$ is a 2-cycle and $\sigma\in S_3$ is a 3-cycle. 

Consequently, an algebra for \GrLie\ in $\GrMod_R$ is a graded
$R$-module $A$ together with a binary operation of degree $-1$, satisfying
the equations
\begin{eqnarray*}
[a,b] + (-1)^{pq}[b,a] 	&=	&0	\\
(-1)^{rp}[[a,b],c] + (-1)^{pq}[ [b,c],a] + (-1)^{qr}[ [c,a],b] &= &0	
\end{eqnarray*}
for $a\in A_p, b\in A_q, c\in A_r$. The signs arise from the symmetry map in
$\GrMod_R$ (see~\ref{sec:prelims-mon-cats}\bref{eg:mon-cat:GrMod}). So as
expected, a \GrLie-algebra is a graded Lie algebra in the usual sense.

From this point on, most of what we have to say about \GrLie\ applies equally
to all \GrAb-enriched symmetric operads. 

Homotopy algebras are defined via the \GrAb-enriched symmetric strict
monoidal category \fmon{\GrLie}, whose objects are the natural numbers and
whose `hom-objects' are the graded abelian groups
\[
\homset{\fmon{\GrLie}}{m}{n}
=
\bigoplus_{f\in \homset{\Phi}{m}{n}}
\GrLie(f^{-1}\{0\}) \otimes\cdots\otimes
\GrLie(f^{-1}\{n-1\}).
\]
(It's not hard to see that $\GrLie(0)=0$, that $\GrLie(n)$ is concentrated in
degree $1-n$ for $n\geq 1$, and that \homset{\fmon{\GrLie}}{m}{n} is therefore
concentrated in degree $n-m$; but this doesn't matter for the present
account.)

Much as in Example~\ref{sec:enr-cats}\bref{eg:enr-cat:ChCx}, $\ChCx=\ChCx_R$
can be viewed as a \GrAb-enriched symmetric monoidal category, with
\ehom{\ChCx}{C}{D} being the graded abelian group whose degree $k$ part is
the abelian group of all degree $k$ chain maps from $C$ to $D$. The ordinary
category underlying this \GrAb-enriched category is the usual \ChCx, in which
\homset{\ChCx}{C}{D} is the set of all degree $0$ chain maps $C\go D$. A
homotopy \GrLie-algebra in \ChCx\ consists of a \GrAb-enriched symmetric
monoidal functor
\[
\pr{X}{\xi}: \fmon{\GrLie} \go \ChCx
\]
such that the components of $\xi$ are homotopy equivalences. Taking this
apart further, a homotopy graded Lie algebra consists of
\begin{itemize}
\item a sequence $X(0), X(1), \ldots$ of chain complexes
\item a chain homotopy equivalence $\xi_0: X(0) \go R$ (where $R$ is the unit
chain complex---see~\ref{sec:prelims-mon-cats}\bref{eg:mon-cat:ChCx})
\item a chain homotopy equivalence
\[
\xi_{n,n'}: X(n+n') \go X(n) \otimes X(n')
\]
for each $n,n' \geq 0$
\item for each map $f:m \go n$ of finite sets and each 
\[
\theta_1 \in (\GrLie(f^{-1}\{0\}))_{p_1}, \ldots, 
\theta_n \in (\GrLie(f^{-1}\{n-1\}))_{p_n},
\]
a chain map 
\[
X(f;\theta_1, \ldots, \theta_n): X(m) \go X(n)
\]
of degree $p_1 + \cdots + p_n$.
\end{itemize}
(Note that $X(n) \eqv X(1)^{\otimes n}$ as usual.) This data satisfies
various axioms. The expression $X(f; \theta_1, \ldots, \theta_n)$ preserves
addition of $\theta_i$'s, and $X$ preserves composites of the morphisms $(f;
\theta_1, \ldots, \theta_n)$ in \fmon{\GrLie}, and similarly preserves
identities. The maps $\xi_0$ and $\xi_{n,n'}$ obey the usual coherence axioms
for a colax monoidal functor (given in~\ref{defn:mon-ftr}). Finally,
$\xi_{n,n'}$ is natural in $n$ and $n'$: if $f:m \go n$ and
$f': m' \go n'$ are maps of finite sets, and if
\[
\begin{array}{l}
\theta_1 \in (\GrLie(f^{-1}\{0\}))_{p_1}, \ldots, 
\theta_n \in (\GrLie(f^{-1}\{n-1\}))_{p_n},	\\
\theta'_1 \in (\GrLie(f'^{-1}\{0\}))_{p'_1}, \ldots, 
\theta'_{n'} \in (\GrLie(f'^{-1}\{n'-1\}))_{p'_{n'}},
\end{array}
\]
then the diagram
\begin{diagram}[width=5em,tight]	
X(m+m')	&\rTo^{\xi_{m,m'}}	&X(m) \otimes X(m')	\\
\dTo<{X(f+f'; \theta_1, \ldots, \theta_n, \theta'_1, \ldots, \theta'_{n'})}	&&
\dTo>{X(f;\theta_1, \ldots, \theta_n) \otimes X(f';\theta'_1, \ldots, \theta'_{n'})}\\
X(n+n')	&\rTo^{\xi_{n,n'}}	&X(n) \otimes X(n')	\\
\end{diagram}
commutes. The horizontal maps here are of degree $0$, and the vertical maps
are of degree $(p_1 +\cdots+ p_n + p'_1 +\cdots+ p'_{n'})$.

(In fact, since $\GrLie(k)$ is concentrated in degree $1-k$, one might as well
consider $X(f;\theta_1, \ldots, \theta_n)$ only when 
\[
p_1 = 1 - |f^{-1}\{0\}|, \ldots, p_n = 1 - |f^{-1}\{n-1\}|.)
\]

This concludes the description of homotopy graded Lie algebras. It would be
interesting to compare these structures with $L_\infty$-algebras, just as
homotopy d.g.\ algebras were compared with $A_\infty$-algebras
in~\ref{sec:A-spaces}.

\subsection{Gerstenhaber algebras}
\label{subsec:Ger}

We have a category \homset{\HtyAlg}{\Ger}{\ChCx} of homotopy Gerstenhaber
algebras, where \Ger\ is the symmetric \GrAb-enriched operad
of~\ref{sec:enr-ops}\bref{eg:enr-op:Ger}. Various other notions of `homotopy
Gerstenhaber algebra' are laid out in the paper of this title by Voronov, and
once again a comparison would be nice but is not attempted. In particular,
Deligne's Conjecture implies that the Hochschild cochain complex
$C^{\blob}(A)$ of an associative algebra $A$ is a `homotopy Gerstenhaber
algebra' in any reasonable sense of the phrase. An interesting challenge,
therefore, is to prove directly that $C^{\blob}(A)$ is a homotopy
Gerstenhaber algebra in our sense of the phrase. A sub-challenge, which does
not involve ideas of enrichment, is to show that $C^{\blob}(A)$ is a homotopy
\CMon-algebra, i.e.\ a homotopy d.g.\ commutative algebra: for any homotopy
Gerstenhaber algebra is certainly a homotopy \CMon-algebra. Roughly speaking,
this means:

\paragraph*{Problem} Given an associative algebra $A$ over $R$, find a
functor $X$ from $\Phi$ (the skeletal category of finite sets) to $\ChCx_R$
(the category of chain complexes over $R$), such that
\begin{itemize}
\item $X(1) \iso C^{\blob}(A)$
\item there is a canonical chain homotopy equivalence
\[
X(0) \go R
\]
\item there is a canonical chain homotopy equivalence
\[
X(m+n) \go X(m) \otimes X(n)
\]
for each $m,n\geq 0$.
\end{itemize}

(For an account of Deligne's Conjecture, see \cite{Kon}. Following a tangled
history of proofs and refutations, it appears that the Conjecture is now a
Theorem.)

\subsection{Monoids with involution}

We have shown that any loop space is a homotopy topological
monoid~(\ref{sec:loops}), in the sense of being a homotopy algebra in \Top\
for the non-symmetric operad \Mon. We have also suggested
(p.~\pageref{page:loop-hty-sym}) that a loop space is a homotopy monoid in
the sense of being a homotopy algebra for the symmetric operad \Sym. Since a
loop can be travelled backwards, any loop space is also acted on (strictly)
by the 2-element group $C_2$. Combining the two structures, one would
therefore expect any loop space to be a homotopy monoid-with-involution: that
is, a homotopy \Inv-algebra in \triple{\Top}{\times}{1}, where \Inv\ is the
symmetric operad defined in~\ref{sec:operads}\bref{eg:operad:Inv}. I do not
know whether this is, in fact, true. It would be enough to show that $S^1$
is a homotopy monoid-with-involution in \triple{\Top_*^{\op}}{\wej}{1}, i.e.\
to construct a suitable colax symmetric monoidal functor
\[
\pr{W}{\omega}: \triple{\fmon{\Inv}}{+}{0} \go \triple{\Top_*^{\op}}{\wej}{1}
\]
with $W(1)=S^1$; for then we could compose with the functor
\homset{\Topstar}{\dashbk}{B}, as in~\ref{sec:loops}.

\subsection{Homotopy homotopy algebras}

I do not know what to make of the following bizarre family of
examples. Consider, for instance, the notion of a homotopy Lie algebra: that
is, a chain complex which is a graded Lie algebra `up to higher
homotopy'. The present work formalizes this idea as a homotopy \GrLie-algebra
in \ChCx; on the other hand, the usual way to formalize it is as an
$L_\infty$-algebra. An $L_\infty$-algebra is an algebra in \ChCx\ for a
certain \GrAb-enriched operad, $L_\infty$. But this means that we can also
consider \emph{homotopy} $L_\infty$-algebras in \ChCx, i.e.\ the category
\[
\homset{\HtyAlg}{L_\infty}{\ChCx}.
\]
A homotopy $L_\infty$-algebra is then a `homotopy homotopy Lie
algebra'. Similarly, one might consider homotopy $A_\infty$-, $B_\infty$-,
$C_\infty$-, $G_\infty$-, \ldots algebras, which are all `homotopy homotopy
algebras'. It would perhaps be desirable to show that in some sense, any
homotopy homotopy algebra is in fact a mere homotopy algebra: `working up to
homotopy is idempotent'.

\section{Iterated Loop Spaces}
\label{sec:iterated}

We have already seen that any loop space is a homotopy monoid in \Top. We
have also seen that a homotopy commutative monoid in \Top\ is the same thing
as a special $\Gamma$-space, and the relation between special $\Gamma$-spaces
and infinite loop spaces is well-established (see e.g.\ \cite{Ad}). One might
therefore ask what structure an $n$-fold loop space carries, in our theory,
when $1<n<\infty$. This section provides an answer.

Before going into detail, here is a sketch of the ideas. Ultimately we will
show that any $n$-fold loop space `has $n$ commuting homotopy monoid
structures on it', or is `an $n$-fold homotopy monoid'. To state this
precisely we need some general definitions. Given operads $(P_i)$, a
`\bftuple{P_1}{P_n}-algebra' in a symmetric monoidal category \emm\ is an
object $A$ of \emm\ which is an algebra for each of the operads $P_i$, in
such a way that the algebra structures commute with each other. An extension
of Theorem~\ref{thm:alt-alg} says that \bftuple{P_1}{P_n}-algebras in \emm\
are essentially the same as multi-monoidal functors
\[
\fmon{P_1} \times\cdots\times \fmon{P_n} \go \emm,
\]
where `multi-monoidal functors' are to monoidal functors as multilinear maps
are to linear maps. A `homotopy \bftuple{P_1}{P_n}-algebra' is defined as a
\emph{colax} multi-monoidal functor \pr{X}{\xi} in which the components of
$\xi$ are equivalences. We will show that the $n$-sphere $S^n$ is a homotopy
\bftuple{\Mon}{\Mon}-algebra in \triple{\Top_*^{\op}}{\wej}{1}, and it
follows easily that any $n$-fold loop space is a homotopy
\bftuple{\Mon}{\Mon}-algebra in \triple{\Top}{\times}{1}.

In detail: let $n\geq 0$, let $\cat{L}_1, \ldots, \cat{L}_n$ be monoidal
categories, and let \emm\ be a symmetric monoidal category. A \emph{colax
multi-monoidal functor\footnote{with apologies for the name}}
\[
\pr{X}{\xi}: \cat{L}_1 \times\cdots\times \cat{L}_n \go \emm
\]
consists of
\begin{itemize}
\item a functor $X: \cat{L}_1 \times\cdots\times \cat{L}_n \go \emm$
\item a map
\[
\xi_{\ldots, L_{i-1}, 0, L_{i+1}, \ldots}:
X(\ldots, L_{i-1}, I, L_{i+1}, \ldots) \go I
\]
for each $i\in \{1, \ldots n\}$ and $L_j \in \cat{L}_j$ $(j\neq i)$
\item a map
\begin{eqnarray*}
\lefteqn{\xi_{\ldots, L_{i-1}, \pr{L_i}{L_{i}'}, L_{i+1}, \ldots}:
X(\ldots, L_{i-1}, L_i \otimes L_{i}', L_{i+1}, \ldots) \go}	\\
&X(\ldots, L_{i-1}, L_i, L_{i+1}, \ldots) \otimes
X(\ldots, L_{i-1}, L_{i}', L_{i+1}, \ldots)
\end{eqnarray*}
for each $i \in \{1,\ldots,n\}$, $L_j\in \cat{L}_j$, and $L_{i}'\in\cat{L}_i$.
\end{itemize}
These maps are required to be natural in all the components $L_j$ and to obey
the usual colax monoidal functor axioms componentwise, so that for each $i$
and $L_1, \ldots, L_{i-1}, L_{i+1}, \ldots, L_n$, the pair
\[
\pr{X(\ldots, L_{i-1}, \dashbk, L_{i+1}, \ldots)}{\xi_{\ldots, L_{i-1},
\dashbk, L_{i+1}, \ldots}}
\]
forms a colax monoidal functor $\cat{L}_i \go \emm$. They are also required
to commute with one another, which means that if $1\leq i < j \leq n$ then
the four diagrams in Figure~\ref{fig:multi-mon} commute. In the figure,
$\mb{X}\pr{J}{K}$ is an abbreviation for 
\[
X(L_1, \ldots, L_{i-1}, J, L_{i+1}, \ldots, L_{j-1}, K, L_{j+1}, \ldots,
L_n),
\]
and the diagrams are required to commute for all $L_k\in \cat{L}_k$ ($1\leq
k\leq n$) and $L_i'\in \cat{L}_i, L_j'\in \cat{L}_j$. All the arrows in the
diagrams are the obvious components of $\xi$, except for the three labelled
as isomorphisms.
\begin{figure}
 
\begin{diagram}
\mb{X}\pr{L_i\otimes L_i'}{L_j\otimes L_j'}	&\rTo	
&\mb{X}\pr{L_i}{L_j\otimes L_j'} \otimes \mb{X}\pr{L_i'}{L_j\otimes L_j'}\\
	&	&\dTo	\\
\dTo	&	&
\mb{X}\pr{L_i}{L_j} \otimes \mb{X}\pr{L_i}{L_j'} \otimes
\mb{X}\pr{L_i'}{L_j} \otimes \mb{X}\pr{L_i'}{L_j'}		\\
	&	&\dTo>{\diso}	\\
\mb{X}\pr{L_i\otimes L_i'}{L_j} \otimes \mb{X}\pr{L_i\otimes L_i'}{L_j'}
&\rTo	&
\mb{X}\pr{L_i}{L_j} \otimes \mb{X}\pr{L_i'}{L_j} \otimes
\mb{X}\pr{L_i}{L_j'} \otimes \mb{X}\pr{L_i'}{L_j'}		\\
\end{diagram}

\[
\begin{diagram}
\mb{X}\pr{L_i \otimes L_i'}{I}	&\rTo	&
\mb{X}\pr{L_i}{I} \otimes \mb{X}\pr{L_i'}{I}	\\
\dTo	&		&\dTo			\\
I	&\rTo_{\diso}	&I\otimes I		\\
\end{diagram}
\diagspace
\begin{diagram}
\mb{X}\pr{I}{L_j \otimes L_j'}	&\rTo	&I	\\
\dTo	&	&\dTo>{\diso}			\\
\mb{X}\pr{I}{L_j} \otimes \mb{X}\pr{I}{L_j'}	&
\rTo	&I\otimes I				\\
\end{diagram}
\]

\begin{diagram}
\mb{X}\pr{I}{I}	&
\pile{\rTo^{\xi_{\ldots,I,\ldots,0,\ldots}} \\ 
 \rTo_{\xi_{\ldots,0,\ldots,I,\ldots}}}		&I	\\
\end{diagram}

\caption{Commutativity axioms for a colax multi-monoidal functor.} 
\label{fig:multi-mon}
\end{figure}
There is an obvious notion of \emph{monoidal transformation} between colax
multi-monoidal functors.

Now let $P_1, \ldots, P_n$ be non-symmetric operads and let \emm\ be a
symmetric monoidal category with equivalences. A \emph{homotopy
\bftuple{P_1}{P_n}-algebra in \emm} is a colax multi-monoidal functor 
\[
\pr{X}{\xi}: \fmon{P_1} \times\cdots\times \fmon{P_n} \go \emm
\]
such that 
each component of $\xi$ is an equivalence. With monoidal transformations as
maps, this gives a category \multihom{\HtyAlg}{P_1, \ldots, P_n}{\emm} of
homotopy \bftuple{P_1}{P_n}-algebras in \emm. 

To see why this is a reasonable definition, consider first `genuine'
\bftuple{P_1}{P_n}-algebras. If $P$ and $Q$ are two operads, \emm\ a
symmetric monoidal category, and $A$ an object of \emm\ endowed with both
$P$-algebra and $Q$-algebra structures, let us say that the two algebra
structures \emph{commute} if for all $\theta\in P(m)$ and $\chi\in Q(k)$, the
diagram
\begin{diagram}
	&	&(A^{\otimes m})^{\otimes k}	&
\rTo^{(\ovln{\theta})^{\otimes k}}	&A^{\otimes k}	\\
	&\ruTo<{\textrm{symmetry}}	&	&	&	\\
(A^{\otimes k})^{\otimes m}	&	&	&	&
\dTo>{\ovln{\chi}}	\\
\dTo<{(\ovln{\chi})^{\otimes m}}	&	&	&	&	\\
A^{\otimes m}	&	&\rTo_{\ovln{\theta}}	&	&A	\\
\end{diagram}
commutes. If $P_1, \ldots, P_n$ are operads and \emm\ a symmetric monoidal
category, a \emph{\bftuple{P_1}{P_n}-algebra in \emm} is an object $A$ of
\cat{M} with the structure of a $P_i$-algebra for each $i\in \{1, \ldots
n\}$, such that the $P_i$-algebra and $P_j$-algebra structures commute
whenever $i\neq j$. One thus obtains a category \multihom{\Alg}{P_1, \ldots,
P_n}{\emm} of \bftuple{P_1}{P_n}-algebras in \emm.

For example, a \pr{\Sem}{\Sem}-algebra in \Set\ is a set $A$ equipped with
two associative binary operations, $\cdot$ and $*$, satisfying the
`interchange law':
\[
(a*b)\cdot (a'*b') = (a\cdot a') * (b\cdot b').
\]

Now, take a symmetric monoidal category \emm\ and define the equivalences in
\emm\ to be just the isomorphisms. Let $P_1, \ldots, P_n$ be non-symmetric
operads. Then \multihom{\HtyAlg}{P_1, \ldots, P_n}{\emm} is
\multihom{\Mon}{P_1, \ldots, P_n}{\emm}, the category of \emph{multi-monoidal
functors} $P_1 \times\cdots\times P_n \go \emm$: that is, those colax
multi-monoidal functors \pr{X}{\xi} in which all the components of $\xi$ are
isomorphisms. Theorem~\ref{thm:alt-alg} can be generalized to give an equivalence of
categories
\[
\multihom{\Mon}{P_1, \ldots, P_n}{\emm}
\eqv
\multihom{\Alg}{P_1, \ldots, P_n}{\emm},
\]
with~\ref{thm:alt-alg} being the case $n=1$; the algebra corresponding to a
multi-monoidal functor \pr{X}{\xi} has $X(1,\ldots,1)$ as its underlying
object. Hence when \emm\ has only trivial equivalences,
\[
\multihom{\HtyAlg}{P_1, \ldots, P_n}{\emm}
\eqv
\multihom{\Alg}{P_1, \ldots, P_n}{\emm},
\]
just as in the case $n=1$.

(All of these definitions can be repeated, with minor modifications, for the
case of \emph{symmetric} operads $P_i$. Moreover, there is no need for the
family $(P_i)$ of operads to be finite; everything above works just as well
for infinite families.)

\small
\paragraph*{Aside} An alternative way of making the definitions is to observe
that \homset{\fcat{(S)Colax}}{P}{\emm} is naturally a symmetric monoidal
category with equivalences, for any (symmetric) operad $P$ and symmetric
monoidal category with equivalences \emm. (Tensor and equivalences are
defined pointwise.) In fact, the subcategory \homset{\HtyAlg}{P}{\emm} is
also a symmetric monoidal category with equivalences. So we could define
\[
\multihom{\HtyAlg}{Q,P}{\emm} =
\homset{\HtyAlg}{Q}{\homset{\HtyAlg}{P}{\emm}}
\]
for any operads $Q$ and $P$; and we could iterate this idea in order to
define homotopy \bftuple{P_1}{P_n}-algebras for $n>2$. This definition is
equivalent to the one given above.

\normalsize
\paragraph*{}

We can now return to iterated loop spaces. Our result is:
\begin{thm}	\label{thm:n-fold-monoid}
Any $n$-fold loop space is an $n$-fold homotopy monoid. That is, if $B$ is a
space with basepoint then there is a homotopy $(\underbrace{\Mon, \ldots,
\Mon}_{n})$-algebra \pr{X}{\xi} in \triple{\Top}{\times}{1} with $X(1,
\ldots, 1) = \homset{\Topstar}{S^n}{B}$.
\end{thm}

\begin{proof}
First we show that $S^n$ is an $n$-fold homotopy comonoid,
i.e.\ a homotopy \bftuple{\Mon}{\Mon}-algebra in
\triple{\Top_*^{\op}}{\wej}{1}. Let 
\[
\pr{W}{\omega}: \triple{\Delta}{+}{0} \go \triple{\Top_*^{\op}}{\wej}{1}
\]
be the colax monoidal functor of Lemma~\ref{lemma:S1-comonoid}, exhibiting
$S^1$ as a homotopy comonoid. Also
let
\[
Z: \Top_*^n \go \Topstar
\]
be the $n$-fold smash product,
\[
Z\bftuple{A_1}{A_n} = A_1 \smsh\cdots\smsh A_n.
\]
Observe that since \smsh\ distributes over \wej, $Z$ naturally becomes a
multi-monoidal functor
\[
\pr{Z}{\zeta}: \triple{\Topstar}{\wej}{1}^n \go \triple{\Topstar}{\wej}{1};
\]
observe moreover that $Z$ preserves homotopy equivalences. Assembling all
of this, we get a composite
\[
\triple{\Delta}{+}{0}^n \goby{\pr{W}{\omega}^n}
\triple{\Top_*^{\op}}{\wej}{1}^n \goby{\pr{Z}{\zeta}}
\triple{\Top_*^{\op}}{\wej}{1},
\]
and this is a colax multi-monoidal functor \pr{S}{\sigma} in which the
components of $\sigma$ are equivalences. Thus \pr{S}{\sigma} defines an
$n$-fold homotopy comonoid, and
\begin{eqnarray*}
S(1,\ldots,1)	&=	&W(1) \smsh\cdots\smsh W(1)	\\
		&=	&S^1 \smsh\cdots\smsh S^1	\\
		&=	&S^n
\end{eqnarray*}
as required.

To finish the proof we simply use the observation made
in~\ref{sec:loops} that \homset{\Topstar}{\dashbk}{B} defines a
homotopy-preserving monoidal functor
\[
\triple{\Top_*^{\op}}{\wej}{1} \go \triple{\Top}{\times}{1}.
\]
Composing this with 
\[
\pr{S}{\sigma}: \triple{\Delta}{+}{0}^n \go 
\triple{\Top_*^{\op}}{\wej}{1}
\]
yields an $n$-fold homotopy monoid
\[
\triple{\fmon{\Mon}}{+}{0}^n =
\triple{\Delta}{+}{0}^n \go \triple{\Top}{\times}{1},
\]
whose value at \bftuple{1}{1} is \homset{\Topstar}{S^n}{B}.
\done
\end{proof}

The Theorem swiftly implies that the higher homotopy groups of a space are
abelian, as we see in~\ref{sec:change-egs}\bref{eg:change:path-components}.

\section{Inside a Monoidal 2-Category}
\label{sec:mon-2-cat}

In~\ref{sec:hty-mon-cats} we saw how a homotopy monoidal category gave rise
to an ordinary monoidal category, and generalized as follows: if \cat{N} is
any monoidal 2-category and \und{\cat{N}} the associated monoidal category
with equivalences, then a homotopy monoid in \und{\cat{N}} gives rise to a
weak monoid in \cat{N}. (This allowed us to deduce comparison results
involving $A_n$-spaces and $A_n$-algebras.)

Here we show that this process works not just for homotopy monoids but for
homotopy algebras for any operad $P$. Thus if \cat{N} is a monoidal
2-category, there is a concept of `weak $P$-algebra' in \cat{N}, and any
homotopy $P$-algebra in \und{\cat{N}} gives rise to one of these. Most of the
section will in fact be devoted to defining weak algebras. Seasoned category
theorists can probably imagine the kind of definition this is. Once the
definition is made, the actual result is quite easily proved. 

So, let $P$ be a non-symmetric operad and \cat{N} a monoidal 2-category. Let
\und{\cat{N}} denote the underlying monoidal category of \cat{N}. If $A$ and
$B$ are objects of \cat{N} then there is a category \homset{\cat{N}}{A}{B},
whose objects are the 1-cells $A\go B$ in \cat{N} and whose morphisms are
2-cells; there is also a mere set \homset{\und{\cat{N}}}{A}{B}, whose
elements are the 1-cells $A\go B$. A $P$-algebra in \und{\cat{N}} consists of
an object $A$ together with a function $P(n) \go
\homset{\und{\cat{N}}}{A^{\otimes n}}{A}$ for each $n$, satisfying axioms; by
weakening the axioms on these functions we arrive at the following
definition. A \emph{weak $P$-algebra in \cat{N}} consists of
\begin{itemize}
\item an object $A$ of \cat{N}
\item for each $n\in \nat$ and $\theta\in P(n)$, a 1-cell
\[
\ovln{\theta}: A^{\otimes n} \go A
\]
in \cat{N}
\item for each 
\[
\theta\in P(n), \theta_1 \in P(k_1), \ldots, \theta_n \in P(k_n),
\]
an invertible 2-cell
\[
A^{\otimes (k_1 +\cdots+ k_n)}%
\ctwo{\ovln{\theta}\of(\ovln{\theta_1} \otimes\cdots\otimes 
\ovln{\theta_n})}%
{\ovln{\theta\of\bftuple{\theta_1}{\theta_n}}}{\diso}%
A
\]
in $N$, called $\gamma_{\theta; \theta_1, \ldots, \theta_n}$
\item an invertible 2-cell
\[
A%
\ctwo{1_A}{\ovln{1_P}}{\diso}%
A
\]
in \cat{N}, called $\iota$.
\end{itemize}
Then $\gamma$ and $\iota$ are required to satisfy coherence axioms looking
like associativity and identity laws, as detailed in
Figure~\ref{fig:coh-wk-alg}.
\begin{figure}

\begin{diagram}
\ovln{\theta} \of
 (\ovln{\theta_1} \otimes\cdots\otimes \ovln{\theta_n}) \of
 (\ovln{\theta_1^1} \otimes\cdots\otimes \ovln{\theta_n^{k_n}})
&\rEquals
&\ovln{\theta} \of
(\theta_1 \of (\ovln{\theta_1^1} \otimes\cdots\otimes \ovln{\theta_1^{k_1}})
 \otimes\cdots\otimes 
 \theta_n \of (\ovln{\theta_n^1} \otimes\cdots\otimes \ovln{\theta_n^{k_n}}))	\\
\dTo		&	&\dTo						\\
\ovln{\theta\of\bftuple{\theta_1}{\theta_n}} \of 
 (\ovln{\theta_1^1} \otimes\cdots\otimes \ovln{\theta_n^{k_n}})
&	&
\ovln{\theta} \of
 (\ovln{\theta_1 \of \bftuple{\theta_1^1}{\theta_1^{k_1}}}
 \otimes\cdots\otimes
 \ovln{\theta_n \of \bftuple{\theta_n^1}{\theta_n^{k_n}}}			\\
\dTo	&	&\dTo							\\
\ovln{\theta \of 
 \bftuple{\theta_1}{\theta_n} \of
 \bftuple{\theta_1^1}{\theta_n^{k_n}}}
&\rEquals
&\ovln{\theta \of
 \bftuple{\theta_1 \of \bftuple{\theta_1^1}{\theta_1^{k_1}}}%
{\theta_n \of \bftuple{\theta_n^1}{\theta_n^{k_n}}}}				\\
\end{diagram}

\[
\begin{diagram}
\ovln{\theta} \of (1_A \otimes\cdots\otimes 1_A)
&\rEquals	&\ovln{\theta}				\\
\dTo		&		&\dEquals		\\
\ovln{\theta} \of (\ovln{1_P} \otimes\cdots\otimes \ovln{1_P})
&\rTo		&\ovln{\theta\of\bftuple{1_P}{1_P}}	\\
\end{diagram}
\diagspace\diagspace\diagspace
\begin{diagram}
1_A \of \ovln{\theta}		&\rEquals	&\ovln{\theta}	\\
\dTo				&		&\dEquals	\\
\ovln{1_P} \of \ovln{\theta}	&\rTo		&\ovln{1_P \of \theta}\\
\end{diagram}
\]

\caption{Coherence axioms for a weak $P$-algebra: these diagrams must
commute. Here $\theta\in P(n)$, $\theta_i \in P(k_i)$, and all arrows are
part(s) of $\gamma$ or $\iota$.}
\label{fig:coh-wk-alg}
\end{figure}

If $A$ and $B$ are weak $P$-algebras in \cat{N} then a \emph{weak map} from
$A$ to $B$ consists of
\begin{itemize}
\item a 1-cell $f: A\go B$
\item for each $\theta\in P(n)$, an invertible 2-cell
\[
A^{\otimes n}%
\ctwo{\ovln{\theta} \of f^{\otimes n}}{f\of\ovln{\theta}}{\diso}%
B
\]
in \cat{N}, called $\phi_\theta$. 
\end{itemize}
Of course, $\phi$ is required to satisfy coherence axioms, as shown in
Figure~\ref{fig:coh-wk-map}. 
\begin{figure}

\begin{diagram}
\ovln{\theta} \of 
 (\ovln{\theta_1} \otimes\cdots\otimes \ovln{\theta_n}) \of
 f^{\otimes (k_1 +\cdots+ k_n)}				&
\rTo							&
\ovln{\theta\of\bftuple{\theta_1}{\theta_n}} \of 
 f^{\otimes (k_1 +\cdots+ k_n)}				\\
\dEquals 	&	&	\\
\ovln{\theta} \of 
 (\ovln{\theta_1} \of f^{\otimes k_1}
  \otimes\cdots\otimes
  \ovln{\theta_n} \of f^{\otimes k_n})	&	&	\\
\dTo		&	&	\\
\ovln{\theta} \of
 (f\of\ovln{\theta_1} \otimes\cdots\otimes f\of\ovln{\theta_n})&
			&	\\
\dEquals	&	&\dTo	\\
\ovln{\theta} \of f^{\otimes n} \of 
 (\ovln{\theta_1} \otimes\cdots\otimes \ovln{\theta_n})	&
			&	\\
\dTo		&	&	\\
f \of \ovln{\theta} \of (\ovln{\theta_1} \otimes\cdots\otimes \ovln{\theta_n})
&\rTo	
&f \of \ovln{\theta\of\bftuple{\theta_1}{\theta_n}}	\\
\end{diagram}

\begin{diagram}
1_{B} \of f	&\rTo	&\ovln{1_P} \of f	\\
\dEquals	&	&			\\
f		&	&\dTo			\\
\dEquals	&	&			\\
f \of 1_A	&\rTo	&f \of \ovln{1_P}	\\
\end{diagram}

\caption{Coherence axioms for a weak map of weak $P$-algebras. Vertical
arrows come from components of $\phi$, and horizontal arrows from $\gamma$ or
$\iota$.}
\label{fig:coh-wk-map}
\end{figure}
We thus arrive at a category \homset{\fcat{WkAlg}}{P}{\cat{N}} of weak
$P$-algebras in \cat{N}.

(Australian category theorists have used this style of definition
extensively in the study of two-dimensional algebra: see \cite{BKP}, for
instance.)

As a motivating example, consider $P=\Mon$ and $\cat{N}=\Cat$. A weak
\Mon-algebra in \Cat\ consists of a category $A$ and a functor
\[
\otimes_n: A^n \go A
\]
for each $n\geq 0$, which we think of as $n$-fold tensor, together with some
coherence data. Morally this is the same thing as a monoidal category in the
traditional sense, the only difference being that the traditional definition
gives a special role to the values $0$ and $2$ of $n$. In fact, the
categories \homset{\fcat{WkAlg}}{\Mon}{\Cat} and \fcat{MonCat} ($=$ monoidal
categories and monoidal functors) are equivalent. I hope to write this result
up sometime soon; meanwhile, some related considerations can be found in
\cite[4.4]{GECM} and \cite[p.~8]{SHDCT}. More generally,
\homset{\fcat{WkAlg}}{\Mon}{\cat{N}} is equivalent to the category of weak
monoids in \cat{N} (as defined on page~\pageref{page:pseudo-monoid}) for any
monoidal 2-category \cat{N}.

We now have the language in which to state and prove the main result. Just as
in the case of homotopy monoids in \Cat, there are some issues concerning
arbitrary choices, which can be dealt with as they were then; we do not give
them further attention here. 

\begin{thm}	\label{thm:hty-to-wk}
Let $P$ be a non-symmetric operad, let \cat{N} be a monoidal 2-category, and
let \und{\cat{N}} be the associated monoidal category with equivalences (as
in~\ref{sec:env}\bref{eg:env:mon-2-cat}). Then there is a functor
\[
\homset{\HtyAlg}{P}{\und{\cat{N}}} \go 
\homset{\fcat{WkAlg}}{P}{\cat{N}}
\]
sending \pr{X}{\xi} to a weak algebra with underlying object $X(1)$.
\end{thm}
\begin{sketchpf}
Take a homotopy $P$-algebra in \und{\cat{N}},
\[
\pr{X}{\xi}: 
\triple{\fmon{P}}{+}{0} \go \triple{\und{\cat{N}}}{\otimes}{I}.
\]
For each $n, k_1, \ldots, k_n \in\nat$, there is a 1-cell
\[
\xi_{k_1, \ldots, k_n} : 
X(k_1 +\cdots+ k_n) \go
X(k_1) \otimes\cdots\otimes X(k_n),
\]
as in Example~\ref{sec:brief}\bref{eg:hty-alg:Obj}. By
Lemma~\ref{lemma:adjt-equiv}, we can choose a 1-cell $\psi_{k_1, \ldots,
k_n}$ and 2-cells $\eta_{k_1, \ldots, k_n}, \epsln_{k_1, \ldots, k_n}$ so
that
\[
\tuplebts{\psi_{k_1, \ldots, k_n}, \xi_{k_1, \ldots, k_n},
\eta_{k_1, \ldots, k_n}, \epsln_{k_1, \ldots, k_n}}
\]
is an adjoint equivalence. When $k_1=\cdots =k_n$, we write $\psi_{k_1,
\ldots, k_n}$ as $\psi^{(n)}$.

A weak $P$-algebra structure on the object $X(1)$ of \cat{N} can now be
defined as follows: if $\theta\in P(n)$, then \ovln{\theta} is the composite
\[
X(1)^{\otimes n} \goby{\psi^{(n)}} X(n) \goby{X(\theta)} X(1),
\]
and the invertible 2-cells $\gamma$ and $\iota$ are built
up from $\eta_{k_1, \ldots, k_n}$'s and $\epsln_{k_1, \ldots, k_n}$'s. The
process for maps works similarly. 
\done
\end{sketchpf}

Once again all of this ought to be repeatable, \emph{mutatis mutandis}, in
the symmetric case.

\section{Inside \Cat}
\label{sec:inside-Cat}

We have looked at how a homotopy algebra gives rise to a weak algebra in
various different contexts (\ref{sec:hty-mon-cats}, \ref{sec:A-spaces},
\ref{sec:A-algs}, \ref{sec:mon-2-cat}), and it is natural to wonder whether
there is a converse process. This section provides a partial answer: we show
how a weak $P$-algebra in \Cat\ naturally gives rise to a homotopy
$P$-algebra in \Cat, for any operad $P$. 

Before explaining how this works, let me say some things about the current
incompleteness of this line of thought. Firstly, I do not know how to repeat
the construction in any environment other than \Cat\ (e.g.\ in an arbitrary
monoidal 2-category). Secondly, I have only tried to make all the proper
checks in the case of non-symmetric operads $P$; the symmetric case is
largely undone. Thirdly, I have not looked seriously at whether a weak map of
weak algebras gives rise to a map of homotopy algebras (and in this
connection, see the remarks just before
Proposition~\ref{propn:hty-mon-ftr}). Finally, and perhaps most importantly,
I have not investigated the relation between the two processes
\[
\mbox{(weak $P$-algebras)} \oppair{}{} \mbox{(homotopy $P$-algebras)}
\]
described in this and the previous section. Do they, for instance, form an
adjunction, or an equivalence, or a `weak equivalence' of some kind?

The idea behind this section goes as follows. Suppose that $P=\CMon$, so that
our task is to define a function
\[
\mbox{(symmetric monoidal categories)} \go 
\mbox{(homotopy symmetric monoidal categories)}. 
\]
For instance, given the symmetric monoidal category
\triple{\Ab}{\otimes}{\integers}, we want to define a homotopy commutative
monoid 
\[
\pr{X}{\xi}: \triple{\Phi}{+}{0} \go \triple{\Cat}{\times}{1}
\]
with $X(1)\eqv \Ab$. Recalling the introduction to
Section~\ref{sec:hty-mon-cats}, let us define:
\begin{itemize}
\item $X(1) = \Ab$
\item $X(2)$ has objects all quadruples \tuplebts{L_1,L_2,u,M} in which
$L_1,L_2,M$ are abelian groups and $u: L_1 \times L_2 \go M$ is a bilinear
map with the universal property for a tensor product
\item the functor $X(!): X(2) \go X(1)$ (induced by the map $!:2\go 1$ in
$\Phi$) sends \tuplebts{L_1,L_2,u,M} to $M$
\item the equivalence $\xi_{1,1}: X(2) \go X(1) \times X(1)$ sends
\tuplebts{L_1,L_2,u,M} to \pr{L_1}{L_2}.
\end{itemize}
Continuing to build the definition in this way, we would obtain a homotopy
symmetric monoidal category \pr{X}{\xi}. More truthfully, we would obtain one
had I not made certain simplifications in the descriptions of $X(1)$ and
$X(2)$: they are inaccurate in various respects, e.g.\ we have not paid
enough attention to the unit for the tensor; but they serve to convey the
main idea. This explanation comes from Segal's paper~\cite{SegCCT}.

Now let me try to explain the idea for general operads $P$; I will
concentrate on the case of non-symmetric operads. Let $A$ be a weak
$P$-algebra, and attempt to construct a homotopy $P$-algebra \pr{X}{\xi}. With
\Ab\ above, an object of $X(2)$ was effectively an object \pr{L_1}{L_2} of
$\Ab\times\Ab$ together with an isomorphism
\[
j: L_1 \otimes L_2 \goiso M
\]
between the `official' tensor product $L_1 \otimes L_2$ (which is part of the
structure of the monoidal category \triple{\Ab}{\otimes}{\integers}) and the
`putative' tensor product $M$ of $L_1,L_2$ (which is the image of
\tuplebts{L_1,L_2,u,M} under $X(!)$). So an object of $X(2)$ consists of an
object \pr{L_1}{L_2} of $\Ab^2$ together with data specifying `what happens
to \pr{L_1}{L_2} under every operation of \CMon'. Inspired by this, an object
of $X(n)$ will consist of an object \bftuple{a_1}{a_n} of $A^n$, together
with an object of $A$ for every substring \bftuple{a_{d+1}}{a_{m+d}} of
\bftuple{a_1}{a_n} and every $m$-ary operation $\theta\in P(m)$. Call this
object \formal{\theta}{a_{d+1}}{a_{m+d}}; it plays the role of $M$ in the
\Ab\ example. Then there must also be isomorphisms such as 
\begin{eqnarray*}
\ovln{\theta}\bftuple{a_{d+1}}{a_{m+d}}
&\iso &
\formal{\theta}{a_{d+1}}{a_{m+d}},	\\
a_{d+1}
&\iso &
1_P \atuplebts{a_{d+1}},
\end{eqnarray*}
and a choice of such isomorphisms is also included in the data for the object
of $X(n)$. A proper description of $X(n)$ is given in the proof of the
Proposition below. 

(Since we are treating the case of \emph{non}-symmetric operads, it is not in
the spirit of things to take arbitrary subsets of $\{ \range{1}{n} \}$
without regard to order, which is why we restrict to `substrings' $\{
\range{d+1}{m+d} \}$. See also the comments at the end of the section.)

We now come to the main result.
\begin{propn}
Let $P$ be a non-symmetric operad. Then any weak $P$-algebra in \Cat\ gives
rise to a homotopy $P$-algebra in \Cat. More precisely, given a weak
$P$-algebra in \Cat\ with underlying category $A$, there is an associated
homotopy $P$-algebra \pr{X}{\xi} in \Cat\ with $X(1)\eqv A$. 
\end{propn}
\paragraph*{Remark} This process is canonical, in that the construction of
\pr{X}{\xi} from $A$ does not involve arbitrary choices. Contrast the
converse process (Theorem~\ref{thm:hty-to-wk}).
\paragraph*{Proof}
Take a weak $P$-algebra in \Cat, consisting of a category $A$, 1-cells
\ovln{\theta}, and 2-cells $\gamma$ and $\iota$, as in the
definition~(\ref{sec:mon-2-cat}). We construct a homotopy $P$-algebra
$\pr{X}{\xi}: \fmon{P} \go \Cat$.

To do this, we first we need to construct a category $X(n)$ for each $n\geq
0$. An \emph{object} of $X(n)$ is a tuple
\[
\tuplebts{a_1, \ldots, a_n, \Omega, g, i},
\]
where
\begin{itemize}
\item $a_1, \ldots, a_n$ are objects of $A$
\item $\Omega$ is a function assigning an object \homset{\Omega}{\theta}{d}
of $A$ to each $m\in\nat$, $\theta\in P(m)$, and $d\in\nat$ with $m+d\leq
n$; it is convenient to write \homset{\Omega}{\theta}{d} as
\formal{\theta}{a_{d+1}}{a_{m+d}} (which must be interpreted as a purely
formal expression)
\item $g$ is a family of isomorphisms
\begin{eqnarray*}
g_{\theta, \theta_1, \ldots, \theta_n, d}:
\ovln{\theta}
\bftuple{\formal{\theta_1}{a_{d+1}}{a_{k_1+d}}}%
{\formal{\theta_m}{a_{k_1+\cdots+k_{m-1}+d+1}}{a_{k_1+\cdots+k_m+d}}} \\
\ \goiso 
\formal{(\theta\of\bftuple{\theta_1}{\theta_m})}%
{a_{d+1}}{a_{k_1+\cdots+k_m+d}}
\end{eqnarray*}
in $A$, one for each $m, k_1, \ldots, k_m \in\nat$, $\theta\in P(m)$,
$\theta_1 \in P(k_1)$, \ldots, $\theta_m \in P(k_m)$, and $d\in\nat$ with
$k_1 + \cdots + k_m + d \leq n$
\item $i$ is  a family of isomorphisms
\[
i_d: a_d \goiso 1_P\atuplebts{a_d},
\]
one for each $d\in\{1,\ldots,n\}$,
\end{itemize}
and $g$ and $i$ satisfy coherence axioms looking like those in
Figure~\ref{fig:coh-wk-alg} (page~\pageref{fig:coh-wk-alg}), replacing some
of the $\gamma$'s and $\iota$'s in the Figure with $g$'s and $i$'s. A
\emph{morphism} 
\[
f: \tuplebts{a_1, \ldots, a_n, \Omega, g, i}
\go \tuplebts{a'_1, \ldots, a'_n, \Omega', g', i'}
\]
in $X(n)$ consists of
\begin{itemize}
\item maps $f_1: a_1 \go a'_1$, \ldots, $f_n: a_n \go a'_n$ in $A$
\item a map
\[
f_{\theta,d}: \formal{\theta}{a_{d+1}}{a_{m+d}} 
\go \formal{\theta}{a'_{d+1}}{a'_{m+d}}
\]
(that is, $f_{\theta,d}: \homset{\Omega}{\theta}{d} \go
\homset{\Omega'}{\theta}{d}$) for each $m\in\nat$, $\theta\in P(m)$ and
$d\in\nat$ with $m+d\leq n$,
\end{itemize}
such that the $f_{\theta, d}$'s commute with the $g$'s and the $i$'s. With
the obvious composition and identities, $X(n)$ forms a category. 

Now that $X(n)$ is defined, there is only one sensible way to define the rest
of the data for \pr{X}{\xi}, and I will just sketch it.

To define $X$ on morphisms, take a map $\Psi: n\go p$ in \fmon{P}, which
consists of an expression $n=n_1 + \cdots + n_p$ together with elements
\[
\psi_1\in P(n_1), \ldots, \psi_p \in P(n_p).
\]
Take an object \tuplebts{a_1,\ldots, a_n, \Omega, g, i} of $X(n)$. Then there
is a resulting object
\[
(X\Psi)\tuplebts{a_1,\ldots, a_n, \Omega, g, i}
=
\tuplebts{b_1, \ldots, b_p, \twid{\Omega}, \twid{g}, \twid{i}}  
\]
of $X(p)$, in which
\[
\bftuple{b_1}{b_p}
=
\bftuple{\formal{\psi_1}{a_1}{a_{n_1}}}%
{\formal{\psi_p}{a_{n_1+\cdots+n_{p-1}+1}}{a_{n_1+\cdots+n_p}}}.
\]

The data for $\xi$ consists of a pair of maps 
\[
X(n) \ogby{\xi_{n,n'}^1} X(n+n') \goby{\xi_{n,n'}^2} X(n')
\]
for each $n,n'\in\nat$. The image under $\xi_{n,n'}^1$ of an object 
\[
\tuplebts{a_1, \ldots, a_n, a_{n+1}, \ldots, a_{n+n'}, \Omega, g, i}
\]
of $X(n+n')$ is of the form \tuplebts{a_1, \ldots, a_n, ?, ?,
?}, and the image under $\xi_{n,n'}^2$ is of the form \tuplebts{a_{n+1},
\ldots, a_{n+n'}, ?, ?, ?}. 

Once all the details are filled in, we arrive at a colax monoidal functor
\[
\pr{X}{\xi}: \fmon{P} \go \Cat.
\]
The remaining tasks are to show that this is in fact a homotopy
$P$-algebra---that is, the maps $\xi_{n,n'}$ and $\xi_0$ are
equivalences---and that $X(1)\eqv A$. We do both these things at once by
considering the forgetful functors
\[
\begin{array}{rccc}
U_n:	&X(n)	&\go	&A^n	\\
	&\tuplebts{a_1, \ldots, a_n, \Omega, g, i}
		&\goesto
			&\bftuple{a_1}{a_n}.
\end{array}
\]
Evidently, the squares
\[
\begin{diagram}
X(n+n')		&\rTo^{\xi_{n,n'}}	&X(n)\times X(n')		\\
\dTo<{U_{n+n'}}	&			&\dTo>{U_n \times U_{n'}}	\\
A^{n+n'}	&\rTo_{\diso}		&A^n \times A^{n'}		\\
\end{diagram}
\diagspace
\begin{diagram}
X(0)		&\rTo^{\xi_0}	&\One		\\
\dTo<{U_0}	&		&\dEquals	\\
A^{0}		&\rTo_{\diso}	&\One		\\
\end{diagram}
\]
both commute, so we will be finished if we can show that each functor $U_n$
is an equivalence. And this in turn is easily accomplished: for $U_n$ has a
canonical pseudo-inverse, sending \bftuple{a_1}{a_n} to an object
\tuplebts{a_1, \ldots, a_n, \Omega, g, i} of $X(n)$. Here 
\[
\homset{\Omega}{\theta}{d} = \ovln{\theta}\bftuple{a_{d+1}}{a_{m+d}},
\]
or put another way,
\[
\formal{\theta}{a_{d+1}}{a_{m+d}} = \ovln{\theta}\bftuple{a_{d+1}}{a_{m+d}},
\]
and $g$ and $i$ are respectively defined by $\gamma$ and $\iota$. 
\done
\paragraph*{}
I believe that the Proposition can be repeated for the case of
symmetric operads. To do this one would replace substrings \bftuple{d+1}{m+d}
of \bftuple{1}{n} with arbitrary non-repeating sequences \bftuple{d_1}{d_m},
with $d_i \in \{ \range{1}{n} \}$; but I have not verified this yet.

\chapter{Change of Environment}
\label{ch:change}

Throughout this work we have discussed homotopy $P$-algebras in \emm\ for a
\emph{fixed} operad $P$ and a \emph{fixed} monoidal category \emm\ with
equivalences. In this short chapter we look at what happens when the
`environment' \emm\ is varied. In other words, we look at how a suitable map
$\cat{L} \go \emm$ induces a functor
\[
\homset{\HtyAlg}{P}{\cat{L}} \go \homset{\HtyAlg}{P}{\emm},
\]
for any operad $P$.

It is not difficult to say precisely how this process works
(see~\ref{sec:principle} below), but first let me try to explain why such a
result is plausible. 

From a formal point of view, a homotopy $P$-algebra in \cat{L} is a certain
kind of map $\fmon{P} \go \cat{L}$, so composing with the right kind of map
$\cat{L} \go \emm$ ought to yield a homotopy $P$-algebra in \emm.

From another point of view, consider, for instance, the path-components
functor 
\[
\pi_0: \Top \go \Set.
\]
Since $\pi_0$ preserves products, and `group objects' can be formed in any
category in which products exist, $\pi_0$ induces a functor
\[
\mbox{(topological groups)} \go \mbox{(groups)}.
\]
More generally, if we have fixed an algebraic theory (groups, in this case)
then any functor $\cat{C} \go \cat{D}$ of the right kind will induce a
functor 
\[
\mbox{(algebras in \cee)} \go \mbox{(algebras in \cat{D})}.
\]
Since a homotopy $P$-algebra is some kind of algebraic structure (albeit a
rather loose kind), we might expect the same principle to apply; it does.

From a third point of view, topologists will expect results such as `the
classifying space of a monoidal category is a homotopy monoid'. This is
indeed the case in our theory, as long as we read `monoidal category' as
`homotopy monoidal category': see
Example~\ref{sec:change-egs}\bref{eg:change:classifying}. 

We could also change operad: a map $P\go Q$ of operads induces a functor
\[
\homset{\HtyAlg}{P}{\emm} \og \homset{\HtyAlg}{Q}{\emm}.
\]
But this will not be discussed here.

Section~\ref{sec:principle} sets out exactly how a `change of environment'
induces a functor between categories of homotopy algebras, and
Section~\ref{sec:change-egs} lists some examples.

\section{The Principle}
\label{sec:principle}

\begin{defn}
\begin{enumerate}
\item	\label{defn-part:non-sym-HMF}
Let \cat{L} and \emm\ be monoidal categories with equivalences. A
\emph{homotopy monoidal functor} $\cat{L}\go\emm$ is a colax monoidal
functor \pr{F}{\phi} such that
\begin{itemize}
\item each component $\phi_0, \phi_{m,n}$ of $\phi$ is an equivalence in
\emm 
\item if $f$ is an equivalence in \cat{L} then $F(f)$ is an equivalence in
\emm. 
\end{itemize}
\item
\emph{Homotopy symmetric monoidal functors} are defined by changing
`monoidal' to `symmetric monoidal' throughout
part~\bref{defn-part:non-sym-HMF}.
\end{enumerate}
\end{defn}
Note that if $P$ is an operad and the (symmetric) monoidal category \fmon{P}
is equipped with isomorphisms as its equivalences, then a homotopy
$P$-algebra in \cat{L} is exactly a homotopy (symmetric) monoidal functor
$\fmon{P} \go \cat{L}$. Note also that the composite of two homotopy
(symmetric) monoidal functors is a homotopy (symmetric) monoidal
functor. Hence a homotopy (symmetric) monoidal functor $\cat{L}\go\emm$
induces a functor
\[
\homset{\HtyAlg}{P}{\cat{L}} \go \homset{\HtyAlg}{P}{\emm}.
\]
This simple piece of theory is the basis of this chapter, the remainder of
which consists of examples.

\section{Examples}
\label{sec:change-egs}

\begin{enumerate}

\item
Suppose that \emm\ is a monoidal category and that \Eee\ and \cat{E'} are
both classes of equivalences in \emm, with $\cat{E'}\sub\Eee$. Then the
identity is a homotopy monoidal functor $\pr{\emm}{\cat{E'}} \go
\pr{\emm}{\Eee}$. Thus if $P$ is any operad then there is an induced functor
\[
\homset{\HtyAlg}{P}{\pr{\emm}{\cat{E'}}} \go 
\homset{\HtyAlg}{P}{\pr{\emm}{\Eee}}
\]
(with what I hope is self-explanatory notation), and this is the obvious
inclusion. In particular, if $\cat{E'} = \{\textrm{isomorphisms}\}$ then this
is the inclusion
\[
\homset{\Alg}{P}{\emm} \go
\homset{\HtyAlg}{P}{\pr{\emm}{\Eee}}.
\]

\item	\label{eg:change:Toph}
Let \fcat{Toph} be the category whose objects are topological spaces and
whose morphisms are homotopy classes of continuous maps. Let $P$ be an
operad. The weakest possible meaning of the phrase `homotopy topological
$P$-algebra' is `$P$-algebra in \fcat{Toph}': e.g.\ a `homotopy topological
semigroup' in this weakest sense is just a space $A$ with a binary operation
which is associative up to homotopy. Any homotopy $P$-algebra in the sense of
this paper certainly gives rise to one of these very weak
structures. Formally, let $Q: \Top \go \fcat{Toph}$ be the quotient functor
(which is the identity on objects). Equip \fcat{Toph} with just the
isomorphisms as its equivalences. Then $Q$ becomes a homotopy symmetric
monoidal functor
\[
\triple{\Top}{\times}{1} \go \triple{\fcat{Toph}}{\times}{1},
\]
so for any operad $P$ there is an induced functor 
\[
\homset{\HtyAlg}{P}{\Top} \go 
\homset{\HtyAlg}{P}{\fcat{Toph}} \eqv \homset{\Alg}{P}{\fcat{Toph}}.
\]

\item
Example~\bref{eg:change:Toph} can be repeated with chain complexes in place
of spaces, or with categories in place of spaces (with natural isomorphism
classes of functors), or indeed with the objects of any monoidal 2-category.

\item	\label{eg:change:classifying}
Let $B: \Cat\go\Top$ be the classifying-space functor (see
\cite{SegCSS}). Then $B$ preserves products and sends equivalences to
homotopy equivalences, so there is an induced functor
\[
B: \homset{\HtyAlg}{P}{\Cat} \go \homset{\HtyAlg}{P}{\Top}
\]
for any operad $P$. (Here \Top\ is equipped with the cartesian monoidal
structure.) For instance, let $C$ be a homotopy monoidal category, i.e.\ a
homotopy monoid in \Cat: then $BC$ is a homotopy topological
monoid. Similarly, the classifying space $BC$ of a homotopy symmetric
monoidal category is a homotopy topological commutative monoid. This
symmetric version is exactly Segal's observation (in \cite[\S 2]{SegCCT})
that the classifying space of a (special) $\Gamma$-category is a (special)
$\Gamma$-space.

\item	\label{eg:change:hom}
If $B$ is a fixed space with basepoint then 
\[
\homset{\Topstar}{\dashbk}{B}:
\triple{\Top_*^{\op}}{\wej}{1} \go \triple{\Top}{\times}{1}
\]
is a homotopy monoidal functor, as observed in~\ref{sec:loops}. So there is
in particular an induced functor
\[
\homset{\HtyAlg}{\Mon}{\triple{\Top_*^{\op}}{\wej}{1}}
\go
\homset{\HtyAlg}{\Mon}{\triple{\Top}{\times}{1}}.
\]
This is effectively the argument we used in~\ref{sec:loops} to show that the
homotopy comonoid structure on $S^1$ gave a homotopy monoid structure on the
loop space \homset{\Topstar}{S^1}{B}. 

\item
In Section~\ref{sec:iterated} we used the $n$-fold smash product
\[
\smsh: \triple{\Topstar}{\wej}{1}^n
\go \triple{\Topstar}{\wej}{1},
\]
which is a `homotopy multi-monoidal functor' in the obvious sense of the
phrase. It translates the homotopy comonoid structure on $S^1$ (or rather,
$n$ copies of this structure) into an $n$-fold homotopy comonoid structure
on $S^1 \smsh\cdots\smsh S^1 = S^n$. 

\item	\label{eg:change:path-components}
The path-components functor $\pi_0: \Top \go \Set$ preserves products and
sends homotopy equivalences to isomorphisms, and so induces a functor
\[
\homset{\HtyAlg}{P}{\Top} \go \homset{\Alg}{P}{\Set}.
\]
For instance, the path-components of any topological monoid form a monoid.

More generally, if $P_1, \ldots, P_n$ are operads then $\pi_0$ induces a
functor
\[
\pi_0: \multihom{\HtyAlg}{P_1, \ldots, P_n}{\Top} 
\go 
\multihom{\Alg}{P_1, \ldots, P_n}{\Set}
\]
(see~\ref{sec:iterated} for the notation). We saw in~\ref{thm:n-fold-monoid}
that any $n$-fold loop space \homset{\Topstar}{S^n}{B} has the structure of
an $n$-fold homotopy monoid: thus 
\[
\pi_n(B) = \pi_0(\homset{\Topstar}{S^n}{B})
\]
is an $n$-fold monoid in \Set. But the Eckmann-Hilton
argument\footnote{$a\cdot b = (a*1)\cdot (1*b) = (a\cdot 1)*(1\cdot b) = a*b
= (1\cdot a)*(b\cdot 1) = (1*b)\cdot (a*1) = b\cdot a$}
says that if a pair of monoid structures on a set commute with each other
then they are identical and commutative: so
\[
\multihom{\Alg}{\underbrace{\Mon,\ldots,\Mon}_n}{\Set}
=
\left\{
\begin{array}{ll}
\Set				&\textrm{if } n=0	\\
\textrm{(monoids)}		&\textrm{if } n=1	\\
\textrm{(commutative monoids)}	&\textrm{if } n\geq 2.	\\
\end{array}
\right.
\]
This means that the $n$th homotopy $\pi_n(B)$ of a space $B$ is a set when
$n=0$, a monoid when $n=1$, and a commutative monoid when $n\geq 2$. Of
course, we know that these monoids are actually groups, and so our result
implies that the higher homotopy groups are abelian. 

\item
Let $\Pi_1: \Top \go \Cat$ be the functor assigning to a space its
fundamental groupoid. Then $\Pi_1$ is a homotopy monoidal functor
\[
\triple{\Top}{\times}{1} \go \triple{\Cat}{\times}{1}
\]
by virtue of preserving products: so, for example, the fundamental groupoid
of a loop space is a homotopy monoidal category. Similarly, the fundamental
groupoid of a special $\Gamma$-space is a special $\Gamma$-category. The
fundamental groupoid of an $n$-fold loop space is an $n$-fold homotopy
monoidal category, i.e.\ a homotopy \bftuple{\Mon}{\Mon}-algebra in \Cat. I
have not investigated $n$-fold homotopy monoidal categories, but it would be
interesting to see how they compare to braided monoidal categories when
$n=2$, and to the iterated monoidal categories of Balteanu, Fiedorowicz,
Schw\"{a}nzl and Vogt for general $n$ (see \cite{JS} and \cite{BFSV}
respectively).

\item
Next is a non-example. One of the main features of the definition of
$A_\infty$-algebra in Stasheff's original paper \cite{HAHII} is that the
chain complex $C_{\blob}(A)$ of an $A_\infty$-space $A$ is an
$A_\infty$-algebra. We can attempt to mimic this here, by trying to give the
singular chains functor $C_{\blob}: \Top\go\ChCx$ the structure of a homotopy
monoidal functor. If this is possible then there is an induced functor
\[
\homset{\HtyAlg}{\Mon}{\Top} \go \homset{\HtyAlg}{\Mon}{\ChCx},
\]
so that the chains of a homotopy topological monoid form a homotopy d.g.\
algebra. However, it appears to be impossible.

To get an idea of the issues at hand, let us see how $C_{\blob}$ is naturally
a \emph{lax} monoidal functor (as defined on
page~\pageref{page:defn-lax}). This basically means that for spaces $X$ and
$Y$ and $p,q\in\nat$ there is a canonical map
\[
C_p(X) \otimes C_q(Y) \go C_{p+q}(X\times Y),
\]
that is,
\[
R\langle \homset{\Top}{\Delta^p}{X} \times \homset{\Top}{\Delta^q}{Y}\rangle
\go
R\langle \homset{\Top}{\Delta^{p+q}}{X\times Y} \rangle
\]
where $R$ is the ground ring, $R\langle S\rangle$ is the free $R$-module on a
set $S$, and $\Delta^r$ is the standard $r$-simplex. This map is induced by
the composite
\[
\homset{\Top}{\Delta^p}{X} \times \homset{\Top}{\Delta^q}{Y}
\goby{\times}
\homset{\Top}{\Delta^p \times \Delta^q}{X\times Y}
\goby{f^*}
\homset{\Top}{\Delta^{p+q}}{X\times Y}
\]
where, in turn, $f: \Delta^{p+q} \go \Delta^p \times \Delta^q$ is defined by
the first degeneracy map $\Delta^{p+q} \go \Delta^p$ and the last degeneracy
map $\Delta^{p+q} \go \Delta^q$. 

(If this is right then it contradicts the suggestion of Kontsevich,
in~\cite[2.2]{Kon}, that one needs to use cubical rather than simplicial
chains in order to make $C_{\blob}$ into a lax monoidal functor.)

Hence $C_{\blob}$ naturally induces a functor
\[
\homset{\Alg}{P}{\Top} \go \homset{\Alg}{P}{\ChCx}
\]
for any operad $P$. So, for instance, the chains of a \emph{genuine}
topological monoid form a \emph{genuine} d.g.\ algebra. But this has come
from $C_{\blob}$ being a \emph{lax} monoidal functor, and what we need in
order to obtain an induced functor on homotopy algebras is its dual, a
\emph{colax} monoidal functor. As far as I know, there is no suitable colax
structure. 

\item
The homology functor $H_{\blob}: \ChCx_R \go \GrMod_R$ sends chain homotopy
equivalences to isomorphisms, for any commutative ring $R$. It is also a
monoidal functor, i.e.\ preserves $\otimes$ and unit up to coherent
isomorphism, provided that $R$ is a field (by the K\"{u}nneth Theorem,
\cite[3.6.3]{Wei}). So if $R$ is a field then there is an induced functor
\[
\homset{\HtyAlg}{P}{\ChCx_R} \go \homset{\Alg}{P}{\GrMod_R}
\]
for any operad $P$. Hence the homology over a field of a homotopy d.g.\
algebra forms a graded algebra, and similarly for commutative algebras,
non-unital algebras, Lie algebras, Gerstenhaber algebras, etc.

\end{enumerate}

\chapter{Final Thoughts}
\label{ch:thoughts}

This has been a long paper, and despite having covered many points, there are
still various loose ends and unanswered questions. I hope it will therefore
be useful for me to give a summary of how things stand.

First is a list of things done, and then things conspicuously undone. Staying
negative in tone, there is next a section on homotopy invariance. More
optimistically, the view is then put forward that our definition of homotopy
algebra is just a 1-dimensional approximation to an infinite-dimensional
ideal, and that the distance between approximation and ideal is what causes
many of our difficulties. Also discussed, briefly, is the matter of how our
definition relates to other definitions of homotopy algebra.

\section*{What We've Done, and What We Haven't}

The main achievements of this paper are as follows.
\begin{description}

\item[General definition] The principal point of the paper is, of course, to
give a definition of homotopy algebra for an operad which works in a very
general context. We have done this, and once one has understood the process
of forming the free monoidal category \fmon{P} on an operad $P$, the
definition is extremely simple.

\item[Special $\Gamma$-spaces and $\Delta$-spaces] We have shown that a
homotopy topological monoid is precisely a `special $\Delta$-space' (or
special simplicial space, in our terminology), and similarly that a homotopy
topological commutative monoid is precisely a special $\Gamma$-space. Indeed,
it was a reformulation of the definition of special $\Gamma$-space which
led me to the general definition. The advantage of this reformulation is that
it allows generalization: to an arbitrary operad $P$ (not just \Mon\ or
\CMon), and to monoidal categories which, unlike \Top, are not cartesian. The
reformulation also clarifies the role of $\Gamma$ from a conceptual point of
view, and clarifies the interplay of $\Delta$ and $\Delta^+$.

\item[Loop spaces] A major example of our definition is that any loop space
is a homotopy topological monoid, and, in fact, that any $n$-fold loop space
is an $n$-fold homotopy topological monoid. To express the latter statement
we had to develop (in brief) a theory of homotopy algebras for several
operads simultaneously. Completing the picture is the fact that any infinite
loop space is a special $\Gamma$-space, i.e.\ a homotopy topological
commutative monoid.

\item[Change of environment] A suitable map $\cat{L}\go\emm$ of monoidal
categories leads, unsurprisingly, to a way of passing from homotopy algebras
in \cat{L} to homotopy algebras in \emm. The most interesting applications
presented here are topological: the classifying space of a homotopy
(symmetric) monoidal category is a homotopy topological (commutative) monoid,
and the fundamental groupoid of an $n$-fold loop space is an $n$-fold
homotopy monoidal category. It also provides an explanation of why the higher
homotopy groups of a space are abelian, and why the first homotopy group is
not, and why the zeroth is only a set.

\item[Comparisons] There are various other notions of weakened or
up-to-homotopy algebraic structure in the literature. We have made some
partial comparisons between a small number of these and our definition. The
result which encompasses most of our comparisons is that a homotopy
$P$-algebra in a monoidal 2-category gives rise to a weak
$P$-algebra~(\ref{sec:mon-2-cat}). It follows that a homotopy monoidal
category gives rise to a monoidal category (non-strict, in the traditional
sense), a homotopy semigroup in the category of based spaces gives rise to an
$A_4$-space, and a homotopy differential graded non-unital algebra to an
$A_4$-algebra. In the opposite direction, we have also shown how to obtain a
homotopy monoidal category from a (traditional) monoidal category.

\item[Clarification] Aside from the specific points listed above, I hope that
this paper has succeeded in clarifying some general points concerning
homotopy algebras. 

Firstly, we have seen that in order to state our definition of homotopy
algebra in \emm, it is only necessary to have knowledge of which morphisms in
\emm\ are `homotopy equivalences'; knowledge of what it means for two maps to
be `homotopic', or what a `homotopy' between maps is, etc., is not
required. 

Secondly, I have tried to draw attention to the distinction between the
canonical and the non-canonical, especially in the sections on loop spaces
and homotopy monoidal categories~(\ref{sec:loops}
and~\ref{sec:hty-mon-cats}). For instance, since there is no canonical recipe
for forming the tensor product of two abelian groups, there is no canonical
functor $\otimes: \Ab^2 \go \Ab$; similarly, there is no canonical way of
composing two based loops in a space. In this connection we repeatedly see
diagrams of the shape
\begin{diagram}[size=1.5em]
	&	&Z	&	&	\\
	&\ldTo	&	&\rdTo	&	\\
X	&	&	&	&Y	\\
\end{diagram}
in which the left-hand map is an equivalence, as a substitute for a map $X\go
Y$. Such diagrams appear in many parts of mathematics: to take a fairly
random selection, \cite{Thomas}, \cite{Mak}, \cite[p.~51]{Ad}, \cite{MayGIL}.
\end{description}

This is what we've done. On the other hand, the ideas presented in this paper
raise many questions crying out to be answered. I believe that the central
definition of homotopy algebra is fundamentally crude and cannot be
formulated satisfactorily until there is a decent theory of weak
$\infty$-categories; but that will be discussed later, and for now I will
limit myself to noting some specific shortcomings.
\begin{description}

\item[Comparisons] The comparison results presented here are blatantly
incomplete. For a start, we showed that homotopy differential graded
non-unital algebras give rise to $A_4$-algebras, but were not able to show
that they give $A_\infty$-algebras; and similarly $A_\infty$-spaces. More
seriously, the comparison results are almost all of the form `a homotopy
algebra in our sense gives rise to a homotopy algebra in someone else's
sense', rarely the other way round. The exception is when we are taking
algebras in \Cat~(\ref{sec:inside-Cat}); but even then, it is not clear
whether the two processes are in any sense mutually inverse or adjoint.

\item[Maps] We have discussed homotopy algebras at length, but not homotopy
maps between homotopy algebras. Just before
Proposition~\ref{propn:hty-mon-ftr} we suggested how a homotopy map between
homotopy algebras might look in a category such as \Top\ where there is a
notion of two maps being homotopic. Another possibility, which makes sense in
any monoidal category with equivalences, is to define a homotopy map of
homotopy $P$-algebras as a homotopy \Map{P}-algebra. Here \Map{P} is the
multicategory (coloured operad) of~\ref{sec:multicats}, for which a genuine
algebra is a pair of $P$-algebras with a map between them. This is a pleasing
definition of homotopy map, but raises further questions when one thinks
about composing them.

\item[Examples] We are a little short on actual examples of homotopy
algebras, mostly for the reasons mentioned under `Comparisons' just
above. Even if one has in mind an object which one suspects ought to be a
homotopy algebra for a certain operad $P$, it takes creative effort to endow
the object with the structure of a homotopy $P$-algebra. In the terminology
of~\cite[p.~60]{Ad}, one has to create a lot of flab. An example of this is
the Problem of~\ref{subsec:Ger}: how to endow the Hochschild cochain complex
of an associative algebra with the structure of a homotopy d.g.\ commutative
algebra.
\end{description}

This completes the summary of things done and undone.

\section*{Homotopy Invariance}

Earlier we came across a disturbing feature of the definition of homotopy
algebra, in Example~\ref{sec:brief}\bref{eg:hty-alg:Act-G}. This was that if
$G$ is a monoid and \Act{G} the non-symmetric operad whose algebras are
objects with a strict action by $G$, then a homotopy \Act{G}-algebra
\pr{X}{\xi} specifies, amongst other things, a \emph{strict} action of $G$ on
the `base object' $X(1)$.

From this example we can see that our homotopy $P$-algebras are not
`homotopy invariant algebraic structures'; at least, not if we regard a
homotopy $P$-algebra \pr{X}{\xi} as being a structure on the
object $X(1)$. That is to say, suppose that \pr{X}{\xi} is a homotopy
algebra for some operad $P$, in some (symmetric or not) monoidal category
\emm\ with equivalences, and suppose that we have a homotopy equivalence
$X(1)\go A$, or $A\go X(1)$, in \emm. `Homotopy invariance' says that there
is an induced homotopy $P$-algebra \pr{W}{\omega} with $W(1)=A$. In the case
$P=\Act{G}$ and $\emm=\Top$, this implies that if $A$ is homotopy equivalent
to a (strict) $G$-space then there is an induced \emph{strict} action of $G$
on $A$. Plainly that is not true, so homotopy invariance fails.

This is worrying: homotopy invariance is an attribute which a good theory of
homotopy-algebraic structures ought to have. (Boardman and
Vogt's book \cite{BV} and Markl's paper \cite{MarHAHA} say much more on why
it is desirable.) In the next section, I will suggest in vague terms an 
$\infty$-categorical version of our definition of homotopy algebra which
\emph{would} be homotopy invariant; but for now, here is a result which is
perhaps the closest we have to homotopy invariance for the definition as it
stands. 
\begin{propn}
Let $P$ be a (symmetric) operad and \emm\ a (symmetric) monoidal
category. Let
\[
\pr{W}{\omega}, \pr{X}{\xi}: \fmon{P} \go \emm
\]
be two colax (symmetric) monoidal functors, and let
\[
\sigma: \pr{W}{\omega} \go \pr{X}{\xi}
\]
be a monoidal transformation which is a `homotopy equivalence', in the sense
that each component $\sigma_n$ is a homotopy equivalence in \emm. Then
\pr{W}{\omega} is a homotopy $P$-algebra if and only if \pr{X}{\xi} is.
\end{propn}
\begin{proof}
Simply apply the axioms for a class of equivalences~(\ref{defn:equivs}) to
the commuting diagrams in the definition of monoidal
transformation~(\ref{defn:mon-transf}).  \done
\end{proof}

\section*{$\infty$-Categories}

A monoidal category with equivalences is a very simple device, and it has
already been argued that it is really too simple (in the introduction to
Chapter~\ref{ch:defn}). \Top\ and \ChCx, for example, naturally form monoidal
$\infty$-categories, but when we treated them as monoidal categories with
equivalences we threw away all information about cells of dimension 2 and
above, except for retaining knowledge of which 1-cells are equivalences in an
$\infty$-categorical sense. In Chapter~\ref{ch:hty-mons} we gave them
marginally more respect by treating them as monoidal 2-categories, but this
is still a long way from appreciating their true $\infty$-categorical
nature. This ignorant behaviour is, of course, excused by the fact that there
is not yet a well-developed theory of (weak) $\infty$-categories.

When such a theory has evolved, it should be possible to make the following
definition. Let $P$ be a symmetric or non-symmetric operad with, for
simplicity, each $P(n)$ being just a set. Let \emm\ be a (symmetric) monoidal
$\infty$-category. Then a `homotopy $P$-algebra in \emm' is simply a
(symmetric) monoidal $\infty$-functor $\fmon{P} \go \emm$. Here the monoidal
category \fmon{P} is made into a monoidal $\infty$-category by saying that
the only $k$-cells for $k\geq 2$ are the identities, and a `monoidal
$\infty$-functor' is meant to preserve tensor, composition and identities up
to equivalence in the weakest $\infty$-categorical sense.

In down-to-earth terms, this higher-dimensional structure would make
differences of the following kind. Take, for example, a monoid $G$. With the
definition of homotopy algebra used in this paper, the base object $X(1)$ of
a homotopy \Act{G}-algebra \pr{X}{\xi} is in fact a strict $G$-object, which
means that there is a map $\alpha_g: X(1)\go X(1)$ for each $g\in G$, such
that the diagrams
\[
\begin{diagram}[height=2em]
X(1)	&\rTo^{\alpha_{g'}}	&X(1)			\\
	&\rdTo<{\alpha_{gg'}}	&\dTo>{\alpha_g}	\\
	&			&X(1)			\\
\end{diagram}
\diagspace
\begin{diagram}
X(1)	&\pile{\rTo^{1}\\ \rTo_{\alpha_1}}	&X(1)	\\
\end{diagram}
\]
commute (strictly). It is not possible that they might only have to commute
`up to homotopy': this simply does not make sense in an arbitrary monoidal
category with equivalences. But with an $\infty$-categorical definition of
homotopy algebra they would not strictly commute: instead, there would be a
2-cell filling in each of these diagrams (e.g.\ a homotopy, if we were
working in \Top). Moreover, these 2-cells would obey coherence laws---not
strictly, but up to a 3-cell; and so on. Thus the effect of an
$\infty$-categorical definition would be to weaken further the original
definition of homotopy algebra.

Another place where this weakening effect would be seen is in maps of
homotopy algebras: such a thing would naturally be defined as a weak monoidal
transformation, with `weak' meant in an $\infty$-categorical sense, and this
would be a respectable notion of a homotopy map of homotopy algebras. Thus
the squares involved in the discussion of maps just before
Proposition~\ref{propn:hty-mon-ftr} would commute only up to coherent
equivalence. Similarly, homotopy invariance should work perfectly well when
we use the $\infty$-categorical definition.

So, we have been using throughout a 1-dimensional approximation to an
infinite-dimensional ideal. This is just about the roughest approximation
possible, but still it has given us plenty to chew on. In particular,
homotopy algebras as we defined them are weak enough to include loop spaces
and monoidal categories, and gain a badge of historical respectability by
having as particular cases special $\Gamma$-spaces and special simplicial
spaces.

A similar situation where one `should' use an $\infty$-category, but instead
substitutes a more crude structure, is described in Hinich's paper
\cite[1.3]{Hin}.

\section*{Other Definitions of Homotopy Algebra}

Various other notions of homotopy algebra exist in the literature, and I have
not attempted anything like a systematic comparison. The notions I am aware
of are the homotopy-invariant algebraic structures in Boardman and Vogt's
book~\cite{BV}, the strong homotopy algebras of~\cite{Lada}, the minimal
models of Markl (\cite{MarHARO}, \cite{MarHAHA}), and the $A_\infty$,
$B_\infty$, $C_\infty$, $E_\infty$, $G_\infty$ and $L_\infty$ structures
developed by many people: see, for instance, \cite{Ad}, \cite{KSV},
\cite{LM}, \cite{HAHI}, \cite{HAHII}, \cite{Vor}. (My knowledge of this
literature is not very thorough, so I hope that no-one will be offended by
omissions.) The definition in this paper has the advantage of working in a
more general context than any of the others, as far as I know; but of course
it has certain disadvantages too.

It seems to me that the comparison between different definitions of homotopy
algebra is very much like the comparison between different definitions of
weak $n$-category, of which there are currently about a dozen. In both
situations there are essentially two approaches. We discussed in
Section~\ref{sec:loops} how the description of a loop space as an
$A_\infty$-space is essentially different from the description of a loop
space as a homotopy monoid in our sense: to describe it as an
$A_\infty$-space we have to make an arbitrary choice concerning how to
compose two loops, and so have a specific but non-canonical composition law;
to describe it as a homotopy monoid we need make no artificial choices at
all, but there is no actual preferred composition law. Similarly, the various
proposed definitions of weak $n$-category split into those which insist that
for suitable cells $f$ and $g$ there should be a definite composite $h=g\of
f$, and those in which one could only ever say `$h$ is \textbf{a} composite
of $g$ and $f$', there being perhaps many possible composites of $g$ and $f$,
all equally valid. We have already seen this distinction for monoidal
categories, at the beginning of~\ref{sec:hty-mon-cats}, and indeed a monoidal
category is precisely a weak 2-category with only one object. Examples of the
first (`algebraic') approach are the definitions of weak $n$-category
proposed by Batanin (\cite{Bat}, \cite{StrRMB}), by Penon~\cite{Pen}, and by
me~\cite{SHDCT}, and the `classical' definitions of monoidal category,
bicategory~\cite{Ben} and tricategory~\cite{GPS}. Examples of the second
approach are the definitions of Baez and Dolan~\cite{BaDoHDA3}, Hermida,
Makkai and Power~\cite{HMP}, Joyal~\cite{Joy}, Street~\cite{StrAOS}, and
Tamsamani~\cite{Tam}. To date no-one knows very much about how the various
definitions relate to one another. It may be that some clues lie in the body
of material on `uniqueness of delooping machines', e.g.\ \cite{MT},
\cite{Thmsn} and \cite{SV}; it may also be that the situation for
$n$-categories is substantially more challenging.

\chapter*{Glossary}
\ucontents{chapter}{Glossary}

A brief description of each term is given; for operads, the description says
what the algebras are. Section or page numbers refer to where the term was
defined or first mentioned. Some terms have meanings at two or more different
levels, in which case there is more than one section number. 

\subsection*{Categories}

\begin{tabbing}
$\GrMod=\GrMod_R$	\=spaces and homotopy classes of maps	\=	\kill
$\cee^\op$	\>opposite category
				\>p.~\pageref{page:misc-notation}	\\
\fmon{P}\>free monoidal category on operad	\>\ref{sec:free-mon}	\\
$|\cat{N}|$	\>underlying monoidal category with	\>		\\
 		\>\ equivalences		\>\ref{sec:mon-2-cat}	\\
\One	\>terminal category	
			\>\ref{sec:prelims-mon-cats}\bref{eg:mon-cat:Cat}\\
\Ab	\>abelian groups		
			\>\ref{sec:enr-ops}\bref{eg:enr-op:Lie}	\\
\homset{\Alg}{P}{\emm}	\>algebras
				\>\ref{defn:algebra}, \ref{sec:enr-ops}	\\
\Cat	\>categories	\>\ref{sec:prelims-mon-cats}\bref{eg:mon-cat:Cat},
	\ref{sec:enr-cats}\bref{eg:enr-cat:Cat},
	\ref{sec:env}\bref{eg:env:mon-2-cat}				\\
$\ChCx=\ChCx_R$	\>chain complexes
			\>\ref{sec:prelims-mon-cats}\bref{eg:mon-cat:ChCx},
			\ref{sec:enr-cats}\bref{eg:enr-cat:ChCx},
			\ref{sec:A-algs}				\\
\homset{\fcat{Colax}}{\cat{L}}{\emm}	\>colax monoidal functors
					\>\ref{sec:Gamma}		\\
\GrAb	\>graded abelian groups	\>\ref{sec:enr-ops}\bref{eg:enr-op:GrLie}\\
$\GrMod=\GrMod_R$\>graded modules	
	\>\ref{sec:prelims-mon-cats}\bref{eg:mon-cat:GrMod},
	\ref{sec:enr-cats}\bref{eg:enr-cat:GrMod}			\\
\homset{\HtyAlg}{P}{\emm}	\>homotopy algebras
			\>\ref{defn:hty-alg}, \ref{defn:enr-hty-alg}	\\
\fcat{HtyMonCat}	\>homotopy monoidal categories
				\>p.~\pageref{page:HtyMonCat}		\\
\twid{\fcat{HtyMonCat}}	\>modified \fcat{HtyMonCat}
				\>p.~\pageref{page:HtyMonCat-twiddle}	\\
$\Mod=\Mod_R$\>modules	\>\ref{sec:prelims-mon-cats}\bref{eg:mon-cat:Mod},
	\ref{sec:enr-cats}\bref{eg:enr-cat:Mod}				\\
\homset{\Mon}{\cat{L}}{\emm}	\>monoidal functors
				\>\ref{defn:mon-ftr}, \ref{sec:free-mon}	\\
\fcat{MonCat}	\>monoidal categories	\>p.~\pageref{page:MonCat}	\\
\homset{\fcat{SColax}}{\cat{L}}{\emm}	\>colax symmetric monoidal functors
					\>\ref{sec:Gamma}		\\
\Set	\>sets	\>\ref{sec:prelims-mon-cats}\bref{eg:mon-cat:Set}	\\
\homset{\SMon}{\cat{L}}{\emm}	\>symmetric monoidal functors
				\>\ref{defn:mon-ftr}, \ref{sec:free-mon}\\
\homset{\fcat{Special}}{\Gamma^\op}{\emm}	
		\>special $\Gamma$-objects
					\>\ref{sec:Gamma}		\\
\homset{\fcat{Special}}{(\Delta^+)^\op}{\emm}	
		\>special simplicial objects
					\>\ref{sec:Gamma}		\\
\Top	\>topological spaces	
		\>\ref{sec:prelims-mon-cats}\bref{eg:mon-cat:Top},
		\ref{sec:enr-cats}\bref{eg:enr-cat:Top}			\\
\Topstar\>based spaces	\>\ref{sec:prelims-mon-cats}\bref{eg:mon-cat:Topstar},
		\ref{sec:A-spaces}					\\
\fcat{Toph}	\>spaces and homotopy classes of maps	
				\>\ref{sec:principle}\bref{eg:change:Toph}\\
\homset{\fcat{WkAlg}}{P}{\emm}	\>weak algebras
				\>\ref{sec:mon-2-cat}			\\
$\Gamma$\>opposite of finite based sets	\>\ref{sec:Gamma}		\\
$\Delta$\>finite totally ordered sets 
			\>\ref{sec:prelims-mon-cats}\bref{eg:mon-cat:Delta}\\
$\Delta_{\mr{inj}}$	\>injective part of $\Delta$
				\>\ref{sec:free-mon}\bref{eg:free-mon:Pt}\\
$\Delta_{\mr{surj}}$	\>surjective part of $\Delta$
				\>\ref{sec:free-mon}\bref{eg:free-mon:Sem}\\
$\Delta^+$	\>non-empty finite totally ordered sets
				\>p.~\pageref{page:Delta-plus}		\\
$\Phi$	\>finite sets		\>\ref{sec:prelims-mon-cats}\bref{eg:mon-cat:Phi}\\
$\Phi_{\mr{inj}}$	\>injective part of $\Phi$
				\>\ref{sec:free-mon}\bref{eg:free-mon:CSem}\\
$\Phi_{\mr{surj}}$	\>surjective part of $\Phi$
				\>\ref{sec:free-mon}\bref{eg:free-mon:CSem}\\
\end{tabbing}

\subsection*{Operads}

\begin{tabbing}
$\GrMod=\GrMod_R$	\=spaces and homotopy classes of maps	\=	\kill
\Act{G}	\>$G$-objects	\>\ref{sec:operads}\bref{eg:operad:Act-G},
	\ref{sec:enr-ops}\bref{eg:enr-op:Act-G},\bref{eg:enr-op:Act-G-alg}\\
\CMon	\>commutative monoids	\>\ref{sec:operads}\bref{eg:operad:CMon}	\\
\CSem	\>commutative semigroups	\>\ref{sec:operads}\bref{eg:operad:CSem}	\\
\Ger	\>Gerstenhaber algebras	\>\ref{sec:enr-ops}\bref{eg:enr-op:Ger}	\\
\GrLie	\>graded Lie algebras	\>\ref{sec:enr-ops}\bref{eg:enr-op:GrLie},
				\ref{subsec:GrLie}			\\
\Inv	\>monoids with involution\>\ref{sec:operads}\bref{eg:operad:Inv}	\\
\Lie	\>Lie algebras	\>\ref{sec:enr-ops}\bref{eg:enr-op:Lie}		\\
\Map{P}		\>$P$-algebra maps	\>\ref{sec:multicats}	\\
\Mon	\>monoids	\>\ref{sec:operads}\bref{eg:operad:Mon}		\\
\Obj	\>objects	\>\ref{sec:operads}\bref{eg:operad:Obj}		\\
\Pt	\>pointed objects	\>\ref{sec:operads}\bref{eg:operad:Pt}	\\
\SAct{G}\>$G$-objects	\>\ref{sec:operads}\bref{eg:operad:Act-G},
	\ref{sec:enr-ops}\bref{eg:enr-op:Act-G}				\\
\Sem	\>semigroups	\>\ref{sec:operads}\bref{eg:operad:Sem}		\\
\SObj	\>objects	\>\ref{sec:operads}\bref{eg:operad:Obj}		\\
\SPt	\>pointed objects	\>\ref{sec:operads}\bref{eg:operad:Pt}	\\
\Sym	\>monoids	\>\ref{sec:operads}\bref{eg:operad:Sym}		\\
\end{tabbing}

\subsection*{Other}

\begin{tabbing}
$\GrMod=\GrMod_R$	\=spaces and homotopy classes of maps	\=	\kill
\eqv		\>equivalence	\>p.~\pageref{page:misc-notation}	\\
\iso		\>isomorphism	\>p.~\pageref{page:misc-notation}	\\
$*$		\>horizontal composite
				\>p.~\pageref{page:misc-notation}	\\
\wej		\>wedge product		
			\>\ref{sec:prelims-mon-cats}\bref{eg:mon-cat:Topstar}\\
\smsh		\>smash product
			\>\ref{thm:n-fold-monoid}			\\
\homset{\cee}{A}{B}	\>hom-set\>p.~\pageref{page:misc-notation}	\\
$S^n$		\>$n$-sphere		\>\ref{sec:iterated}		\\
$S_n$		\>$n$th symmetric group	\>\ref{defn:operad}		\\
$\Delta^n$	\>$n$-simplex		\>\ref{sec:loops}		\\
$\Omega$	\>loop space functor	\>\ref{sec:loops}		\\
\end{tabbing}
\end{document}